\documentclass{amsart}
\usepackage{amssymb}

\usepackage{ifthen}

\input diagrams

\newcommand{\trt}{2}
\newcommand{\lbl}[1]
{\ifthenelse{\trt=3}{\quad\label{#1}\ref{#1}}{\label{#1}}}

\newfont{\rmm}{cmr10 scaled 1000}
\newfont{\bff}{cmbx10  scaled 1000}
\newfont{\itt}{cmsl10 scaled 1000}
\newfont{\nazad}{cmff10 scaled 1000}

\begin{document}

%ARTICLES IN FOREIGN LANGUAGES

\newcommand{\Arno}{ V. Arnold,
\emph{Dynamics of intersections},
        Proceedings of a Conference
 in Honour of
J.Moser, edited by
 P.Rabinowitz and R.Hill,
 Academic Press,
 1990
pp. 77--84.  }

\newcommand{\farrell}
{T.Farrell,
\emph{The obstruction to
fibering a manifold over a circle},
Indiana Univ.~J.~
{\bf 21}
(1971),
    315--346.
}

\newcommand{\fried}{  D.Fried,
\emph{Homological Identities for closed orbits},
Inv. Math.,
{\bf 71},
(1983)
419--442.
}

\newcommand{\laudsiko}
{   F.Laudenbach, J.-C.Sikorav,
\emph{
Persistance d'intersection avec la section
nulle au cours d'une isotopic hamiltonienne
 dans un fibre cotangent
},
Invent.~Math.~
{\bf 82}
    (1985     ),
pp. 349--357.
 }

\newcommand{\milnWT}
{J.Milnor,
\emph{ Whitehead Torsion},
Bull. Amer. Math. Soc.
{\bf 72}
(1966),
358 - 426.
}

\newcommand{\noviquasi}
{S.P.Novikov,
\emph{Quasiperiodic Structures in topology},
in the  book:  Topological Methods in Modern Mathematics,
 Publish or Perish,
 1993,
pp. 223 -- 235.  }

\newcommand{\patou}
{ A.V.Pazhitnov,
\emph{
On the Novikov
complex for rational Morse forms},
 Annales de la Facult\'e de Sciences de Toulouse,
{\bf 4}      (1995), 297 -- 338.
 }

\newcommand{\pasur}
{ A.V.Pajitnov,
\emph{
Surgery on the Novikov Complex},
K-theory
{\bf 10} (1996),      323-412.
 }

\newcommand{\pamrl}
{  A.V.Pajitnov,
\emph{Rationality and exponential growth
properties of the boundary operators in the Novikov
Complex},
Mathematical Research Letters,
{\bf 3}
(1996),
  541-548.
 }

\newcommand{\sikoens}
{ J.-Cl.~Sikorav,
\emph{Un probleme de
disjonction par isotopie symplectique
dans un fibr\'e cotangent},
 Ann.~Scient.
~Ecole~Norm.~Sup.,
{\bf 19}
 1986,
 543--552.
 }

\newcommand{\smale}
{  S.~Smale,
\emph{On the structure of manifolds},
 Am.~J.~Math.,
{\bf 84} (1962)
 387--399.}

\newcommand{\smdyn}
{  S.~Smale,
\emph{Differential dynamical systems},
 Bull. Amer. Math. Soc. {\bf 73} (1967)
747--817.}

\newcommand{\smapoi}
{  S.~Smale,
\emph{Generalized Poincare's conjecture in dimensions
greater than four},
 Ann.~Math.,
{\bf 74} (1961)
 391--406.}

\newcommand{\thom}
{ R.Thom,
\emph{Sur une partition en cellules associ\'ee
\`a une fonction sur une vari\'et\'e},
     Comptes Rendus de l'Acad\'emie de Sciences,
{\bf 228}
(1949),
 973 -- 975.
}

\newcommand{\witt}
{ E.Witten,
 \emph{Supersymmetry and Morse theory},
 Journal of Diff.~Geom.,
{\bf 17} (1985)
 no. 2.
}

%BOOKS IN FOREIGN LANGUAGES

\newcommand{\abrob}
{   R. Abraham, J.Robbin, 
\emph{Transversal mappings and flows
},
 Benjamin, New York,  1967.
}

\newcommand{\dold}
{   A.Dold 
\emph{Lectures on Algebraic Topology
},
Springer,  1972.
}

\newcommand{\kling}
{   W.~Klingenberg
\emph{Lectures on closed geodesics
},
 Springer, 1978.
}

\newcommand{\milnhcob}
{   J.~Milnor,
\emph{Lectures on the
h-cobordism theorem},
 Princeton University
Press,
 1965.
}

\newcommand{\milnstash}
{   J.Milnor and J.Stasheff,
\emph{ Characteristic Classes},
 Princeton University Press,
 1974.}

\newcommand{\morse}
{  M.Morse,
\emph{Calculus of Variations in the Large},
  American Mathematical Society Colloquium Publications,
Vol.18,
 1934.}

\newcommand{\cohn}
{ P.M.Cohn,
\emph{Free rings and their relations},
  Academic press
( 1971)}

\newcommand{\lam}
{    T.Y.Lam,
\emph{Serre's Conjecture,         }
Lecture Notes in Mathematics
{\bf 635},
(1978)
227 p.
}

%PREPRINTS

\newcommand{\huli}{   M.Hutchings, Y.J.Lee
\emph{              Circle-valued Morse theory, Reidemeister torsion
and Seiberg-Witten invariants of 3-manifolds},
e-print dg-ga/9612004 3  Dec 1996.}

\newcommand{\hulihuli}{ M.Hutchings, Y.J.Lee,
\emph{ Circle-valued Morse theory and Reidemeister Torsion
}, e-print dg-ga/9706012 23 June 1997}

\newcommand{\mengt}
{ G.Meng, C.H.Taubes
\emph{SW=Milnor Torsion }
Preprint
(1996)
}

\newcommand{\paodense}
{ A.V.Pazhitnov,
\emph{
On the Novikov
complex for rational Morse forms
},
Preprint of Odense University,
 Odense,
  1991.
 }

\newcommand{\paura}
{  A.V.Pajitnov,
\emph{Surgery on the Novikov Complex},
Rapport de Recherche CNRS URA 758,
 Nantes,
  1993.
}

\def\paepr{  A.V.Pajitnov,
\emph{The incidence coefficients in the Novikov Complex
are generically rational functions},
e-print dg-ga/9603006 14 Mar 96.
 }

\def\paepri{ A.V.Pajitnov,
\emph{Incidence coefficients in the Novikov Complex
for Morse forms: rationality and exponential growth properties},
e-print dg-ga/9604004 20 Apr 96.}

\def\pazeta{ A.V.Pajitnov,
\emph{Simple homotopy type of Novikov complex for closed 1-forms
and Lefschetz $\zeta$-functions of the gradient flow},
e-print dg-ga/970614 9 July 1997.}

\newcommand{\pator}{ A.V.Pajitnov,
\emph{Simple homotopy type of Novikov 
Complex for closed 1-forms and
Lefschetz $\zeta$-functions of the gradient flow},
e-print dg-ga/9706014 26 Jun 97.}

\newcommand{\paepradv}{ A.V.Pajitnov,
\emph{$C^0$-generic properties of boundary operators 
in Novikov Complex
},
e-print math.DG/9812157, to appear in "Advances in Mathematical 
Sciences"
}

\newcommand{\ranprepr}
{  A.Ranicki,
      \emph{
Finite domination and Novikov rings},
preprint,
 1993 }

%ARTICLES IN RUSSIAN

\def\novidok{S.P.Novikov,
\emph{Mnogoznachnye funktsii i funktsionaly.
Analog teorii Morsa}, Doklady AN SSSR
{\bf 260}
(1981),
 31-35.
   English translation:
 S.P.Novikov. \emph{ Many-valued functions
 and functionals. An analogue of Morse theory},
 Sov.Math.Dokl.
{\bf 24}
(1981),
       222-226. }

\newcommand{\farber}
{M.Farber, \emph{Tochnost' neravenstv Novikova},
Funktsional'nyi analiz i ego prilozheniya
 {\bf 19}
 (1985),
     49--59
(in   Russian).
English translation:M.Farber, \emph{
 Exactness of Novikov inequalities },
Functional Analysis and Applications  {\bf
 19},
 1985.
 }

\newcommand{\novshme}
{  S.P.Novikov, I.Schmeltzer,
\emph{ Periodicheskie resheniya uravneniya tipa Kirhgofa dlya dvizheniya
tverdogo tela v zhidkosti i rasshirennaya teoriya
Lusternika-Schirel'mana-Morsa LSchM 1},
Funktsional'nyi analiz i ego prilozheniya
{\bf 15}
no. 3
 (1981).
English translation: S.P.Novikov, I.Shmel'tser, \emph{
 Periodic Solutions of
Kirchhof Equation for the free motion of a rigid body
in a fluid and the extended Theory of Lyusternik - Shnirelman -
Morse (LSM) 1                                        },
 Functional Analysis and Applications
{\bf 15} (1981),
 197 -- 207.
}

\newcommand{\novikirh}
{  S.P.Novikov, \emph{
Variatsionnye metody i periodicheskie resheniya uravneniya tipa
Kirhgofa},
Funktsional'nyi analiz i ego prilozheniya,
{\bf 15},
no. 4
 (1981).
English translation: S.P.Novikov,
\emph{ Variational Methods and Periodic Solutions
of Kirchhof type equations 2},
   Functional Analysis and applications,
{\bf 15}
no. 4
   (1981)
}

\newcommand{\noviuspe}
{ S.P.Novikov,
\emph{ Gamil'tonov formalizm i mnogoznchnyi
analog teorii Morsa}
Uspekhi matematicheskih nauk,
{\bf 37}
   (1982),
 3-49.
English translation: S.P.Novikov,\emph{
The hamiltonian formalism and a
multivalued analogue of
Morse theory                          },
 Russ. Math. Surveys,
{\bf 37} (1982),
 1 -- 56.}

\def\padok{  A.Pajitnov
\emph{O tochnosti neravenstv tipa Novikova
 dlya sluchaya $\pi_1(M)=\ZZZ^m$
i klassov kogomologii v obschem polozhenii}
Doklady ANSSSR
{\bf 306},
 1989
no. 3.
English translation: A.V.
Pazhitnov, \emph{On the sharpness
of the inequalities of Novikov
type in the case
 $\pi_1(M)={\bold Z}^m$
for Morse forms whose cohomology classe are in general
position           },
 Soviet Math. Dokl.
{\bf 39}
 (1989),
no. 3.}

\newcommand{\pasbor}
{  A.Pajitnov
\emph{ O tochnosti neravenstv tipa Novikova dlya sluchaya
$\pi_1(M)=\ZZZ^m$
dlya mnogoobrazii so svobodnoi abelevoi fundamental'noi gruppoi}
(1989)
no. 11
(in   Russian).
English translation: A.V.Pazhitnov,\emph{
 On the sharpness of
 Novikov-type inequalities for
manifolds with free abelian fundamental group.},
 Math. USSR Sbornik,
{\bf 68}
 (1991),
 351 - 389.
}

\newcommand{\pastpet}
{ A.Pajitnov
\emph{
      Ratsional'nost' granichnyh operatorov v komplekse Novikova
v sluchae obschego polozheniya}, 
Algebra i Analiz, {\bf 9}, no.5 (1997), p. 92--139.
English translation:
A.Pajitnov,
\emph{The incidence coefficients in the Novikov complex are 
generically rational
functions}, Sankt Petersburg Mathematical Journal,
{\bf 9} (1998) No.5 
 }

         \newcommand{\tur}
{ V.Turaev
\emph{Kruchenie Raidemaistera v teorii uzlov}
Uspekhi matematicheskih nauk,
{\bf 41:1}
(1986),
  119 - 182.
(in Russian).
English Translation:
V.G.Turaev,
\emph{
Reidemeister Torsion in knot theory,}
Russian Math. Surveys, 41:1 (1986), 119 - 182.}

\newcommand{\turtur}
{V.Turaev,
\emph{ Eilerovy structury, neosobye vektornye polya
i krucheniya tipa Raidemaistera}
Izv ANSSSR , 53:3
(1989),
(in Russian)
English translation:
V.G.Turaev,
\emph{ Euler structures, nonsingular vector fields
and torsions of Reidemeister type},
Math. USSR Izvestia 34:3 (1990), 627 - 662.}

\renewcommand{\a}{\alpha}
\renewcommand{\b}{\beta}
\newcommand{\g}{\gamma}
\renewcommand{\d}{\delta}
\newcommand{\e}{\epsilon}
\newcommand{\ve}{\varepsilon}
\newcommand{\z}{\zeta}
\renewcommand{\t}{\theta}
\renewcommand{\l}{\lambda}
\renewcommand{\k}{\varkappa}
\newcommand{\m}{\mu}
\newcommand{\n}{\nu}
\renewcommand{\r}{\rho}
\newcommand{\vr}{\varrho}
\newcommand{\s}{\sigma}
\newcommand{\vp}{\varphi}
\renewcommand{\o}{\omega}

%GrekBig

\newcommand{\G}{\Gamma}
\newcommand{\D}{\Delta}
\newcommand{\T}{\Theta}
\renewcommand{\L}{\Lambda}
\renewcommand{\P}{\Pi}
\newcommand{\Si}{\Sigma}
\renewcommand{\O}{\Omega}

\renewcommand{\AA}{{\mathcal A}}
\newcommand{\BB}{{\mathcal B}}
\newcommand{\CC}{{\mathcal C}}
\newcommand{\DD}{{\mathcal D}}
\newcommand{\EE}{{\mathcal E}}
\newcommand{\FF}{{\mathcal F}}
\newcommand{\GG}{{\mathcal G}}
\newcommand{\HH}{{\mathcal H}}
\newcommand{\II}{{\mathcal I}}
\newcommand{\JJ}{{\mathcal J}}
\newcommand{\KK}{{\mathcal K}}
\newcommand{\LL}{{\mathcal L}}
\newcommand{\MM}{{\mathcal M}}
\newcommand{\NN}{{\mathcal N}}
\newcommand{\OO}{{\mathcal O}}
\newcommand{\PP}{{\mathcal P}}
\newcommand{\QQ}{{\mathcal Q}}
\newcommand{\RR}{{\mathcal R}}
\renewcommand{\SS}{{\mathcal S}}
\newcommand{\TT}{{\mathcal T}}
\newcommand{\UU}{{\mathcal U}}
\newcommand{\VV}{{\mathcal V}}
\newcommand{\WW}{{\mathcal W}}
\newcommand{\XX}{{\mathcal X}}
\newcommand{\YY}{{\mathcal Y}}
\newcommand{\ZZ}{{\mathcal Z}}

\renewcommand{\aa}{{\mathbb{A}}}
\newcommand{\bb}{{\mathbb{B}}}
\newcommand{\cc}{{\mathbb{C}}}
\newcommand{\dd}{{\mathbb{D}}}
\newcommand{\ee}{{\mathbb{E}}}
\newcommand{\ff}{{\mathbb{F}}}
\renewcommand{\gg}{{\mathbb{G}}}
\newcommand{\hh}{{\mathbb{H}}}
\newcommand{\ii}{{\mathbb{I}}}
\newcommand{\jj}{{\mathbb{J}}}
\newcommand{\kk}{{\mathbb{K}}}
\renewcommand{\ll}{{\mathbb{L}}}
\newcommand{\mm}{{\mathbb{M}}}
\newcommand{\nn}{{\mathbb{N}}}
\newcommand{\oo}{{\mathbb{O}}}
\newcommand{\pp}{{\mathbb{P}}}
\newcommand{\qq}{{\mathbb{Q}}}
\newcommand{\rr}{{\mathbb{R}}}
\renewcommand{\ss}{{\mathbb{S}}}
\newcommand{\ttt}{{\mathbb{T}}}
\newcommand{\uu}{{\mathbb{U}}}
\newcommand{\vv}{{\mathbb{V}}}
\newcommand{\ww}{{\mathbb{W}}}
\newcommand{\xx}{{\mathbb{X}}}
\newcommand{\yy}{{\mathbb{Y}}}
\newcommand{\zz}{{\mathbb{Z}}}

\newcommand{\AAA}{{\mathbf{A}}}
\newcommand{\BBB}{{\mathbf{B}} }
\newcommand{\CCC}{{\mathbf{C}} }
\newcommand{\DDD}{{\mathbf{D}} }
\newcommand{\EEE}{{\mathbf{E}} }
\newcommand{\FFF}{{\mathbf{F}} }
\newcommand{\GGG}{{\mathbf{G}}}
\newcommand{\HHH}{{\mathbf{H}}}
\newcommand{\III}{{\mathbf{I}}}
\newcommand{\JJJ}{{\mathbf{J}}}
\newcommand{\KKK}{{\mathbf{K}}}
\newcommand{\LLL}{{\mathbf{L}}}
\newcommand{\MMM}{{\mathbf{M}}}
\newcommand{\NNN}{{\mathbf{N}}}
\newcommand{\OOO}{{\mathbf{O}}}
\newcommand{\PPP}{{\mathbf{P}}}
\newcommand{\QQQ}{{\mathbf{Q}}}
\newcommand{\RRR}{{\mathbf{R}}}
\newcommand{\SSS}{{\mathbf{S}}}
\newcommand{\TTT}{{\mathbf{T}}}
\newcommand{\UUU}{{\mathbf{U}}}
\newcommand{\VVV}{{\mathbf{V}}}
\newcommand{\WWW}{{\mathbf{W}}}
\newcommand{\XXX}{{\mathbf{X}}}
\newcommand{\YYY}{{\mathbf{Y}}}
\newcommand{\ZZZ}{{\mathbf{Z}}}

\newcommand{\gA}{{\mathfrak{A}}}
\newcommand{\gB}{{\mathfrak{B}}}
\newcommand{\gC}{{\mathfrak{C}}}
\newcommand{\gD}{{\mathfrak{D}}}
\newcommand{\gE}{{\mathfrak{E}}}
\newcommand{\gF}{{\mathfrak{F}}}
\newcommand{\gG}{{\mathfrak{G}}}
\newcommand{\gH}{{\mathfrak{H}}}
\newcommand{\gI}{{\mathfrak{I}}}
\newcommand{\gJ}{{\mathfrak{J}}}
\newcommand{\gK}{{\mathfrak{K}}}
\newcommand{\gL}{{\mathfrak{L}}}
\newcommand{\gM}{{\mathfrak{M}}}
\newcommand{\gN}{{\mathfrak{N}}}
\newcommand{\gO}{{\mathfrak{O}}}
\newcommand{\gP}{{\mathfrak{P}}}
\newcommand{\gQ}{{\mathfrak{Q}}}
\newcommand{\gR}{{\mathfrak{R}}}
\newcommand{\gS}{{\mathfrak{S}}}
\newcommand{\gT}{{\mathfrak{T}}}
\newcommand{\gU}{{\mathfrak{U}}}
\newcommand{\gV}{{\mathfrak{V}}}
\newcommand{\gW}{{\mathfrak{W}}}
\newcommand{\gX}{{\mathfrak{X}}}
\newcommand{\gY}{{\mathfrak{Y}}}
\newcommand{\gZ}{{\mathfrak{Z}}}

\newcommand{\gota}{{\mathfrak{a}}}
\newcommand{\gotb}{{\mathfrak{b}}}
\newcommand{\gotc}{{\mathfrak{c}}}
\newcommand{\gotd}{{\mathfrak{d}}}
\newcommand{\gote}{{\mathfrak{e}}}
\newcommand{\gotf}{{\mathfrak{f}}}
\newcommand{\gotg}{{\mathfrak{g}}}
\newcommand{\goth}{{\mathfrak{h}}}
\newcommand{\goti}{{\mathfrak{i}}}
\newcommand{\gotj}{{\mathfrak{j}}}
\newcommand{\gotk}{{\mathfrak{k}}}
\newcommand{\gotl}{{\mathfrak{l}}}
\newcommand{\gotm}{{\mathfrak{m}}}
\newcommand{\gotn}{{\mathfrak{n}}}
\newcommand{\goto}{{\mathfrak{o}}}
\newcommand{\gotp}{{\mathfrak{p}}}
\newcommand{\gotq}{{\mathfrak{q}}}
\newcommand{\gotr}{{\mathfrak{r}}}
\newcommand{\gots}{{\mathfrak{s}}}
\newcommand{\gott}{{\mathfrak{t}}}
\newcommand{\gotu}{{\mathfrak{u}}}
\newcommand{\gotv}{{\mathfrak{v}}}
\newcommand{\gotw}{{\mathfrak{w}}}
\newcommand{\gotx}{{\mathfrak{x}}}
\newcommand{\goty}{{\mathfrak{y}}}
\newcommand{\gotz}{{\mathfrak{z}}}

\newcommand{\krest}{\begin{picture}(14,14)
\put(00,04){\line(1,0){14}}
\put(00,02){\line(1,0){14}}
\put(06,-4){\line(0,1){14}}
\put(08,-4){\line(0,1){14}}
\end{picture}     }

\newcommand{\tret}{{\frac 13}}
\newcommand{\dvet}{{\frac 23}}
\newcommand{\polt}{{\frac 32}}
\newcommand{\polo}{{\frac 12}}

\renewcommand{\leq}{\leqslant}
\renewcommand{\geq}{\geqslant}

\renewcommand{\th}{therefore}
\newcommand{\ata}{almost~ transversality~ assumption}
\newcommand{\gr}{gradient}
\newcommand{\Mf}{Morse~ function}
\newcommand{\iis}{it is sufficient}
\newcommand{\sut}{~such~that~}
\newcommand{\sufsm}{~sufficiently~ small~}
\newcommand{\wrt}{~with respect to}
\newcommand{\sm}{\setminus}
\newcommand{\ems}{\emptyset}
\newcommand{\sbs}{\subset}
\newcommand{\ho}{homomorphism}
\newcommand{\ma}{manifold}
\newcommand{\nei}{neighborhood}

\newcommand{\ifff}{if and only if}

\newcommand{\FR}{{\mathcal{F}}r}
\newcommand{\gt}{{\mathcal{G}}t}

\newcommand{\aand}{\quad\text{and}\quad}
\newcommand{\wwhere}{\quad\text{where}\quad}
\newcommand{\ffor}{\quad\text{for}\quad}
\newcommand{\for}{~\text{for}~}
\newcommand{\iif}{\quad\text{if}\quad}
\newcommand{\iiif}{~\text{if}~}

\newcommand{\en}{enumerate}

\newcommand{\Prf}{{\it Proof.\quad}}

\newcommand{\Wkr}{W^{\circ}}

\newcommand{\Ker}{\text{\rm Ker }}
\newcommand{\ind}{\text{\rm ind}}
\newcommand{\rk}{\text{\rm rk }}
\renewcommand{\Im}{\text{\rm Im }}
\newcommand{\supp}{\text{\rm supp }}
\newcommand{\Int}{\text{\rm Int }}
\newcommand{\grad}{\text{\rm grad }}
\newcommand{\Id}{\text{\rm Id}}

\newcommand{\nr}{\Vert}
\newcommand{\smo}{C^{\infty}}

\newcommand{\fpr}[2]{{#1}^{-1}({#2})}
\newcommand{\sdvg}[3]{\widehat{#1}_{[{#2},{#3}]}}
\newcommand{\disc}[3]{B^{({#1})}_{#2}({#3})}
\newcommand{\Disc}[3]{D^{({#1})}_{#2}({#3})}
\newcommand{\desc}[3]{B_{#1}(\leq{#2},{#3})}
\newcommand{\Desc}[3]{D_{#1}(\leq{#2},{#3})}
\newcommand{\komp}[3]{{\bold K}({#1})^{({#2})}({#3})}
\newcommand{\Komp}[3]{\big({\bold K}({#1})\big)^{({#2})}({#3})}

\newcommand{\ran}{\{(A_\lambda , B_\lambda)\}_{\lambda\in\Lambda}}
\newcommand{\rran}{\{(A_\lambda^{(s)},
 B_\lambda^{(s)}  )\}_{\lambda\in\Lambda, 0\leq s\leq n }}
\newcommand{\brs}{\rran}
\newcommand{\rans}{\{(A_\sigma , B_\sigma)\}_{\sigma\in\Sigma}}

\newcommand{\fmin}{F^{-1}}
\newcommand{\vh}{\widehat{(-v)}}

\newcommand{\chart}{\Phi_p:U_p\to B^n(0,r_p)}
\newcommand{\atlas}{\{\Phi_p:U_p\to B^n(0,r_p)\}_{p\in S(f)}}
\newcommand{\flow}{{\VV}=(f,v, \UU)}

\newcommand{\Rn}{\bold R^n}
\newcommand{\Rk}{\bold R^k}

\newcommand{\fcob}{f:W\to[a,b]}

\newcommand{\phicob}{\phi:W\to[a,b]}

\newcommand{\vphi}{\varphi}

\newcommand{\crr}{p\in S(f)}
\newcommand{\nrv}{\Vert v \Vert}
\newcommand{\nrw}{\Vert w \Vert}
\newcommand{\nru}{\Vert u \Vert}

\newcommand{\obb}{\cup_{p\in S(f)} U_p}
\newcommand{\proob}{\Phi_p^{-1}(B^n(0,}

\newcommand{\stind}[3]{{#1}^{\displaystyle\rightsquigarrow}_
{[{#2},{#3}]}}

\newcommand{\indl}[1]{{\scriptstyle{\text{\rm ind}\leqslant {#1}~}}}
\newcommand{\inde}[1]{{\scriptstyle{\text{\rm ind}      =   {#1}~}}}
\newcommand{\indg}[1]{{\scriptstyle{\text{\rm ind}\geqslant {#1}~}}}

\newcommand{\obbi}{\cup_{p\in S_i(f)}}
\newcommand{\vem}{\text{Vectt}(M)}
\newcommand{\pr}{\partial}

\newcommand{\id}{\text{id}}

\newcommand{\lau}[1]{{\xleftarrow{#1}}}
\newcommand{\rau}[1]{{\xrightarrow{#1}}}
%newcommand{\lad}[1]{{\xleftarrow[#1]{}}}
\newcommand{\rad}[1]{ {\xrightarrow[#1]{}} }

\newcommand{\xit}{\tilde\xi_t}

\newcommand{\st}[1]{\overset{\rightsquigarrow}{#1}}
\newcommand{\bst}[1]{\overset{\displaystyle\rightsquigarrow}
\to{\boldkey{#1}}}

\newcommand{\stexp}[1]{{#1}^{\rightsquigarrow}}
\newcommand{\bstexp}[1]{{#1}^{\displaystyle\rightsquigarrow}}

\newcommand{\bstind}[3]{{\boldkey{#1}}^{\displaystyle\rightsquigarrow}_
{[{#2},{#3}]}}
\newcommand{\bminstind}[3]{\stind{({\boldkey{-}\boldkey{#1}})}{#2}{#3}}

\newcommand{\Tb}{\text{ \rm Tb}}

\newcommand{\VODIN}{V_{1/3}}
\newcommand{\VDVA}{V_{2/3}}
\newcommand{\VM}{V_{1/2}}
\newcommand{\ddd}{\cup_{p\in S_i(F_1)} D_p(u)}
\newcommand{\dddmin}{\cup_{p\in S_i(F_1)} D_p(-u)}
\newcommand{\where}{\quad\text{\rm where}\quad}

\newcommand{\kr}[1]{{#1}^{\circ}}

\newcommand{\vew}{\text{\rm Vect}^1 (W,\bot)}

\newcommand{\Imm}{\text{\rm Im}}
\newcommand{\hrrr}{\text{\rm Vect}^1(M)}
\newcommand{\vemm}{\text{\rm Vect}^1_0(M)}
\newcommand{\ver}{\text{\rm Vect}^1(\RRR^ n)}
\newcommand{\verr}{\text{\rm Vect}^1_0(\RRR^ n)}

\newcommand{\mods}{\vert s(t)\vert}
\newcommand{\exd}{e^{2(D+\alpha)t}}
\newcommand{\exmin}{e^{-2(D+\alpha)t}}

\newcommand{\intt}{[-\theta,\theta]}

\newcommand{\ffmin}{f^{-1}}
\newcommand{\vvol}{\text{\rm vol}}
\newcommand{\Mat}{\text\rm Mat}
\newcommand{\Tub}{\text{\rm Tub}}
\newcommand{\vxi}{v\langle\vec\xi\rangle}
\newcommand{\tr}{~trajectory }

\newcommand{\grs}{~gradients}
\newcommand{\trs}{~trajectories}

\newcommand{ \co}{~cobordism}
\newcommand{
\sma}{~submanifold}
\newcommand{
\hos}{~homomorphisms}
\newcommand{
\Th}{~Therefore}

\newcommand{
\tthen}{\text \rm ~then}

\newcommand{
\wwe}{\text \rm ~we  }
\newcommand{
\hhave}{\text \rm ~have}
\newcommand{
\eevery}{\text \rm ~every}

\newcommand{
\noconf}{~there~is~no~possibility~of~confusion~}

\newcommand{
\ST}
{\stexp}

\newcommand{\qt}{\hfill\triangle}
\newcommand{\qs}{\hfill\square}

\newcommand{\Vect}{\text{\rm Vect}}

\newcommand{\pa}{\vskip0.1in}

\newcommand{\wi}{\widetilde}

\newcommand{\ove}{\overline}
\newcommand{\unde}{\underline}
\newcommand{\ptf}{\pitchfork}

\renewcommand{\(}{\big(}
\renewcommand{\)}{\big)}

\newcommand{\grd}{{\text{\rm grd}}}

\newcommand{\RA}{\Rightarrow}
\newcommand{\LA}{\Leftarrow}

\newcommand{\emp}{\emptyset}
\newcommand{\wh}{\widehat}

\newcommand{\GC}{\GG\CC}
\newcommand{\GCT}{\GG\CC\TT}
\newcommand{\GT}{\GG\TT}
\newcommand{\GA}{\GG\AA}
\newcommand{\GRP}{\GG\RR\PP}

\newcommand{\GgC}{\GG\gC}
\newcommand{\GgCT}{\GG\gC\TT}

\newcommand{\GgCY}{\GG\gC\YY}
\newcommand{\GgCYT}{\GG\gC\YY\TT}

\newcommand{\GCCT}{\GG\gC\CC\TT}
\newcommand{\GCC}{\GG\gC\CC}

\newcommand{\stv}{\stexp {(-v)}}
\newcommand{\stu}{\stexp {(-u)}}
\newcommand{\stw}{\stexp {(-w)}}

\newcommand{\strv}[2]{\stind {(-v)}{#1}{#2}}
\newcommand{\strw}[2]{\stind {(-w)}{#1}{#2}}
\newcommand{\stru}[2]{\stind {(-u)}{#1}{#2}}

\newcommand{\stvv}[2]{\stind {v}{#1}{#2}}
\newcommand{\stuu}[2]{\stind {u}{#1}{#2}}
\newcommand{\stww}[2]{\stind {w}{#1}{#2}}

\newcommand{\ATA}{Almost~ Transversality~ Condition}
\newcommand{\mx}{\mbox}

\newcommand{\vbsm}{V_b^{\{\leq s-1\}}    }
\newcommand{\vasm}{V_a^{\{\leq s-1\}}    }
\newcommand{\vbs}{V_b^{\{\leq s\}}    }
\newcommand{\vas}{V_a^{\{\leq s\}}    }

\newcommand{\Vbsm}{V_b^{[\leq s-1]}(\d)    }
\newcommand{\Vasm}{V_a^{[\leq s-1]}(\d)    }
\newcommand{\Vbs}{V_b^{[\leq s]}(\d)    }
\newcommand{\Vas}{V_a^{[\leq s]}(\d)    }

\newcommand{\vass}{V_{a_{s+1}}}

\newcommand{\vovo}{\stexp {(-v0)}}
\newcommand{\vov}{\stexp {(-v1)}}

\newcommand{\vbkm}{V_b^{\{\leq k-1\}}    }
\newcommand{\vakm}{V_a^{\{\leq k-1\}}    }
\newcommand{\vbk}{V_b^{\{\leq k\}}    }
\newcommand{\vak}{V_a^{\{\leq k\}}    }

\newcommand{\Vbkm}{V_b^{[\leq k-1]}(\d)    }
\newcommand{\Vakm}{V_a^{[\leq k-1]}(\d)    }
\newcommand{\Vbk}{V_b^{[\leq k]}(\d)    }
\newcommand{\Vak}{V_a^{[\leq s]}(\d)    }

\newcommand{\cob}{~cobordism}

\newcommand{\lc}{\lceil}
\newcommand{\rc}{\rceil}

\newcommand{\sps}{\supset}

\newcommand{\bere}{\begin{rema}}
\newcommand{\bede}{\begin{defi}}

\renewcommand{\beth}{\begin{theo}}
\newcommand{\bele}{\begin{lemm}}
\newcommand{\bepr}{\begin{prop}}
\newcommand{\bega}{\begin{gather}}
\newcommand{\been}{\begin{enumerate}}
\newcommand{\beco}{\begin{coro}}

\newcommand{\bethh}{\begin{ttheo}}
\newcommand{\enthh}{\end{ttheo}}

\newcommand{\beal}{\begin{aligned}}

\newcommand{\enre}{\end{rema}}
\newcommand{\enpr}{\end{prop}}
\newcommand{\enth}{\end{theo}}
\newcommand{\enle}{\end{lemm}}
\newcommand{\enen}{\end{enumerate}}
\newcommand{\enga}{\end{gather}}
\newcommand{\enal}{\end{aligned}}

\newcommand{\enco}{\end{coro}}

\newcommand{\subs}{\subsection}
\newcommand{\hot}{homotopy}
\newcommand{\ity}{\infty}

\newcommand{\TA}{Transversality Condition}

\newcommand{\hog}{homology}

\newcommand{\lb}{\label}

\newcommand{\Vv}{{\boldsymbol{v}}}

\newcommand{\tens}[1]{\underset{#1}{\otimes}}

\newcommand{\babs}{~by the abuse of notation}

\newcommand{\fg}{~finitely generated}

\newcommand{\Hom}{{\text{\rm Hom}}}

\newcommand{\dfm}{diffeomorphism}

\newcommand{\dow}{\pr_0 W}
\newcommand{\daw}{\pr_1 W}

\newcommand{\tidow}{\pr_0 \wi W}
\newcommand{\tidaw}{\pr_1 \wi W}

\newcommand{\su}{Subsection}
\newcommand{\sus}{Subsections}

\newcommand{\she}{simple~homotopy~equivalence}

\newcommand{\sht}{simple~homotopy~type}

\newcommand{\coh}{cohomology}

\newcommand{\Lxi}{\wh\L_\xi}

\newcommand{\Lx}{\L_{(\xi)}}

\newcommand{\iso}{isomorphism}

\newcommand{\bu}{\bullet}

\newcommand{\nics}{\nic s}
\newcommand{\nic}{Novikov incidence coefficient}

\newcommand{\Tr}{\text{Tr}}

\newcommand{\Cl}{\text{Cl}}
\newcommand{\Fix}{\text{Fix}}
\newcommand{\GKSC}{\GG\KK\SS\CC}

\newcommand{\da}{\downarrow}

\newcommand{\scc}[1]{|{\scriptscriptstyle{#1}}}\newcommand{\rrr}{\{\wi r\}}

\newcommand{\celkup}{\pr_1 W^{[k]}}
\newcommand{\celkmup}{\pr_1 W^{[k-1]}}
\newcommand{\celkupfact}{\celkup/ \celkmup}
\newcommand{\celkuppair}{(\celkup/, \celkmup)}

\newcommand{\celkdo}{\pr_0 W^{[k]}}
\newcommand{\celkmdo}{\pr_0 W^{[k-1]}}
\newcommand{\celkdofact}{\celkdo/ \celkmdo}
\newcommand{\celkdopair}{\(\celkdo, \celkmdo\)}

\newcommand{\afcobv}{~Morse function on a cobordism $W$, and $v$ be an
$f$-gradient}

\newcommand{\afcv}{~Morse function on a cobordism $W$, and $v$ is an
$f$-gradient}

%PREDUPREZHD: po sravneniyu s defspe iz advances, 
%izmeneno: \wmok eto teper' s \prec, \succ
%a est' escho \Wmok eto staryi \wmok

%phi-proobrazy itd
\newcommand{\mi}[3]{{#1}^{-1}\([{#2},{#3}]\)}

\newcommand{\fii}[2]{\mi {\phi}{a_{#1}}{a_{#2}} }

\newcommand{\fifi}[2]{\mi {\phi}{#1}{#2} }

\newcommand{\pf}[2]{\mi {\phi_1}{\a_{#1}}{\a_{#2}} }

\newcommand{\mf}[2]{\mi {\phi_0}{\b_{#1}}{\b_{#2}}}

%KOBORDISM  W
\newcommand{\wa}[2]{ W_{[a_{#1}, a_{#2}]}}

\newcommand{\waa}[1]{ W_{[a, a_{#1}]}}

\newcommand{\Wa}[2]{ W_{[{#1}, {#2}]}}

\newcommand{\WS}[1]{ W^{\{\leq {#1}\}}}

\newcommand{\ws}{\WS {s}}

\newcommand{\wsm}{\WS {s-1}}

\newcommand{\wsmm}{\WS {s-2}}

\newcommand{\wk}{\WS {k}}

\newcommand{\wkm}{\WS {k-1}}

\newcommand{\wkmm}{\WS {k-2}}

\newcommand{\tiws}{\wi W^{\{\leq s\}}}

\newcommand{\tiwk}{\wi W^{\{\leq k\}}}

\newcommand{\tiwsm}{\wi W^{\{\leq s-1\}}}

\newcommand{\tiwkm}{\wi W^{\{\leq k-1\}}}

\newcommand{\wsn}{ W^{[\leq s]}(\nu)}

\newcommand{\wsmn}{ W^{[\leq s-1]}(\nu)}

\newcommand{\wsk}{ W^{[\leq k]}(\nu)}

\newcommand{\wkwk}{ W^{[ k]}}

\newcommand{\wmok}{W^{\prec k\succ}}

\newcommand{\tiwmok}{\wi W^{\prec k\succ}}

\newcommand{\womok}{W_0^{\prec k\succ}}

\newcommand{\wmokp}{W^{\prec k+1\succ}}

\newcommand{\wmokm}{W^{\prec k-1\succ}}
\newcommand{\wmokmm}{W^{\prec k-2\succ}}

\newcommand{\wmos}{W^{\prec s\succ}}

\newcommand{\wmosm}{W^{\prec s-1\succ}}

\newcommand{\wmoo}{W^{\prec 0\succ}}

\newcommand{\Wmok}{W^{\langle k\rangle}}

\newcommand{\Womok}{W_0^{\langle k\rangle}}

\newcommand{\Wmokp}{W^{\langle k+1\rangle}}

\newcommand{\Wmokm}{W^{\langle k-1\rangle}}
\newcommand{\Wmokmm}{W^{\langle k-2\rangle}}

\newcommand{\Wmos}{W^{\langle s\rangle}}

\newcommand{\Wmosm}{W^{\langle s-1\rangle}}

\newcommand{\Wmoo}{W^{\langle 0\rangle}}

\newcommand{\hWmok}{\wh W^{\langle k\rangle}}
\newcommand{\hWmokm}{\wh W^{\langle k-1\rangle}}

\newcommand{\tiWmok}{\wi W^{\langle k\rangle}}
\newcommand{\tiWmokm}{\wi W^{\langle k-1\rangle}}

\newcommand{\wwk}{\( \wmok , \wmokm \)}

\newcommand{\wwkp}{\( \wmokp , \wmok \)}

\newcommand{\wwkm}{\( \wmokm , \wmokmm \)}

\newcommand{\wws}{\bigg( \wmos , \wmosm \bigg)}

\newcommand{\wasn}{W^{\lc s\rc}(\nu)}

\newcommand{\wakn}{W^{\lc k\rc}(\nu)}

%verhnii i nizhnii kraya W
\newcommand{\pws}{(\pr_1 W)^{\{\leq s\}}}

\newcommand{\mws}{(\pr_0 W)^{\{\leq s\}}}

\newcommand{\pwk}{(\pr_1 W)^{\{\leq k\}}}

\newcommand{\mwk}{(\pr_0 W)^{\{\leq k\}}}

\newcommand{\pwkm}{(\pr_1 W)^{\{\leq k-1\}}}

\newcommand{\tipwk}{(\pr_1 \wi W)^{\{\leq k\}}}

\newcommand{\timwk}{(\pr_0 \wi W)^{\{\leq k\}}}

\newcommand{\tipwkm}{(\pr_1 \wi W)^{\{\leq k-1\}}}

\newcommand{\pwkmm}{(\pr_1 W)^{\{\leq k-2\}}}

\newcommand{\mwkmm}{(\pr_0 W)^{\{\leq k-2\}}}

\newcommand{\mwkm}{(\pr_0 W)^{\{\leq k-1\}}}

\newcommand{\mwsm}{(\pr_0 W)^{\{\leq s-1\}}}

\newcommand{\mwkp}{(\pr_0 W)^{\{\leq k+1\}}}

\newcommand{\dqr}{\pr_- Q_r}

\newcommand{\ds}{\pr_s}

\newcommand{\dsm}{\pr_{s-1}}

\newcommand{\hwm}{H_*\( \wmok, \wmokm\)}

\newcommand{\hkwm}{H_k\( \wmok, \wmokm\)}

\newcommand{\hwmp}{H_*\( \wmokp, \wmok\)}

\newcommand{\hkwmm}{H_k\( \wmokm, \wmokmm\)}

\newcommand{\hkwmp}{H_k\( \wmokp, \wmok\)}

\newcommand{\yz}{Y_k(v)\cup Z_k(v)}

\newcommand{\Gama}{{\nazad{ \Gamma}}}

\newcommand{\Fi}[1]{\phi^{-1}(a_{#1})}

%v minus itd
\newcommand{\vk}{V_{\langle k\rangle}^-}

\newcommand{\vm}{\wh V^-}
\newcommand{\tivm}{\wi V^-}

\newcommand{\vkm}{V_{\langle k-1\rangle}^-}
\newcommand{\vkp}{V_{\langle k+1\rangle}^-}

\newcommand{\hvk}{\wh V_{\langle k\rangle}^-}

\newcommand{\hvkm}{\wh V_{\langle k-1\rangle}^-}
\newcommand{\hvkp}{\wh V_{\langle k+1\rangle}^-}

\newcommand{\tivk}{\wi V_{\langle k\rangle}^-}

\newcommand{\tivkm}{\wi V_{\langle k-1\rangle}^-}
\newcommand{\tivkp}{\wi V_{\langle k+1\rangle}^-}

\newcommand{\vkvk}{V_{\prec k\succ}^-}

\newcommand{\vkvkm}{V_{\prec k-1\\succ}^-}

\newcommand{\tivkvk}{\wi V_{\prec k\succ}^-}

\newcommand{\tivkvkm}{\wi V_{\prec k-1 \succ}^-}

\newcommand{\cmm}{C_*^\D(\wi M)}

\newcommand{\cvm}{C_*^\D(\wi V^-)}

\newcommand{\ey}{\wi E_*}
\newcommand{\eey}{\wi \EE_*}

\newcommand{\eky}{\wi E(k)_*}

\newcommand{\eti}{\wi{\wi\EE_*}}

\newcommand{\etik}{\wi{\wi\EE}_k}

\newcommand{\etikp}{\wi{\wi\EE}_{k+1}}

\newcommand{\ctiu}{\wi{\wi C}_*(u_1)}

\newcommand{\ctiuk}{\wi{\wi C}_k(u_1)}

\newcommand{\ctiv}{\wi{\wi C}_*(v)}

\newcommand{\ctiukm}{\wi{\wi C}_{k-1}(u_1)}

%Numeration of theorems: section.theorem
%Numeration of formulas: throughout the text as stand Latex
%This is an env-list for showkeys: after "Theorem" etc 
%it leaves a blank line

%V remark i definition avtomaticheski vstavlyaetsya treugolnik

\newcounter{e}[section]

\renewcommand{\thee}{\thesection.\arabic{e} }

\newenvironment{prop}{\vskip0.1in     \noindent
\par\refstepcounter{e}%
   {\bf Proposition
                   \arabic{section}.\arabic{e}. }
\quad\it\vskip0.01in
  }{\vskip0.1in}

\newenvironment{theo}{\vskip0.1in     \noindent
\par\refstepcounter{e}%
   {\bf Theorem
                   \arabic{section}.\arabic{e}. }
\quad\it\vskip0.01in
 }{\vskip0.1in}

\newenvironment{ttheo}{\vskip0.1in     \noindent
   {\bf Theorem
                   }
\quad\it\vskip0.01in
 }{\vskip0.1in}

\newenvironment{coro}{\vskip0.1in     \noindent
\par\refstepcounter{e}%
   {\bf Corollary
                   \arabic{section}.\arabic{e}. }
\quad\it\vskip0.01in
  }{\vskip0.1in}

\newenvironment{lemm}{\vskip0.1in     \noindent
\par\refstepcounter{e}%
   {\bf Lemma
                   \arabic{section}.\arabic{e}. }
\quad\it\vskip0.01in
  }{\vskip0.1in}

\newenvironment{defi}{\vskip0.1in     \noindent
\par\refstepcounter{e}%
   {\bf Definition
                   \arabic{section}.\arabic{e}. }
\quad\vskip0.01in
  }{$\qt$\vskip0.1in}

\newenvironment{defis}{\vskip0.1in     \noindent
\par\refstepcounter{e}%
   {\bf Definitions
                   \arabic{section}.\arabic{e}. }
\quad\vskip0.01in
  }{\vskip0.1in}

\newenvironment{rema}{\vskip0.1in     \noindent
\par\refstepcounter{e}%
   {\bf Remark
                   \arabic{section}.\arabic{e}. }
\quad\vskip0.01in
  }{$\qt$\vskip0.1in}

%Numeration of theorems: section.theorem
%Numeration of formulas: throughout the text as stand Latex
%This is an env-list for showkeys: after "Theorem" etc 
%it leaves a blank line

\newenvironment{propcit}[1]{\vskip0.1in     \noindent
\par\refstepcounter{e}%
   {\bf Proposition
                   \arabic{section}.\arabic{e}. (#1)}
\quad\it\vskip0.01in
  }{\vskip0.1in}

\newcommand{\beprcit}{\begin{propcit}}

\newcommand{\enprcit}{\end{propcit}}

\newenvironment{theocit}{\vskip0.1in     \noindent
\par\refstepcounter{e}%
   {\bf Theorem
                   \arabic{section}.\arabic{e}. }
\quad\it\vskip0.01in
 }{\vskip0.1in}

\newenvironment{corocit}{\vskip0.1in     \noindent
\par\refstepcounter{e}%
   {\bf Corollary
                   \arabic{section}.\arabic{e}. }
\quad\it\vskip0.01in
  }{\vskip0.1in}

\newenvironment{lemmcit}{\vskip0.1in     \noindent
\par\refstepcounter{e}%
   {\bf Lemma
                   \arabic{section}.\arabic{e}. }
\quad\it\vskip0.01in
  }{\vskip0.1in}

\newenvironment{deficit}{\vskip0.1in     \noindent
\par\refstepcounter{e}%
   {\bf Definition
                   \arabic{section}.\arabic{e}. }
\quad\vskip0.01in
  }{\vskip0.1in}

\newenvironment{defiscit}{\vskip0.1in     \noindent
\par\refstepcounter{e}%
   {\bf Definitions
                   \arabic{section}.\arabic{e}. }
\quad\vskip0.01in
  }{\vskip0.1in}

\newenvironment{remacit}{\vskip0.1in     \noindent
\par\refstepcounter{e}%
   {\bf Remark
                   \arabic{section}.\arabic{e}. }
\quad\vskip0.01in
  }{\vskip0.1in}

\title[Novikov complex and $\zeta$-functions]
{ Simple homotopy type of the Novikov \\ Complex and
Lefschetz $\zeta$-function \\
of the gradient flow}
\author{A.V.Pajitnov}
\dedicatory{ To  my teacher S.P.Novikov
for his 60th birthday}
\address{Universit\'e de Nantes\\
D\'epartement de Math\'ematiques\\
2, rue de la Houssini\`ere, 44072, Nantes Cedex France}
\email{pajitnov@math.univ-nantes.fr}
\begin{abstract}
Let $f:M\to S^1$ be a Morse map from a closed 
connected  manifold to a circle.
S.P.Novikov constructed an analog of the Morse complex
for $f$. The Novikov complex is 
a chain complex defined over the ring of
Laurent power series with integral coefficients and
finite negative part. As its classical predecessor this complex 
depends on the choice of a gradient-like vector field.
The homotopy type of the Novikov complex
is the same as the homotopy type of the completed
complex of the simplicial chains of the cyclic covering
associated to $f$.

In  the present paper we prove that for every $C^0$-generic
$f$-gradient there is a homotopy equivalence between these
two chain complexes, such that its torsion equals the Lefschetz
zeta-function of the gradient flow. 
For these gradients the Novikov complex is defined
over the ring of rational functions and the Lefschetz zeta-function
is also rational.
The main theorem of the paper contains a more general
 statement concerning the Lefschetz zeta function with 
twisted coefficients and the version of the Novikov complex 
defined, respectively, over the completion of the group ring of $H_1(M)$.

The paper contains also a survey of Morse-Novikov theory
and of  previous results of the author on the 
$C^0$-generic properties of the Novikov complex and the 
Novikov exponential growth conjecture.
\end{abstract}

\maketitle
\newcommand{\spa}{\hspace{0,5cm}}

\newcommand{\sspa}{\hspace{1cm}}
\newcommand{\df}{\dotfill}

\pa\pa
\centerline{\bf Table of Contents}
\pa

\spa \S\ref{s:intro}
Introduction and statement of results
\df\pageref{s:intro}

\spa \S\ref{s:isr}
Morse-Novikov theory
\df\pageref{s:isr}

\spa \S \ref{s:bssr}
Brief survey of some results of \cite{pastpet} 
\df\pageref{s:bssr}

\spa \S \ref{s:ccohgd}
Condition $(\gC)$ and operator of homological gradient descent
\df\pageref{s:ccohgd}

\spa \S \ref{s:mtfc}
Morse-type filtrations of cobordisms
\df\pageref{s:mtfc}

\spa \S \ref{s:aprel}
Algebraic lemmas\df\pageref{s:aprel}

\spa \S \ref{s:ccmtfrat}
Condition $(\gC\CC)$, Morse-type filtrations for
 Morse maps $M\to S^1$, and $C^0$-generic
rationality of Novikov incidence coefficients
\df\pageref{s:ccmtfrat}

\spa \S \ref{s:proof}
Proof of the main theorem
\df\pageref{s:proof}

\pa
 References
\df\pageref{refer}

\newpage

\section{Introduction and statement of results}
\lb{s:intro}

We begin  with Subsection \ref{su:inintro}
which presents to the reader the basic ideas of the paper.
The next three subsections contain the minimal set of definitions
necessary for the statement of the main theorem,
given in
Subsection \ref{su:stat}.
In the text  of this section the reader will find  references
 to other parts of the text which contain further explanations, details
 etc.

\subsection{Introduction}
\label{su:inintro}

In this paper we continiue the study of $C^0$-generic
 properties of Novikov complex, began in
the papers \cite{pamrl}, \cite{pastpet}.
Let $M$ be a closed connected manifold.
Recall that for a Morse map
$f:M\to S^1$
and an $f$-gradient $v$, satisfying
\TA, an analog of the classical Morse complex --
the Novikov complex $\CC_*(v)$ -- can be associated. It is a free
finitely generated chain complex over the ring
$\ZZZ[[t]][t^{-1}]$
of the formal power series with finite negative part.
In the cited papers we constructed for a given Morse map $f:M\to S^1$
a special class of $f$-gradients which form a subset,
 $C^0$-open-and-dense
in the subset of all $f$-gradients satisfying \TA.
For any $f$-gradient $v$ in this class the boundary  operators in Novikov
complex, associated to $v$, are not merely power series, but rational
 functions.
The main construction in the proof is a certain group \ho~
$h(v)$, introduced in \S 4 of \cite{pastpet},
which we call here   {\it homological gradient descent}.
(To explain the origins of this notion, consider a regular
value $\l\in S^1$
of $f$, and consider the integral
curve $\g$ of $(-v)$ starting at
$x\in V=f^{-1}(\l)$.
If this curve does not converge to a critical point of $f$, then
it will intersect again the \sma~$V$ at some point, say, $\rho(x)$.
The map $\rho$ is thus a (not everywhere defined ) smooth map
of $V$ to $V$, and the {\it homological gradient descent operator}
is a substitute for the \ho~induced by $\rho$ in \hog.
This operator is defined for a $C^0$-generic $v$-gradient satisfying \TA.
)
In the present paper we develop further the properties of the
homological gradient descent.
We show that this \ho~carries enough information to compute
the  Lefschetz zeta-function of the gradient flow.
This leads to the main theorem of the present paper.
See Subsection \ref{su:stat} for the precise statement.
To explain intuitively the idea of the theorem,
recall (\cite{patou}) that there is a chain \hot~ equivalence
between the Novikov complex and the
completed
chain simplicial complex of the corresponding cyclic covering.
The main theorem of the present paper says that
for $C^0$-generic $f$-gradient $v$
this
chain \hot~equivalence can be chosen
so that its torsion equals to Lefschetz zeta-finction
of the flow generated by $(-v)$.

\subsection{Algebraic notions: $K_1$, Novikov rings etc}
\label{su:algk1}
Let $R$ be a commutative ring with a unit.
Recall from \cite{milnWT}
the group
$\ove{K}_1 R=K_1 R/\{0, [-1]\}$
where $[-1]$ denotes  the element of order 2 corresponding to the unit
$(-1)\in GL(R,1)\sbs GL(R)$.
Further,
$K_1R=R^\bu\oplus SK_1(R)$
where $R^\bu$
is the full group of units of $R$.
Let $U$ be a subgroup in $R^\bu$.
Set $K_1(R\vert U)=
\ove{K}_1 R/U=R^\bu/(\pm U)\oplus SK_1(R)$
(here $\pm U=\{\pm u\vert u\in U\}$ and we identify the units
in $R^\bu$ with their images in $\ove{K}_1 R$).
For $S\in GL(R)$ the image of $S$ in $\ove{K}_1 R$
or in $K_1(R\vert U)$
will be denoted by $[S]$.

A free based finitely generated chain complex $C_*$ of right $R$-modules
(see\cite{pasur}, def. 1.3)
is called {\it $R$-complex}.

Recall from \cite{milnWT}, \S 3
that for an
acyclic $R$-complex $C_*$
the torsion $\tau(C_*)\in \ove{K}_1 R$
is defined.

If $\phi:F_*\to D_*$
is a \hot~equivalence
of $R$-complexes, then its {\it torsion}
$\tau(\phi)\in \ove{K}_1 R$
is defined to be the torsion
of the chain cone $C_*(\phi)$.
The image of
$\tau(\phi)$
in the group
$K_1(R|U)$ will be denoted by $\tau(\phi|U)$.

Proceeding to the definitions related to Novikov completions, let
$G$ be an abelian group.
Let $\xi:G\to \RRR$
be a \ho.
Set $\L=\ZZZ G$
and denote by
$\wh{\wh\L}$
the abelian group of all functions $G\to\ZZZ$.
Equivalently,
$\wh{\wh\L}$
is the set of all formal linear combinations 
$\sum_{g\in G} n_g g$
(not necessarily finite)
of the elements of $G$ with integral coefficients.
For $\l\in \L$
set
$\supp\l=
\{g\in G\mid n_g\not= 0\}$.
Denote by $\Lxi$ the subgroup of
$\wh{\wh \L}$
formed by all the $\l$ \sut~
for every $C>0$ the set
$\xi^{-1}([C,\infty[)\cap\supp\l$
is finite. Then $\Lxi$ has a natural structure of  a ring.
We shall also need an analog of this ring with $\ZZZ$  replaced
by $\QQQ$.
That is, let
 $\wh{\wh \L}_\QQQ$
 be the set of all the functions $G\to\QQQ$
and set
$$
\wh\L_{\xi,\QQQ}=
\{\l\in
\wh{\wh \L}_\QQQ
\mid \forall C\in\RRR \quad \supp\l\cap\xi^{-1}([C,\infty[)
\mbox{ is finite }\}
$$

In the sequel we shall work with two subgroups of units of $\Lxi$.
 First is the group
$G$ itself, the second is
introduced in the following formula:
$$
U_\xi=\{\l\in\Lxi\mid\l=\pm g(1+\mu), \mbox{ where } g\in G \mbox{ and }
\supp\mu\sbs \xi^{-1}\(]-\infty, 0[\)\}
$$
In \cite{patou} the group $K_1(\Lxi|U_\xi)$ was denoted by $Wh(G,\xi)$.

In the present paper we shall  also use a ring which is smaller than
$\Lxi$ but still has a lot of essential properties of $\Lxi$.
Set $S_\xi=\{\l\in\ZZZ G\mid \l=1+\mu,
\supp\mu\sbs \xi^{-1}\(]-\infty, 0]\)$.
Denote $S_\xi^{-1}\L$ also by $\Lx$, note that we
have  natural inclusions
$\L\rInto\Lx\rInto\Lxi$.

\subsection{Novikov complex}
\label{su:ncintro}

Here we just say {\it what it is}, 
leaving the details of definition and 
the construction together with motivations for the section
\ref{su:simnt}.

Let $f:M\to S^1$
be a Morse map.
We assume that $f_*:H_1(M)\to H_1(S^1)=\ZZZ$
is epimorphic.
 Let $\PP:\wi M\to M$ be a connected regular
covering with structure group $G$, \sut~$f\circ \PP$ is \hot~to zero.
\footnote{
Thus we use the symbol $\tilde{}$, normally reserved for the universal
 covering,
for another purpose,
but there will be no  possibility for confusion.}
In this paper we consider only the case when $G$ is abelian.
Then we obtain the natural epimorphism
$\pi:H_1(M)\rOnto G$
and $f_*:H_1(M)\to\ZZZ$ factors through a \ho~$\xi: G\to\ZZZ$.
Set $\L=\ZZZ G$, then the corresponding Novikov ring $\Lxi$ is defined.
Let $v$ be an $f$-gradient satisfying \TA~ (see Subsection \ref{su:mftgbt}
for definition).
Choose for every critical point $p$ of $f$ an orientation of the
stable \ma~of $p$
and a lifting of $p$ to $\wi M$. To this data one associates a
$\Lxi$-complex $\wi \CC_*(v)$ \sut~the number of free 
generators of $\wi \CC_k(v)$
equals to the number of critical points of $f$ of index $k$.

\subsection{Terminology: flows and orbits}
\label{su:termorb}

Let $w$ be a $\smo$ vector field on a closed manifold $M$.
{\it A closed orbit} of $w$ is a non constant 
trajectory $\g$ of $w$, defined on a closed interval
$[a,b]\sbs\RRR$
and satisfying 
$\g(b)=\g(a)$.
We shall identify two such trajectories
$\g:[a,b]\to\RRR,
\g':[a',b']\to\RRR
$
if there is $C$, such that
$a'=a+C, b'=b+C$
and
$\g(t)=\g'(t+C)$.

The set of all closed orbits is denoted by 
$Cl(w)$.
For 
$\g\in Cl(w),
\g:[a,b]\to M$
and $m\in \NNN$
define 
$\g^m\in Cl(w)$
as the result of successive gluing of $m$ copies of
$\g$ to each other.
{\it The multiplicity}
$m(\g)\in \NNN$
is the maximal number $k\in\NNN$
such that 
$\g=\theta^m$
for some 
$\theta\in Cl(w)$.
For every 
$\g\in Cl(w)$ its homology class 
$[\g]\in H_1(M)$
is well defined.
A vector field $w$ will be called {\it Kupka-Smale} if
every zero and every closed orbit is hyperbolic, and
$w$
satisfies the \TA.
For every Kupka-Smale vector field $w$ and every 
$\g\in Cl(w)$
the index 
$\ve(\g)\in \{-1,1\}$
is defined as the index of the corresponding Poincare map.

\subsection{Lefschetz zeta-functions}
\label{su:lfzeta}

We shall consider zeta-functions only for a very particular class
of flows: flows on closed maniflods generated by
gradients of Morse maps
$f:M\to S^1$.

The set of all Kupka-Smale gradients of $f$ is denoted by
$\GG\KK\SS(f)$.
Let
$v\in
\GG\KK\SS(f)$.
Consider $\pi([\g])$ as an element of the group ring
$\L=\QQQ G$, multiply it by the rational number
$       \frac {\ve(\g)}{m(\g)}      $
and consider an infinite series
$$
\eta_L(-v)=
\sum_{\g\in Cl(-v)}
\frac {\ve(\g)}{m(\g)}
\pi([\g])
$$
The Kupka-Smale property implies that there is at most finite set
of $\g\in Cl(-v)$ with given 
$[\g]$. Thus the above series is a well-defined
 element of the abelian group
$\wh{\wh\L}_\QQQ$.
It is not difficult to show that
$\eta_L(-v)$
belongs to the Novikov ring
$\wh{\L}_{\xi, \QQQ}$.
Note also that
$\supp \eta_L(-v)\sbs \xi^{-1}(]-\infty, -1])$.
Therefore
$\supp (\eta_L(-v))^k\sbs \xi^{-1}(]-\infty, -k])$
and the power series
$$
\zeta_L(-v)=\exp(\eta_L(-v))
$$
is well defined and belongs again to
$\wh{\L}_{\xi, \QQQ}$.
This element will be called
{\it Lefschetz zeta-function}.

\subsection{Statement of the main theorem}
\label{su:stat}

Let $M$ be a closed connected \ma,
$f:M\to S^1$ be a Morse map. 
We assume here the terminology of 
Sect. \ref{su:ncintro}.
For a unit $\l$ of the ring $\L_{(\xi)}$
we shall denote the image
of $\l$
in
$\ove{K}_1(\L_{(\xi)}\mid G)$
by $\bar \l$.
We denote by
$\cmm$
the simplicial complex of $\wi M$ associated with a smooth triangulation
of $M$, so that $\cmm$
is a $\ZZZ G$-complex.

\bethh
There is a subset
$\GG\KK\SS\gC (f)
\sbs\GG\KK\SS (f)$
with the following properties:
\been
\item
$\GG\KK\SS\gC (f)$ is open-and-dense subset of
$\GG\KK\SS (f)$ \wrt~ $C^0$-topology.
\item For every $v\in
\GG\KK\SS\gC (f)$
the Novikov complex $\wi \CC_*(v)$ is defined over $\L_{(\xi)}$, and
$\zeta_L(-v)\in\Lx$.
\item
For every $v\in\GG\KK\SS\gC (f)$
there is a \hot~equivalence
$$
\phi:\wi \CC_*(v)\rTo^{\sim}
\cmm\tens{\L}\Lx
$$
\sut
~$\tau(\phi|G)=\ove{\zeta_L(-v)}$.
\enen
\enthh

\subsection{Contents of the paper section by section
and further remarks}
\label{su:contents}

The results of this paper first appeared in the e-print 
\cite{pator}. The present paper contains the detailed 
exposition of the main part of the results of
\cite{pator}; due to the lack of time we left aside the results of 
\cite{pator} concerning the irrational Morse forms, as well as the 
relations to Seiberg-Witten 
invariants of 3-manifolds.

The
subsections \ref{su:mftgbt} -- \ref{su:mcf} of Section \ref{s:isr}
establish the terminology we need.
In the subsections \ref{su:simnt} and \ref{su:negcohgd}
the reader will find  an introduction to Morse-Novikov 
theory
and an exposition of the author's results on the Novikov
Exponential Growth Conjecture (\cite{pamrl}, \cite{pastpet}).
This part is written with the emphasis on the ideas 
rather than on the technical details.

The section \ref{s:bssr}
contains a compressed survey of the techniques and results of 
\cite{pastpet}
used in the sequel.
We included this part in order to make the exposition self-contained 
in a sense
(still I did not include the proofs, for which the reader is invited
 to consult
\cite{pastpet}).
The contents of the section \ref{s:ccohgd} is a version of the Section 4
of \cite{pastpet}; the condition 
$(\gC)$ which we introduce here is of the same origin as the condition 
$(RP)$ of the \S 4 of \cite{pastpet}, but better suited for our present 
purposes.

The section  \S \ref{s:mtfc}
  is one of the central parts
of the paper.
See the beginning of \S \ref{s:mtfc}
for an introduction to the main ideas of this section.

\S \ref{s:aprel} contains some simple algebraic computations.
The section \ref{s:ccmtfrat}
contains the most important step of the proof of the main theorem:
we compute here the simple homotopy type of the Novikov complex in terms
of the Morse filtrations of
cobordisms
(introduced in \S \ref{s:mtfc}) and related invariants.
The section \ref{s:proof}
contains the end of the proof of the main theorem.

We end the present introduction by a remark on the origins of our paper.
Using the methods of \cite{pastpet}
it is not difficult to prove that the Lefschetz zeta-function of a 
$C^0$-generic 
Kupka-Smale gradient is rational, but one does not see immediately 
how to compute it
in terms of usual invariants of the \ma.
The formula of Hutchings-Lee (\cite{huli}, Theorem 1.1)
was the first example of such computation. It suggested that in general 
case there should be a formula relating the simple \hot~ type
of the Novikov complex, the simple \hot~ type of the \ma~ 
itself and the zeta-function of the gradient flow.

I am grateful to S.P.Novikov, Bai-Ling Wang and  V.Turaev 
for valuable discussions.

\newpage

\section{Morse-Novikov theory}
\lb{s:isr}

\subsection{Morse functions  and their gradients: basic terminology}
\lb{su:mftgbt}

Before going into the details of Morse-Novikov theory
and Novikov Exponential Growth Conjecture, we  first
 list some basic notions
in order to establish the Morse-theoretic language
which we use.

Let $v$ be a $C^1$ vector field on a \ma~ $M$.
The value of the integral curve of $v$ passing by $x$ at $t=0$
will be denoted by $\g(x,t;v)$.

\pa
We call {\it cobordism}
a compact manifold $W$
together with a presentation
$\pr W=\pr_0W\sqcup \pr_1W$
where
$\pr_1W$ and $\pr_0W$
are compact manifolds
without
boundary
of dimension $\dim W-1$
(one or both of them can be empty). The \ma
~$W\sm\pr W$
will be denoted by $\Wkr$.

\pa
A {\it Morse function}
$f:W\to[a,b]$ on a cobordism $W$
is a $\smo$ map $f:W\to\RRR$, \sut ~
$f(W)\subset[a,b], f^{-1}(b)=\pr_1W,
f^{-1}(a)=\pr_0W$,
all the critical points of $f$ are
non-degenerate and
belong to
$\Wkr$.

The set of all critical points of a Morse
 function $f$ will be denoted by $S(f)$; the set
 of all critical points of $f$ of index $k$
will be denoted by $S_k(f)$.

\pa
Let $f$  be a Morse function on a cobordism
$W, \dim W=n$.
A vector field $v$ is called {\it gradient } of $f$,
or $f$-gradient, 
if for every critical point $p$ of $f$ there is a chart $\phi:U\to V$
around $p$, \sut~ $\phi_*(v)$ is a standard vector field
$(-x_1,..., -x_k, x_{k+1},..., x_n)$ on $\RRR^n$ (where $k=\ind p,
n=\dim M$), and $f\circ  \phi^{-1}$ is a quadratic form
$\sum_i\a_ix_i^2$ with $\a_i>0$ for $i>k$ and $\a_i<0$ for $i\leq k$.

\pa
The set of all $f$-gradients will be denoted by $\GG(f)$.
\bere
\been\item
The reader will check easily that
if we  demand $\| a_i \|=1$ in the definition above,
we obtain an equivalent definition.
\item
Riemannian gradient of $f$ (\wrt~ some riemannian
metric)
 is not included in this definition. There is a natural way
to extend the definition so as to cover the case of Riemannian
gradients and to carry over  our results to this general framework
(see 
\cite{paepradv}).
\enen\enre

\pa
From now on up to the end of this subsection
$\fcob$ is a Morse function on a cobordism of
dimension $n$, and $v$ is an $f$-gradient.
\pa
Denote by $K_1$ the set of all
$x\in \pr_1 W$, \sut~ the $(-v)$-trajectory
starting at $x$ converges to a critical point of $f$.
Similarly, $K_0$ denotes the set of all $x\in \pr_0W$,
\sut~
$\g(x,t;v)$ converges to a critical point of $f$.
Note that $K_0, K_1$ are compacts.
The shift along the trajectories of $(-v)$ defines
a diffeomorphism $\pr_1 W\sm K_1\to \pr_0W\sm K_0$, which will
be denoted by $\stv$.
\label{stind}
The set
$\stv (A\sm K_1)$ will be denoted (\babs) by
$\stv(A)$.

Let $\l<\m$ be regular values of $f$. Set $W'=f^{-1}\([\l,\m]\)$, and
$w=v\mid W'$. The diffeomorphism $\stw$ will be denoted by
 $\stind {(-v)}\mu\l$.

We shall need a construction which is quite standard in Morse theory
(see \cite{milnhcob}, p.62).
Let $\Psi:\dow\times[0,l]\to W$ be
a map, given by
$\Psi(x,t)=\g(x,t;v)$, and choose $l$ small
enough so that $\Psi$ is a diffeomorphism onto its image.
Let $h$ be a $\smo$ positive function on $[0,l]$,
$\supp h\sbs]0,l[$. Set $T=\int_0^lh(\tau)d\tau$.
Let $u$ be a vector field on $\dow$.
Define a vector field $w$ on $\dow\times [0,l]$
by $w(x,t)=h(t)u(x)$.
Define a vector field $w'$
on $W$ setting 
$w'=\Psi_*(w)$
in $\Im\Psi$
and $w'=0$
in $W\sm\Im\Psi$.
Set $v'=v+w'$.
If $h$ and $u$ are sufficiently small then
$\Vert v-v'\Vert$
is small, and $v'$ is still an $f$-gradient. Note that
$\stexp {(-v')}=\Phi(-u, T)\circ \stv$.
We shall call this construction 
{\it
adding to $v$ a horizontal
component $u$ nearby $\pr_0 W$.
}

\pa
 We say that
 $v$ satisfies
           {\it Transversality Condition }, if
$$\big(
x,y \in S(f) \big)
\Rightarrow \big( D(x,v)\cap \Wkr \pitchfork D(y,-v)\cap
\Wkr \big)
$$
\pa

 We say that
 $v$ satisfies {\it Almost Transversality Condition},
if
   $$\big( x,y \in S(f) ~\&~
\ind x \leq \ind y \big) \Rightarrow \big(
D(x,v)\cap\Wkr \pitchfork D(y,-v)\cap\Wkr\big)$$
\vskip0.1in

\pa
The set of all $f$-gradients satisfying
\TA, resp. \ATA~
will be denoted by $\GT(f)$, resp. by $\GA(f)$.

\pa

 A Morse function $\phi:W\to[\a,\beta]$
is
called {\it adjusted to $(f,v)$}, if:

 1) $S(\phi)=S(f)$, and $v$ is also a $\phi$-gradient.

 2) The function $f-\phi$ is constant in a \nei~
of $\pr_0W$, in a \nei~ of $\pr_1 W$, and in a \nei~ of each
point of
$S(f)$.

\pa

We say that $f$ is {\it ordered} Morse function
with an {\it ordering sequence} $(a_0,...,a_{n+1})$, if
$a=a_0<a_1<...<a_{n+1}=b$ are regular values of $f$
\sut~ $S_i(f)\subset
 f^{-1}(]a_i,a_{i+1}[)$.

\subsection{More terminological conventions}
\lb{su:mtc}

We adopt the convention that structure
 groups of regular coverings act on the cover {\it from
the right}.
All the regular coverings
considered in this paper will have abelian structure groups
(a large part of our results can be carried over to
the non-abelian case without problems).

\subsection{Morse complex for functions: recollections}
\lb{su:mcf}

Before proceeding to Morse-Novikov theory, we shall recall briefly
the classical notion
of Morse complex of a Morse function.
Let $g:M\to\RRR$ be a Morse function on a closed manifold $M$.
The classical Morse theory \cite{morse}
implies that
$m_p(g)\geq b_p(M)$, where $m_p(g)$
is the number of critical points of index $p$, and $b_p(M)$ is the
$p$-th Betti number of $M$.
This can be refined as to construct a chain complex
$C_*$ of free finitely generated abelian groups, \sut~the
number $\mu(C_p)$
of free generators of $C_p$
equals to
$m_p(g)$ and $H_*(C_*)\approx H_*(M)$.
The construction of this complex historically
was done  in several steps,
the main progress being done in the papers of R.Thom
\cite{thom} and S.Smale \cite{smapoi}.
See the book \cite{milnhcob}
of J.Milnor for systematic exposition of this step.
(In this book one does not find
the definition of the complex, but
it can be extracted from the \S 6 of this book.)
Later this construction was reconsidered by Witten \cite{witt}
from the entirely new point of view, involving DeRham cohomology.

Here is the main idea of the construction.
Choose a $g$-gradient $v$
satisfying \TA.
Define  $C_s$ to be a free abelian group generated by
the critical points of $g$ of index $s$.
Let $p,q$ be critical points of $g$ of
indices $k$, resp. $k-1$.
Since $v$ satisfies  \TA, the set $L(p,q)$
of all orbits, joining $p$ to $q$
is finite (exercise for the reader).
For each critical point of $g$  choose an orientation
of the stable \ma~ of this point
(\wrt~the vector field $v$). Then to each orbit
$\g\in L(p,q)$
a sign $\ve(\g)$ can be attributed in a natural way.
Set
$n(p,q)=\sum_{\g\in L(p,q)}\ve(\g)$.
The  boundary
$\pr_k: C_k\to C_{k-1}$
is defined by
$\pr_kp=\sum_q n(p,q)q$
(where $q$ ranges over the set of critical points of $g$ of index $k-1$).

It turns out  that $\pr_k\circ\pr_{k+1}=0$ for every $k$
and that the \hog~ of $C_*$
\wrt~ $\pr_*$
is isomorphic to $H_*(M)$
(the details can be found in \cite{patou}, Appendix). This complex is
 called {\it Morse complex}.

\subsection{Short introduction to Morse-Novikov theory}
\lb{su:simnt}

In the paper \cite{novidok} Novikov have laid the foundation of
Morse theory for closed 1-forms. One of the basic concepts of
 this theory
is that of {\it Novikov complex}.
We shall give here a brief outline
of the construction
(see \cite{patou} for more details).

Consider a Morse map
$f:M\to S^1$.
The construction from \ref{su:mcf}  can not be carried over to this
 case as it is, since for two
critical points $x,y$
of adjacent indices the number of the trajectories of a gradient
$v$ joining $x$ to $y$ can well be infinite.
Here is a  way round this difficulty:
split the set of all trajectories joining $x$ and $y$
into some number of disjoint finite subsets indexed by some family $I$,
count separately the trajectories in each class and obtain then
an incidence coefficient as a function
$I\to\ZZZ$.
Thus the base ring of the resulting complex will be much larger
than $\ZZZ$.
We  shall describe the most natural way of such splitting
(introduced in \cite{novidok})
which leads to
the base ring $\ZZZ((t))$).

\bede
The ring $\ZZZ((t))$ consists of
 power series in $t$
with integral coefficients and finite
negative part.
That is,
$\l=\sum_{-\infty}^{\infty} a_it^i$ is in $\ZZZ((t))$,
if $a_i\in\ZZZ$ for all $i$ and there is
$N=N(\l)$, \sut~
$a_i=0$ if $i<N(\l)$.

\end{defi}

We shall assume that our Morse map $f$ is not homotopic to zero,
otherwise $v$ is a gradient
of an $\RRR$-valued Morse function,
and the  construction  from Subsection \ref{su:mcf} applies.
Let $\CC_p(f)$ be the free $\ZZZ((t))$-module generated by
critical points of $f$ of index $p$.
Consider the infinite cyclic cover
$\CC:\bar M\to M$ \sut~ $f\circ \CC$ is homotopic to zero.
Lift the function $f:M\to S^1$ to a function $F:\bar M\to\RRR$.
Let $t$ be the generator of
the structure group of
 $\CC$ \sut~ $F(xt)< F(x)$.
Choose a gradient $v\in \GT(f)$ and lift it to
$\bar M$ (we shall keep for this lifting the same notation
$v$).
As above choose for every critical point $x\in S(f)$
an orientation of the stable manifold
of $x$ \wrt~$v$.
Moreover,
for each critical point $x$ of $f:M\to S^1$
choose a lifting $\bar x$ of $x$ to $\bar M$.
Let $x,y$ be a pair of critical points of $f$ with $\ind x=\ind y+1$.
It is not difficult to show that
for any $k\in\ZZZ$ the set of all $v$-orbits joining $\bar x$ and
$\bar yt^k$
is finite.
Sum them up with the corresponding signs (as above)
and obtain an integer $n_k(x,y)$.
Set $n(x,y)=\sum_{k} n_k(x,y)t^k\in\ZZZ((t))$
and set
$\pr_k x=\sum_{y}y\cdot n(x,y)$.
(The power series $n(x,y)$ are called
{\it novikov incidence coefficients}.
Note that $n(x,y)$ and $n_k(x,y)$ depend on $v$ and
sometimes we shall write $n(x,y;v)$, resp.
$n_k(x,y;v)$ to stress this dependance.)
Note the Novikov complex is for us a complex of {\it right modules}
(see Subsection \ref{su:mtc}).

One can show that that $\pr_k\circ\pr_{k+1}=0$.
Thus we obtain a chain complex of free finitely generated
 $\ZZZ((t))$-modules.
We shall denote it by $\CC_*(v)$ in order to stress its dependance on
$v$, which is essential for our purposes.\footnote{We have chosen the
italic $\CC$ in order to distinguish between
the Novikov complex and the ordinary Morse complex.}
(This complex actually depends also on the choice of
orientations of stable manifolds of critical points
of $f$, and on the liftings of the critical points to $\bar M$,
but the influence of these choices on the result is less important:
the ambiguity is reduced to the multiplication
of the lines and columns of the matrices of the boundary operators
$\pr_k$
by some units of of the form $\pm t^m$ of the base ring.
Note, on the other hand that this complex does not depend on
the particular choice
of a function $f$ for which $v$ is a gradient.
)
To describe the \hog~  of the complex, recall that $H_*(\bar M)$ is a
$\ZZZ[t,t^{-1}]$-module, thus the tensor product
$\wh H_*(\bar M)=
H_*(\bar M)\tens{\ZZZ[t,t^{-1}]}\ZZZ((t))$
makes sense. One can prove that
$H_*(\CC_*(v))\approx \wh H_*(\bar M)$.

All the cited results on the Novikov complex are stated in \cite{novidok}.
The detailed proofs
can be found in \cite{patou}.

\subsection{Novikov exponential growth conjecture
 and the operator of homological gradient descent}
\label{su:negcohgd}

Let $f:M\to S^1$ be a Morse map, and $v$ be a generic $f$-gradient,
satisfying \TA.
We have seen in the previous subsection that to this data the
 Novikov complex $\CC_*(v)$ is associated.
It is a free based chain complex over $\ZZZ((t))$,
therefore each boundary operator
$\pr_k$ is represented by a matrix $D^{(k)}$
with the entries
$D_{ij}^{(k)}$
are in $\ZZZ((t))$ (these entries are nothing else than the
 Novikov incidence coefficients).
Each of $D_{ij}^{(k)}$
is then a power series
$a(t)=\sum_k a_kt^k$
with a finite negative part.
The Novikov exponential growth conjecture says that the
coefficients $a_k$ grow
at most exponentially when $k\to\infty$
(maybe with some restrictions on $v$ of the type analyticity
or general position imposed).
See the Introduction of \cite{pastpet}
for more details on this conjecture.
This conjecture was proved in \cite{pastpet} for any Morse
map $f:M\to S^1$
and a set $\GT_0(f)$ of $f$-gradients, which is $C^0$
open and dense in $\GT(f)$.
For the gradients in $\GT_0(f)$ the Novikov incidence coefficients
are actually Taylor series of some rational functions.
The proof of the rationality of the incidence coefficients
for $v\in \GT_0(f)$
is based on the construction
of what we call
{\it operator of homological gradient descent}.
This construction is essential for our work so we shall
spend the rest of this subsection explaining the main idea of the
construction.
See the detailed exposition in Subsections \ref{su:rshgd}, \ref{su:hgdsv}.

Assume that the \hot~ class of $f$ in $[M,S^1]\approx H^1(M,\ZZZ)$
is indivisible.
Let $p,q$ be critical points of $f$,
$\ind p=\ind q+1$; we shall work with the Novikov incidence coefficient
$n(p,q;v)$.
Assume for simplicity of notation that $1\in S^1$
is a regular value of $f$, that the first critical level after $f(q)$ in
the counter-clockwise direction is $f(p)$
and that $1$ is between $f(q)$ and $f(p)$.
Cut $M$ along $f^{-1}(1)$  to obtain a cobordism $W$ with two parts of
the boundary:
$\pr_1 W$ and $\pr_0 W$.
There is the identification  diffeomorphism $\Phi:\pr_0 W\to \pr_1 W$.
We obtain also a Morse function $f_0:W\to[0,1]$
and its gradient, which will be denoted by the same letter $v$.
Note that $p$ is the lowest critical point of $f_0$, and $q$ is the
 highest critical point of $f_0$.
The descending disc of $p$ intersects $\pr_0 W$ by
an embedded sphere $S(p)$. The ascending
disc of $q$ intersects $\pr_1 W$ by an embedded sphere $S(q)$.
The Novikov incidence coefficient
$n_k=n_k(p,q;v)$
is then the algebraic intersection number of
$(\stv\circ\Phi)^k (S(p))$
with $\Phi^{-1}(S(q))$.
\footnote{See page \pageref{stind} for the definition of
$\stexp{(-v)}$.}
Denote
$\stv\circ\Phi $
by $\phi$.
If $\phi$ were an everywhere defined diffeomorphism, then,
denoting the \hog~ classes
of $S(p)$ and of $\Phi^{-1}(S(q))$
by, respectively, $[p]$ and $[q]$,
we would obtain the following formula
\begin{equation}
n_k(p,q;v)=(\phi_*)^k([p])\krest [q]\qquad
\lb{f:ifif}
\end{equation}
(where $\krest$ stands for the algebraic intersection index
of \hog~ classes; we assume for the moment
that $\pr_0 W$ is oriented, although actually the orientability is
not necessary).
Then the rationality of the power series $\sum_kn_kt^k$
would be a simple fact of linear algebra.
( For the sake of completeness we shall indicate here this
linear-algebraic argument.\lb{linalg}
Let $L$ be a free \fg~ abelian group, $A:L\to L$
be a \ho.
We shall identify $A$ with its matrix.
 Then the power series
$\sum_{k=0}^{\infty} A^k t^k$ represents the matrix
$(1-At)^{-1}$.
It is easy to deduce from this that

\begin{equation}\lb{f:lalg}
\begin{gathered}
\mbox{ for each } x\in L \mbox{ and each } \xi\in \Hom (L, \ZZZ)
\mbox{ the power series }\\
\sum_{k=0}^{\infty} \xi (A^k x) t^k
\mbox{ is a rational function with denominator } 1-t\det A.
\end{gathered}
\end{equation}

By the way,
the similar argument proves that the homological Lefschetz
$\zeta$-function of a diffeomorphism
is always rational see \cite{smdyn}, p.767 -- 768).
Unfortunately $\phi$ is not everywhere defined, the  soles of
the ascending
discs form exactly the set of indetermination of  $\phi$. \footnote{
Novikov writes in \cite{noviquasi}:   {\it
In our case we do not have a mapping, but a cobordism with
two equal boundaries generating the $\ZZZ$-covering
over the compact manifold. }}
The
iterations of $\phi$ only accumulate this indetermination.
But it is {\it generically }
possible to construct some objects which will serve
as the substitutes of
$[p]$, $[q]$ and  $\phi_*$
in the formula (\ref{f:ifif}).

Now I must prevent the reader that the exposition in the following
 paragraphs up to the end of the section \ref{s:isr}
is very informal.
One should consider it as a {\it program}
which was realized in \cite{pastpet}.
Nevertheless I think that the understanding
 of these few lines
is rather helpful for understanding of the proofs in
\cite{pastpet} and of the present paper.

Start with a "nice" cellular decomposition of $\pr_1 W$. By
 "nice" I mean in particular that
every cell $e$ is represented by
a $\smo$ map $\rho$ of $D^k$ to $\pr_1 W$
\sut~$\rho\mid \Int D^k$
is an embedding onto a submanifold of $\pr_1 W$ ($k$ is the
dimension of the cell).
One expects that such nice decomposition can be constructed
from the stratification
of $\pr_1 W$ by the stable manifolds of critical points of some
Morse function. Carry over this decomposition to $\pr_0 W$ by
means of the diffeomorphism
$\Phi^{-1}$.
The map
$\stv$
is not a continuous map
of $\pr_1 W$ to $\pr_0 W$
but nevertheless one can try to perturb $v$
 so that the map
$\stv$ corresponding to the perturbed gradient be cellular.
To explain why it should work, let $p$ be a critical point
of $f$ of index $k$. Denote
$D(p,v)\cap\pr_0 W$
by
$S_-(p,v)$
and
$D(p,-v)\cap\pr_1 W$
by
$S_+(p,v)$.
Perturbing $v$ we can assume that all the
manifolds $S_+(p,v)$
are transversal to all the cells
of $\pr_1 W$.
In particular the cells of $\pr_1 W$
of dimension $k$ will intersect
only the manifolds
$S_+(p,v)$
with $\ind p\leq k$.
Let $e$ be a cell of $\daw$ of dimension $k$.
The part of $e$ where $\stv$
is {\it not}
defined is then the union of $S_+(p,v)\cap e$,
where $p$ ranges over critical points of $f$ of index
$\leq k$.
The  geometric picture of behavior
of the descending trajectories of the gradient
is as follows.
The gradient descent shrinks the set $S_+(p,v)\cap e$
down to the point $p$
(during infinite time).
If we try to descend still lower to $\pr_0 W$, then
we shall obtain a map which is multivalued in the sense that to every
point of $S_+(p,v)\cap e$
corresponds the {\it set}
$S_-(p,v)$.
Note however that
$\dim S_-(p,v)\leq k-1$.
Assume for a moment that 
for every $s$
all the sets $S_-(p,v)$ with $\ind p\leq s+1$
were in the $s$-skeleton of
$\pr_0 W^{[s]}$
(which is plausible in view of the cellular approximation procedure).
Then our "map"
$\stv$
will be defined as a continuous map
$$W_k:
C_k^{(1)}=\celkupfact
\to
C_k^{(0)}=\celkdofact$$
or, equivalently, as a map from $C_k^{(0)}$
to itself, since $\Phi$ identifies $\pr_0 W$
with $\pr_1 W$.
The space $C_k^{(1)}$
is of course a wedge of
spheres of dimension $k$. Denote its \hog~ by $H_k$;
it is a free abelian group, and the map $W_k$ induces an endomorphism
$\xi$ of $H_k$.
Let $r,q$ be critical points of $f$ of indices resp. $k+1, k$.
By the argument similar to the above the manifold
$S_-(r,v)$ defines a relative cycle
in
$\celkdopair$,
hence an element $[r]\in H_k\celkdopair$.
Intersection index with $S_-(q,v)$
defines a \ho~
$\l(q):H_k\to\ZZZ$ and we have:
$n_k(r,q;v)=\l(q)\(\xi^k([r])\)$.
By the  argument above from linear algebra  (page \pageref{linalg}), the
Novikov incidence coefficient $n(r,q;v)$ is a rational function.

Now I shall make two remarks concerning the above program.
\begin{enumerate}
\item
The resulting perturbation of our vector field $v$
can be chosen $C^0$-small, but in general {\it not} $C^\infty$-small.
\item
We have consciously simplified the exposition, working with
cellular decompositions of
$\pr_0 W$ and $\pr_1 W$. In fact this does not work this way; 
instead of cellular decompositions we should consider 
{\it handle decompositions} of these manifolds
(the sets $D_\delta (\indl k; v)$ in the notation of
\cite{pastpet}).
\end{enumerate}

\newpage

\section{Brief survey of some results of \cite{pastpet} }
\lb{s:bssr}

We shall use the techniques and the results of \cite{pastpet}
in the following sections; so we collected in this section
some basic ideas and results of \cite{pastpet}.

\subsection{ $\d$-thin handle decompositions}
\lb{su:dthd}

In this subsection $W$ is a riemannian cobordism
of dimension $n$, $\fcob$ is a Morse function on $W$, and $v$ is an
$f$-gradient.

Let $x\in W$. Let $\d>0$.
Assume that for some $\d_0>\d$ the restriction of
the exponential map
$\exp_q:T_qW\to W$ to the disc $B^n(0,\d_0)$
is a \dfm~on its image.
Denote by $B_\d(p)$ (resp. $D_\d(p)$) the riemannian open ball
(resp. closed ball)
of radius $\d$
centered in $p$. We shall use the notation
$B_\d(p), D_\d(p)$ only when the assumption above on $\d$ holds.

Set
\begin{gather*}
B_\delta(p,v)=\{x\in W~\vert~ \exists
t\geq 0 :
 \gamma (x,t;v)\subset B_\delta (p)\}\\
D_\delta(p,v)=\{x\in W~\vert ~\exists
t\geq 0 :
\gamma (x,t;v)\subset D_\delta (p)\},\\
D(p,v)=\{x\in W\vert \lim\limits_{t\to\infty}
\gamma(x,t;v)=p\}
\end{gather*}

We  denote by $D (\indl s;v)$
the union of   $
D(p,v)$ where $p$ ranges over critical
points of $f$ of index $\leq s$.
We denote by $B_\d (\indl s;v)$, resp.
by $D_\d(\indl s;v)$
the union of
$B_\d(p,v)$, resp. of
$D_\d(p,v)$, where $p$ ranges over critical
points of $f$ of index $\leq s$.
We shall use
similar
notation like $D_\d(\inde s;v)$ or
$B_\d(\indg s;v)$, which is now clear
without special definition.
The union of all $D(p,v)$ is denoted by
 $D(v)$, the union of all $D_\d(p,v)$ is denoted by
$D_\d(v)$.
For the uniformity of notation we shall also denote
$D(p,v)$ by $D_0(p,v)$, and $D(v)$ by $D_0(v)$ etc.

Let $\phi:W\to[a,b]$ be an ordered Morse function with an
ordering sequence
$(a_0<a_1... <a_{n+1})$. Let $w$ be a $\phi$-gradient.
Denote $\phi^{-1}\([a_i,a_{i+1}])$ by $W_i$.

\begin{defi}
We say that $v$ is {\it $\d$-separated \wrt~
$\phi$ } (and the ordering sequence $(a_0,...,a_{n+1})$), if

i) for every $i$ and every $p\in S_i(f)$ we have
$D_\d(p)\sbs \Wkr_i$;

ii) for every $i$ and every $p\in S_i(f)$
there is a Morse function
$\psi:W_i\to[a_i,a_{i+1}]$, adjusted to
$(\phi\mid W_i, w)$
and a regular value $\l$ of $\psi$ \sut~
$$D_\d(p)\sbs\psi^{-1}(]a_i,\l[)$$
and for every $ q\in S_i(f), q\not= p$ we have
$$D_\d(q)\sbs \psi^{-1}(]\l,a_{i+1}[)$$ 

\end{defi}

We say that $v$ is {\it $\d$-separated }if
it is $\d$-separated \wrt~some ordered Morse
function $\phi:W\to[a,b]$, adjusted to $(f,v)$.

It is obvious that each $f$-gradient satisfying \ATA~ is $\d$-separated
for
some $\d>0$.

\begin{prop}\lb{p:ddd}
If $v$ is  $\delta_0$-separated, then
$\forall
\delta\in[0,\delta_0]$ and
$\forall
s: 0\leq s\leq n$
\been
\item $D_\delta (\indl s;v)$         is compact.
\item ${\bigcap} _{\theta>\delta} B_\theta (\indl s;v)
= D_\delta (\indl s;v)$
\item  If $\d>0$     then $\overline{B_\delta (\indl s;v)}
 = D_\delta (\indl s;v)$. $\qs$
\enen
\end{prop}

Thus the  collection of descending discs $D(p,v)$ form
a sort of stratified manifold, and
the $\d$-thickened descending discs $B_\d(v)$
 form  a collection of
neighborhoods of this manifold.
In the preceding proposition we have listed some
natural properties of these objects.
These properties deserve to be formulated  and
studied separately. Thus we are lead to
the notions of $s$-submanifold and of
$ts$-submanifolds,  which will be introduced
  in the following subsection.

\subsection{ $s$-submanifolds and $ts$-submanifolds}
\lb{su:ssttss}

Let $\aa =\{A_0,...,A_k\}$
be a finite sequence of subsets of
a topological space $X$.
We  denote  $A_s$ also by $\aa_{(s)}$, and
we denote  $A_0\cup...\cup A_s$
by $\aa_{(\leq s)}$.
We say that $\aa$ is a
{\it compact family }
if  $\aa_{(\leq s)}$
is compact for every $s$.
\begin{defi}
Let $M$ be a manifold without boundary.
 A finite sequence $\xx = \break \{X_0,...,X_k\} $
 of subsets of $M$ is called  {\it
$s$-submanifold of $M$ }
($s$ for stratified) if
\begin{enumerate}
\item  Each
 $X_i$ is a submanifold of $M$ of dimension $i$
with the trivial normal bundle.
\item $\xx$ is a compact family.
\item For $i\not= j$ we have: $X_i\cap X_j=\emp$
\end{enumerate}
\end{defi}
 For a diffeomorphism $\Phi :M\to N$ and
 an $s$-submanifold
 $\xx$ of $M$ we denote by $\Phi (\xx)$
the $s$-submanifold
of $N$ defined by $\Phi (\xx)_{(i)} =
\Phi(\xx_{(i)})$.
\vskip0.1in

Let $\xx , \yy$ be  two $s$-submanifolds
of $M$. We say, that $\xx$
 is {\it transversal to $\yy$ }
 (notation: $\xx \pitchfork \yy$) if
$\xx_{(i)}\pitchfork \yy_{(j)}$ for every
$i,j$; we say that $\xx$
is {\it almost transversal to }
   $\yy$ (notation: $\xx\nmid\yy$) if
$\xx_{(i)}\pitchfork \yy_{(j)}$ for
$i+j<\dim M$. Note, that
$\xx\nmid\yy$ if and only if $ \xx _{(\leq i)}
\cap  \yy _{(\leq j)} =
 \emptyset$
whenever $i+j<\dim M$.

There is an analog of Thom transversality theorem:
given two $s$-submanifolds
$\xx, \yy$ of $M$, there is a small isotopy $\Phi_t$
of $M$, \sut~
$\Phi_1(\xx)\nmid \yy$
(see \cite{pastpet}, Theorem 2.3.

\begin{defi}
Let $X$ be a topological space,
$\aa = \{A_0,...,A_k\}$ be a  compact  family  of  subsets  of
$X$~,~   $I$  be an  open  interval  $]0,\delta_0[$.  A  {\it  good
fundamental system of neighborhoods of $\aa$ }  (abbreviation:
$gfn$-system  for  $\aa$)  is  a  family  $\AAA=\{  A_s(\delta
)\}_{\delta\in  I,0\leq  s\leq  k}$  of  open  subsets  of  $X$,
satisfying the following conditions:
\been\item
     For every $s$ and
 every $\delta_1<\d_2$ we have $A_s\subset A_s(\d_1)\sbs A_s(\d_2)$

 \item  For every $\delta$  and every $j$
we have
$\overline{A_{\leq j}(\delta)} =
{\bigcap}_{\theta >\delta}\big( A_{\leq j}(\theta)\big)$

     \item  For every   $j$
 we have
$ A_{(\leq j)} =
{\bigcap}_{\theta >0}\big( \AAA_{(\leq i)}(\theta)\big)$.
\enen
$I$ is called {\it interval of definition
}
of the system. $\aa$ is called the {\it core } of $\AAA$, and
$\AAA$ is called {\it thickening } of $\aa$.
We shall denote $A_s(\d)$ also by $\AAA_{(s)}(\d)$ and
$A_{\leq i}(\d)$ also by $\AAA_{(\leq i)}(\d)$.
Sometimes we shall denote $A_s$ by $\AAA_s(0)$ for uniformity
of notation. $\qt$
\end{defi}

The term "good fundamental system" is justified by the
 following easy consequence of the properties 1 -- 3 :
if $X$ is compact, then
for every $\d\geq 0$ and every $s$ the family
$\{A_{\leq s}(\theta)\}_{\theta>\d}$
is a fundamental system of neighborhoods of
$\ove{\AAA_{(\leq s)}(\d)}$.

   If $M$ is a manifold without
boundary, $\xx$ is an $s$-submanifold of $M$,
and $\XXX$ is a $gfn$-system for $\xx$, then
we say that $\XXX$ is a
{\it $ts$-submanifold  with the core $\xx$. }
We say also that $\XXX$ is the {\it thickening} of $\xx$.
For a $ts$-submanifold
$\XXX =\{ X_s(\delta )\}_{\delta\in I,0\leq s\leq k}$
we shall denote $X_i(\delta)$ by $\XXX _{(i)}(\delta)$,
and $X_{\leq i}(\delta)$
  by
$\XXX _{(\leq k)}(\delta)$.

Here is the basic example of $gfn$-system. Let $\fcob$
be a \afcobv. Assume that $v$ is $\e$-separated.
Then
 the family
$\{ B_\delta (\inde s;v)\}_{\delta\in ]0,\epsilon[,0\leq s\leq n}$
is a $gfn$-system for $\dd (v)$.
If $W$ is a closed manifold, this family
is a $ts$-submanifold with the core $\dd(v)$.
 This $gfn$-system will be denoted
by $\DDD (v)$.
(If $\pr W\not= \emp$
then $\dd(v)$ is a family of manifolds with boundary, but
 we shall neither use it,
nor formulate the corresponding generalization of the
 notion of $ts$-submanifold.)

Let $\fcob$ be a \afcobv. For a subset $X\sbs \pr_1 W$
the set $T(X,v)$ formed by all the $(-v)$-trajectories
starting in $X$ will be called {\it track} of $X$.
The set $T(X,v)$ is not necessarily compact even if $X$ is.
But one can show that if
$X$ is compact, then $T(X,v)\cup D_\d(v)$ is compact for every $\d$
(see \cite{pastpet}, Lemma 2.7).
Next we shall extend the notion of track to $s$-submanifolds of $\daw$.
Let
$\aa$ be an $s$-submanifold of $\pr_1W$.
Intuitively, the
$k$-th component of the track-family
$\ttt(\aa,v)$
is the track $T(A_k,v)$ plus all the discs 
$D(p,v)$ with $\ind p\leq k$.
Proceeding  to precise definitions,
 we shall assume
from here to the end of the section
that $v$ satisfies the \ATA, and that
$\aa\nmid\dd _b(-v)$.
(The condition $\aa\nmid\dd_b(-v)$ is not very restrictive,
since given $\aa$ and $v$
one can always find a $\smo$-small
perturbation of $v$ \sut~this condition holds,
see \cite{pastpet}, Lemma 4.2.)

 \begin{defi}
Set  $TA_i(v) = T(A_{i-1} ,v)\cup D(\inde i ;v)$
($A_{-1}=\emptyset$ by definition).
Denote by $\ttt (\aa ,v)$ the
 family
$\{ TA_i(v) \}_{0\leq i\leq k+1}$
 of subsets of $W$
and by $\stind {(-v)}b\lambda (\aa)$ the family
$\{ TA_{i+1}
(v)
\cap f^{-1}(\lambda)\}_{0\leq i\leq k}$
of subsets of $f^{-1}(\lambda)$.
If the values $b,\lambda$ are clear from the context
 we shall abbreviate
$\stind {(-v)}b\lambda (\aa)$ to $\st {(-v)} (\aa)$.
The family $\ttt (\aa ,v)$ will be called {\it track }
of $\aa$, and the family $\stind {(-v)}b\lambda (\aa)$
will be called
{\it $\st {(-v)}$-image } of $\aa$.
\end{defi}

One can prove that $\ttt(\aa,v)$ and
$\stind {(-v)}b\l(\aa)$
 are compact families;
further, if $\l$ is a regular value of $f$, then
$\stind {(-v)}b\l(\aa)$
is an $s$-submanifold of $f^{-1}(\l)$ (see \cite{pastpet}, Lemma 2.9).

Now we can introduce tracks of $ts$-submanifolds.
If $\AAA$ is a $ts$-submanifold of $\dow$
with the core $\aa$, then the track of $\AAA$ will be
a $gfn$-system with the core $\ttt(\aa,v)$.

\begin{defi}
 Let
 $\AAA =\{ A_s(\delta )\}_{\delta\in ]0,\delta_0[ ,0\leq s\leq k}$
be a $ts$-submanifold of $\daw$ with the core
$\aa$.
Assume that $v$ is $\delta_1$-separated.
 For $0<\delta<\min (\delta_0,\delta_1)$
set
$TA_s(\delta , v) = T(A_{s-1} (\delta),v)
\cup
 B_\delta (\inde s; v)$ ($A_{-1}(\delta)=\emptyset$ by definition).$\qt$
\end{defi}

\beprcit{\cite{pastpet}, Prop. 2.12}
\lb{p:trackgfn}
There is $\epsilon \in ]0,\min (\delta_0,\delta_1)[$
such that
$\{TA_s(\delta , v) \}_{\delta\in]0,\epsilon[,~0\leq s\leq k+1}$
is a gfn-system for $\ttt (\aa ,v)$.
\enprcit

The $gfn$-system, introduced in the Proposition
\ref{p:trackgfn}
will be denoted by $\TTT(\AAA,v)$
and called {\it track of $\AAA$}.

\subsection{ Rapid flows}
\lb{su:rf}

Let $M$ be a closed \ma~ of dimension $m$, $f$ be a Morse
function on $M$, $v$ be an $f$-gradient
satisfying \ATA. Let $k$ be an integer and $U$ be an open \nei~ of
$D(\indl k;-v)$, $V$ be an open \nei~ of $D(\indl {m-k-1};v)$.
It is clear that
the flow $\Phi_t$, generated by
$v$ will carry $M\sm V$ to $U$
if only we wait sufficiently long.
(In precise terms: denote the diffeomorphism
$t\mapsto \g(x,t;v)$
of $M$ by $\Phi(v,t)$, then
for $T$ large enough we have:
$\Phi(v,T)(M\sm V)\sbs U$.)
If the flow $\Phi_t$ does this operation in a small time
and moreover if $v$ has small norm, we shall say that $\Phi_t$ is
{\it rapid}.
Here is the precise definition.

\begin{defi}\label{d:rapid}
Let $v$ be an $f$-gradient, satisfying \ATA.
Let $\e>0, t\geq 0$.
We say that the flow, generated by $v$ is {\it $(t,\e)$-rapid},
if
for every $s:0\leq s\leq m$ we have:
\begin{equation}\label{qpu1}
\Phi(v,t)\(M\sm B_\e(\indl s;v)\)\sbs
B_\e\(\indl {m-1-s}; -v\)
\end{equation}
\begin{equation}\label{qpu2}
\Phi(-v,t)\(M\sm B_\e(\indl s;-v)\)\sbs
B_\e\(\indl {m-1-s}; v\)
\end{equation}
(Sometimes for the sake of brevity we shall say
that $v$ is $(t,\e)$-rapid.) $\qt$
\end{defi}

It is not difficult to check that each flow $v$ satisfying \ATA~
is $(t,\e)$-rapid for {\it some} numbers $t,\e$.

{\bf Remark on the Terminology.}
The reader will not find the notion of
rapid flow in \cite{pastpet}.
In \cite{pastpet}
we used the notion of {\it quick flows}
which is better suited for the proof of the
theorem
\ref{t:quickpush} below.

\begin{defi}
Let $\fcob$ be a Morse function on
a  cobordism, and $u,v$ be $f$-gradients.
We say that $u$ and $v$ are
{\it equivalent},
if $u(x)=\phi(x)v(x)$,
where $\phi:W\to\RRR$ is a $\smo$ function, \sut~$\phi(x)>0$
for every $x$
and $\phi(x)=1$
for $x$ in a \nei~ of $S(f)$.
\end{defi}

\begin{theo}\label{t:quickpush}

Let $M$ be a closed riemannian manifold. Let $C>0, t>0$.

There is a Morse function
$f:M\to\RRR$ and an $f$-gradient $v$, satisfying \ATA~
such that for every $\e>0$ there is
a $(t,\e)$-rapid $f$-gradient $u$
equivalent to $v$, and
satisfying  $\Vert u\Vert\leq C$.
\end{theo}

This theorem is proved in \cite{pastpet} (Corollary 1.14).
\pa
\subsection{ Condition $(RP)$}
\lb{su:cnrp}

In this subsection $W$ is a riemannian cobordism, $f:W\to [a,b]$
a Morse function, $v$ an $f$-gradient, satisfying \ATA,  $n=
 \dim W$.
\begin{defi}\lb{d:crp}
We shall say that $v$ has a ranging pair
(or, equivalently, that $v$ satisfies condition (RP))
if there are Morse functions
$\phi_0:\dow\to\RRR, \phi_1:\daw\to\RRR$,
and their gradients $u_0$, resp. $u_1$,
satisfying \ATA, and a number $\d>0$, \sut~ for every
$s$ the following conditions hold:
\begin{equation}
\d~ \text{\rm is~in the interval of the definition of
}~\TTT(\DDD(u_1),v) \text{\rm~and of~}
\TTT(\DDD(-u_0), -v)\tag{RP1}
\end{equation}
\begin{equation}
\text{The gradients}~~
v,u_0\aand u_1 \quad \text{are}\quad \delta-
\text{separated}
\tag{RP2}
\end{equation}
\begin{equation}
\overline{\TTT (\DDD (u_1), v) _{(\leqslant s+1)}
 (\delta)} \cap \dow
\subset
\DDD (u_0)_
{(\leqslant s)}(\delta)\tag{RP3}
\end{equation}
\begin{equation}
\overline{\TTT    (\DDD (-u_0), -v )
_{(\leqslant s+1)}     (\delta)}
\cap \daw
\subset
\DDD (-u_1)_{(\leqslant s)}(\delta)
 \tag{RP4}
\end{equation}
\end{defi}

The set of $f$-gradients of $v$ satisfying
(RP) will be denoted by $\GRP(f)$.

\beth\lb{t:rpdense}
The set $\GRP(f)$ is $C^0$ dense in $\GG(f)$.
\enth
\Prf
Let $v\in \GA(f)$.
Let $\phi_1:\daw\to\RRR,
\phi_0:\dow\to\RRR$
be Morse functions, $u_1, u_0$ be their gradients, satisfying \ATA.
Here if the plan of our proof.
To obtain an $f$-gradient
$w\in\GRP(f)$
we shall modify $v$ by adding the horizontal components:
$\xi_1$ nearby $\daw$ and $\xi_0$ nearby $\dow$.
The vector field $\xi_1$ will be equivalent to $u_1$
and $\xi_0$ will be equivalent to $u_0$.
We shall first do this modification without making any assumptions on
$u_1, u_0$,
and we shall obtain the condition (RP) for the modified gradient.
In the end we shall show how to choose $u_1, u_0$
so that the resulting $f$-gradient $w$ be $C^0$ close to $v$.

Proceed now to the proof.
Perturbing $v$ if necessary we can assume that
$\stv(\dd(u_1))$
is almost transversal to $\dd(-u_0)$,
that is,
$\stv(\dd_{(\leq k)}(u_1))$
does not intersect
$\dd_{(\leq m)}(-u_0)$
if $k+m<n-1$.
Therefore by the properties of $gfn$-systems
we have for small $\d$:
$$
\stv (D_\d(\indl k;u_1))
\cap
D_\d(\indl m; -u_0)=\emp
\mbox{ if } k+m<n-1
$$

To make our notation more brief, set:
\begin{gather*}
A_k^-(\d)=D_\d(\indl k; u_1); \qquad
A_k^+(\d)=D_\d(\indl k; -u_1); \\
B_k^-(\d)=D_\d(\indl k; u_0); \qquad
B_k^+(\d)=D_\d(\indl k; -u_0); \\
\end{gather*}

(Thus $A_k^-(\d)$ is "the $\d$-thickened $k$-skeleton of $\daw$",
and
$A_k^+(\d)$ is
"the $\d$-thickened $k$-skeleton of the dual cell decomposition of
$\daw$". Similar interpretation is valid for $B_k(\d)^{\pm}$.)

Let $\xi_0$ be a $\phi_0$-gradient, equivalent to
$u_0$.
Let $\xi_1$ be a $\phi_1$-gradient equivalent to $u_1$.
For $T>0$ sufficiently large and any $m$
we have:

\begin{gather}
\Phi(\pm\xi_0, T)(\dow\sm B_k^\mp(\d))\sbs
B_{n-k-2}^\pm(\d)\lbl{f:nul}\\
\Phi(\pm\xi_1, T)(\daw\sm A_k^\mp(\d))\sbs
A_{n-2-k}^\pm(\d)\lbl{f:odin}
\end{gather}

Set $\Phi_0=\Phi(\xi_0, T), \quad \Phi_1=\Phi(\xi_1, T)$.
We shall add to $v$ two horizontal components:
 $\xi_0$ nearby $\dow$, and
$\xi_1$ nearby $\daw$.
If we choose the parameter function $h$ with
$\int_0^lh(\tau)d\tau=T$, then
 for the resulting
vector field $v'$ we have:
\begin{equation}
\stexp {(-v')}=
\Phi_0\circ \stv\circ\Phi_1^{-1}
\lbl{f:compo}
\end{equation}
It is obvious now that
(RP3) and (RP4)
hold for $v'$. Indeed, the set
$A_m^-(\d)$
is sent to
itself by $\Phi_1^{-1}$, then the result is sent to
$\dow\sm B_{n-2-m}^+(\d)$
by
$\stv$ and $\Phi_0$ pushes the
$\dow\sm B_{n-2-m}^+(\d)$
to
$B_m^-(\d)$. This proves (RP3), similar for (RP4).
(Of course we should take $\d$ small enough in order to
satisfy (RP1) and (RP2).
One must also make sure that $v'$ be an $f$-gradient.
We refer to \cite{pastpet}, \S 4 for details.)
Now
I shall indicate how to make the norm  $\Vert v'-v\Vert$ small.
If we choose $u_0, u_1$ as to satisfy the conclusions of the theorem on
 the rapid flows, we can obtain vector fields
 $\xi_0, \xi_1$ with small norm  for
which the properties (\ref{f:odin}), (\ref{f:nul}) hold already for $T$
 sufficiently small.
Therefore the modified vector field $v'$
will be $C^0$ close to $v$.
$\qs$

\subsection{Ranging systems and
homological gradient descent (first version)}
\lb{su:rshgd}

As we have already mentioned, the map $\stv$ is not everywhere defined.
Still it is sometimes possible to define an analog of
"the homomorphism, induced by $\stv$
in homology".
More precisely, in some cases (for example, when
$v$ satisfies (RP))
one can define
 a homomorphism
$H(v):H_*(\daw\sm B_1, A_1)\to H_*(\dow\sm B_0, A_0)$,
where $A_i, B_i$ are some disjoint subsets of $\pr_i W$.
The construction of $H(v)$ was given in \cite{pastpet}, \S 4,
 and we shall briefly
expose it here.
We would like to stress that this homomorphism is {\it not}
induced by a continuous map of pairs.

\begin{defi}\lb{d:ransys}
Let $\Lambda =\{\lambda_0,...,\lambda_k\}$ be a finite set of
regular values of $f$, such that $\lambda_0 =a,\lambda_k =b$,
and for every $i$ we have $\lambda_i <\lambda_{i+1}$
and there is exactly one critical value of $f$ in
$[\lambda_i ,\lambda_{i+1}]$. The values $\lambda_i ,\lambda_{i+1}$
will be called {\it  adjacent}. The set of pairs $\ran$ is called
\it ranging
system
for $(f,v)$
\rm
 if
\begin{enumerate}
\item[(RS1)] $\forall \lambda\in \Lambda$ we have: $A_\lambda$ and
$B_\lambda$ are disjoint compacts in $\fpr f{\lambda}$.

\item[(RS2)] Let $\lambda,\mu\in\Lambda$ be adjacent. Then for
 every $p\in S(f)\cap \fpr f{[\lambda,\mu]}$ either

i) $D(p,v)\cap \fpr f{\lambda}\subset \Int A_\lambda$

 or

ii) $D(p,-v)\cap \fpr f{\mu}\subset \Int B_\mu$.

\item[(RS3)]  Let $\lambda,\mu\in \Lambda$ be adjacent. Then
$\stind {(-v)}\mu\lambda (A_\mu)\subset \Int A_\lambda ~\text{and}~
\stind v\lambda\mu (B_\lambda)\subset
\Int B_\mu.
 \qt$
\end{enumerate}
\end{defi}

Ranging systems have the following properties
(see \cite{pastpet}, \S 4.3).

\been\item[(1)]
If $\ran$ is a ranging system for $(f,v)$, then
 for every $f$-gradient $w$, sufficiently $C^0$-close to
$v$, $\ran$ is a ranging system for $(f,w)$.

\item[(2)]
 Let $N$ be a submanifold of $\daw\setminus B_b$
   such that
$    N\setminus \text{\rm Int}~
    A_b$ is compact.

Then $N'=\stind vba (N)$ is a submanifold of $\dow\setminus B_a$
  such that
$     N' \setminus \text{\rm Int}~ A_a$     is compact.
\item[(3)]
There is a homomorphism
$\qquad H(v):H_*(\daw\setminus B_b,A_b)\longrightarrow
\break
H_*(\dow\setminus B_a,A_a)$, such that
\begin{enumerate}
\item If $N$ is an oriented submanifold of $\daw$, satisfying the
hypotheses of (2), then
$H(v)([N])=[\stind vba (N)]$.

\item There is an $\epsilon>0$ such that for
 every $f$-gradient $w$ with $\nr w-v\nr <\epsilon$
we have $H(v)=H(w)$.
     \enen
\enen

We shall explain here the main idea of the construction of $H(v)$
in the following simple case: $f$ has only one critical point $p$,
and $D(p,v)\cap \dow\sbs\Int A_a$.
Moreover, we shall assume $B_a=\emp, B_b=\emp$.
Let
$U$ be a small compact \nei~ of $D(p,-v)\cap\daw$, \sut~
$\stv(U)\sbs\Int A_a$.
Let $x\in H_*(\daw\sm B_b, A_b)$.
Consider the  image of $x$ in
$H_*(\daw, A_1)$
and reduce it further to
$x'\in H_*(\daw, A_b\cup U)$.
Applying excision, we obtain an element
$\bar x\in H_*\(\daw\sm D(p,-v), A_b\cup U\sm D(p,-v)\)$.
Any singular chain representing $\bar x$ is in the domain
of definition of
$\stv$.
Apply $\stv_*$ and obtain an element in $H_*(\dow, A_0)$,
which is by definition $H(v)(x)$.
(The reader will  recognize here the basic idea of the
construction
of homological gradient descent from Subsection \ref{su:negcohgd}.)
\newpage

\section{Condition $(\gC)$ and operator of homological gradient descent}
\lb{s:ccohgd}

We need  some more terminology, we introduce it in Subsection
\ref{su:hlfc}.
In the second part of this Subsection
we recall the classical construction of Morse complex.
The material of this subsection is
not necessary for understanding of the rest of
the present Chapter, we placed it here just because it is natural
from the logical point of view.

In Subsection \ref{su:cc} we introduce the condition
$(\gC)$, which is a close
analog of
(RP), but it has some technical advantages over (RP). In particular,
the set of all $f$-gradients satisfying $(\gC)$ is
$C^0$-open in the set of all gradients.
In \su~ \ref{su:hgdsv} we introduce a new version of
 the operator of homological
gradient descent.

{\it Terminology}:
In the rest of \S 4
$W$ is a riemannian cobordism of dimension $n$.
\pa

\subsection{Handle-like filtrations of cobordisms
and Morse complexes}
\lb{su:hlfc}

In this subsection $\phi:W\to[a,b]$ is an
ordered Morse function with ordering sequence
$a_0=a<... <a_{n+1}=b$, and $v$ is a $\phi$-gradient.
Assume that $v$ is $\d$-separated.
Let $\nu\in]0,\d]$ and
$s$ be an integer between $0$ and $n+1$. Set

\begin{gather*}
W^{\{\leq s\}}
=\phi^{-1}([a_0, a_{s+1}]); \quad W^{\{\geq s\}}=
\phi^{-1}([a_s,a_{n+1}])\\
W^{[\leq s]} (\nu)
= D_\nu(\indl s ;v)\cup \pr_0 W;\\
W^{\lc s\rc}(\nu)=
 \fii 0s
\cup D_\nu(\inde s;v)
\end{gather*}

We shall reserve for  all these  filtrations of
$W$ the generic name "handle-like filtrations".
Of course these definitions are valid with larger
 assumptions on $v$ and $\phi$,
but only with our assumptions these filtrations
have the nice homological properties.
These properties
 are collected in the following proposition.
For a  critical point $p\in S_k(\phi)$ we denote by
$d(p)$ the embedded disc
$D(p,v)\cap\phi^{-1}([a_k,a_{k+1}])$.
For each critical point $p\in S(\phi)$
choose an orientation of the descending disc $D(p,v)$. Let
$[p]$ denote the image in the group
$H_k(\wk, \wkm)$
of the fundamental class of the pair
$(d(p), \pr d(p))$.

\bepr\lb{p:fil}
Let $s\in\NNN$ and $\nu\in ]0,\d]$.
\begin{enumerate}
\item The  inclusions
$W^{[\leq s]} (\nu) \sbs
\wasn\sbs
W^{\{\leq s\}} $
are homotopy equivalences.
\item
The inclusions of pairs
$$\bigg(W^{[\leq s]} (\nu),  W^{[\leq s-1]} (\nu)\bigg)
\sbs
\bigg(W^{\lc s\rc},W^{\{\leq  s-1\}}\bigg)
\sbs
\bigg(
W^{\{\leq s\}} , W^{\{\leq s-1\}}\bigg)
 $$
are homotopy equivalences.
\item
$H_*\(W^{\{\leq s\}} , W^{\{\leq s-1\}}\)
=0$ if $*\not= s$.
The group
$H_s\(W^{\{\leq s\}} , W^{\{\leq s-1\}}\)$
is a free abelian group generated by the elements
$[p], p\in S_s(\phi)$.
\end{enumerate}
\end{prop}

Denote $H_s(\ws,\wsm)$ by $C_s$.
The boundary operator of the exact sequence of the triple
$
(\ws, \wsm, \wsmm)$
gives a \ho~
$\pr_s:C_s\to C_{s-1}$;
we have obviously
$\pr_s\circ \pr_{s+1}=0$
for every $s$.
Thus the graded group
$C_*$
endowed with the boundary operator $\pr_*$
is a free chain complex.
In the case when $v$ satisfies the \TA~one can give another
construction of this complex, which is known as
{\it Morse complex}.
Namely, let $k$ be a positive integer, and let
$p\in S_k(\phi), q\in S_{k-1}(\phi)$.
Denote by
$\G(p,q;v)$
the set of the integral curves of $(-v)$
joining $p$ with $q$, where each integral curve is considered up
to a reparametrization.
It is not difficult to deduce from \TA, that
$\G(p,q;v)$
is a finite set.
The chosen orientations of the descending discs allow to attribute
to each
$\g\in \G(p,q;v)$
a sign
$\ve(\g)\in\{-1,1\}$.
Now let $C_k(v)$ be the free abelian group freely generated by
critical points of $\phi$ of index $k$.
Define a \ho~
$\pr'_k(v):C_k(v)\to C_{k-1}(v)$
setting
$\pr_k'(v)(p)=\sum_{q\in S_{k-1}(\phi)} n(p,q;v) q$.
It turns out that the \ho
~$J_k:C_k(v)\to C_{k}$
sending $p$ to $[p]$
satisfies
$J_k\circ\pr'_{k+1}(v)=\pr_{k+1}\circ J_{k+1}$;
thus
$(C_*(v), \pr_k(v))$
is a chain complex, isomorphic to $C_*$.

To define an "equivariant version" of the Morse complex, consider
a regular covering
$\PP:\wi W\to W$
with structure group $G$.
The inverse image of a subset $A\sbs W$
in $\wi W$
will be denoted by $\wi A$.
In addition to the orientation of descending discs, choose for
 every point
$p\in S(\phi)$
a lifting $\wi p$ of $p$ to $\wi W$, that is a point
$\wi p$ with
$\PP(\wi p)=p$.
Then the lifting of the disc $d(p)$
to $\wi W$ is defined.
Denote by
$[\wi p]$
the image in the group
$H_k(\tiwk, \tiwkm)$
of the fundamental class of the pair
$(d(\wi p), \pr d(\wi p))$.

\bepr\label{p:fileq}
$H_*(\tiws, \tiwsm)=0$, if $*\not= s$.
The group
$H_s(\tiws, \tiwsm)$
is a free $\ZZZ G$-module freely generated by
the elements $[\wi p], p\in S_s(\phi)$.
\enpr

Similarly to the above one gives a construction of chain complex
$(\wi C_*(v), \wi\pr_*(v))$
of free $\ZZZ G$-modules, which is freely generated in degree $p$
 by the critical points
of $\phi$
of index $p$, and the boundary operator is defined via counting
the $v$-trajectories joining the critical points.
The reader will find the details and the proofs of the cited
statements in
\cite{patou}, Appendix.

\subsection{Condition $(\gC)$}
\lb{su:cc}
\pa
{\it Terminology}
\lb{ssu:term}

Let $g:W\to[a,b]$ be a Morse function on a cobordism $W$
of dimension $n$, $v$ be a $g$-gradient.

\begin{defi}\label{d:desce}
Let $\l\in[a,b]$ be a regular value of $f$. Let $A\sbs\daw$.
We say, that {\it $v$ descends $A$ to $f^{-1}(\l)$},
if every
$(-v)$-trajectory starting in $A$ reaches $f^{-1}(\l)$.
Let $X\sbs W$. We say that
{\it $v$ descends $A$ to $f^{-1}(\l)$
without intersecting $X$}, if $v$ descends $A$ to
$f^{-1}(\l)$
and
$T(A, v)\cap X=\emp$.

We leave to the reader to the meaning of
{\it "$v$ lifts $B$ to $f^{-1}(\l)$" }
and
{\it "$v$ lifts $B$ to $f^{-1}(\l)$
without intersecting $Y$"}, where $B\sbs \dow, Y\sbs W).
$
\end{defi}

Recall the sets
$B_\d(\indl s;v), D_\d(\indl s;v)$
from Subsection \ref{su:dthd}
and introduce one more set
with similar properties:
set
$C_\d(\indl s; v)=W\sm B_\d(\indg {n-s-1};v)$.
This terminology is justified by the following lemma
which is proved by standard Morse-theoretical arguments.

\bele\lb{l:cdelta}
Assume that $v$ is $\d$-separated
\wrt~$\phi$ and the ordering sequence
$(a_0,...,a_{n+1})$.
Then for every $s$
we have:
$$
\dow\cup B_\d(\indl s;v)\sbs
\phi^{-1}\([a_0,a_{s+1}]\)
\sbs
C_\d(\indl s;v)\cup \dow
$$
and these inclusions are  homotopy equivalences.
$\qs$\enle

Now we can formulate the condition $(\gC)$.
\bede\lb{d:cc}

We say, that
{\it $v$ satisfies condition $(\gC)$}
if there are objects 1) - 4), listed below,
with the properties (A), (B) below.
\pa
{\it Objects:}
\pa
\begin{enumerate}
\item[1)] An ordered Morse function $\phi_1$
on $\pr_1 W$ with ordering sequence
$(\a_0,...,\a_n)$, and a $\phi_1$-gradient $u_1$.
\item[2)]  An ordered Morse function $\phi_0$
on $\pr_0 W$ with ordering sequence
$(\b_0,...,\b_n)$, and a $\phi_0$-gradient $u_0$.
\item[3)]  An ordered Morse function $\phi$
on $ W$ with ordering sequence
$(a_0,...,a_{n+1})$,
adjusted to $(f,v)$.
\item[4)]  A number $\d>0$.
\end{enumerate}
\pa
{\it Properties:}
\pa
(A)\quad  $u_0$ is $\d$-separated \wrt~ $\phi_0$,
$u_1$ is $\d$-separated \wrt~ $\phi_1$, $v$ is
$\d$-separated \wrt~ $\phi$.
\pa

 \begin{equation*}\stv \bigg(C_\d(\indl j;u_1)\bigg)
\cup
\bigg(D_\d\(\indl {j+1};v\)\cap \pr_0 W\bigg)
\sbs
B_\d(\indl j,u_0)
\mbox{ for every } j\tag{B1}
\end{equation*}

\begin{equation*}\st v\bigg(C_\d(\indl j;-u_0)     \bigg)
\cup
\bigg( D_\d(\indl {j+1}; -v)\cap \pr_1 W\bigg)
\sbs
B_\d(\indl j;-u_1)
\mbox{ for every } j\tag{B0}
\end{equation*}
$\qt$
\end{defi}

The set of all $f$ gradients satisfying $(\gC)$ will be denoted by
$\GgC(f)$. Recall that the set of all $f$-gradients satisfying
Transversality
Condition is denoted by $\GT(f)$. The intersection $\GgC(f)\cap \GT(f)$
will be denoted by $\GgCT(f)$.
\pa
{\it Comments:}
\begin{enumerate}\item
As we have already mentioned in the introduction,
the condition (B) is an analog of cellular approximation condition.
Indeed, the sets $C_\d(\indl j;u_1)$ and $B_\d(\indl j;u_0)$
are "thickenings" of
$D(\indl j;u_1)$, resp. of $D(\indl j; u_0)$, the first one
being "thicker" than the second. Thus the condition
(B1)
requires that $\stv$ sends a certain thickening of a
$j$-skeleton of $\daw$
to
a certain thickening of a $j$-skeleton of $\dow$.
Warning: (B1) actually requires more than that:
for every $j$
the soles of $\d$-handles of $W$ of indices $\leq j+1$
must belong to the
set
$B_\d(\indl j;u_0)$.

\item
The condition $(\gC)$ is stronger than the condition
$(\CC)$ from \cite{pator}, \S 1.2. This is obvious, since
$(\gC)$
is formulated similarly to $(\CC)$, except that in $(\CC)$ we
have demanded
$\stv (\phi_1^{-1}\([a_0, a_{j+1}]\)
\sbs
B_\d(\indl j;u_0)$. The condition $(\CC)$ is in turn stronger, than
the condition (RP) from the previous Section, this is a bit less obvious,
and we shall not prove it here.

\item
If we are given a $\phi_0$-gradient $u_0$, a $\phi_1$-gradient $u_1$
a $\phi$-gradient $v$, which satisfy Almost
Transversality Condition,
then the condition (A) is always true if
$\d$ is sufficiently small.

\end{enumerate}

\begin{theo}\lb{t:cc}
$\GgC(f)$ is open and dense in $\GG(f)$ \wrt~  $C^0$ topology.
Moreover, if $v_0$ is any $f$-gradient then one can choose a $C^0$
 small perturbation $v$
of $v_0$ \sut~ $v\in\GgC(f)$ and
$v=v_0$ in a \nei~ of $\pr W$.
\end{theo}

The proof  occupies the rest of the section and
is subdivided into 2 parts: openness and density.
\pa

{\it $C^0$-openness of $\GgC(f)$}\lb{ssu:cogct}
\pa
Let $v$ be an $f$-gradient satisfying $(\gC)$.
Choose the corresponding functions
$\phi_0, \phi_1, \phi$, their gradients
$u_0, u_1$, and a number $\d>0$ satisfying (A) and (B0), (B1).
Fix now the string
$\SS=(\d,u_0, u_1,\phi_0, \phi_1,\phi)$.
For a $\phi$-\gr~ $w$ we shall denote the condition (A) with
respect to $w$ and $\SS$
by (A)($w$). Similarly, we shall denote
 the conditions  (B1), (B0) \wrt~ $w$ and $\SS$
by  (B1)($w$), (B0)($w$). We know that
(A)($v$)$\&$ (B1)($v$) $\&$ (B0)($v$)
 holds; we shall prove that
(A)($w$)$\&$ (B1)($w$) $\&$ (B0)($w$)             holds
for every $w$ sufficiently $C^0$ close to $v$.

 Let $w$ be another $f$-gradient. Lemma 1.6 and Corollary 
1.7 of \cite{pastpet} imply that
if $\Vert w-v\Vert$ is sufficiently small,
then $w$ is
a $\phi$-gradient and is
 also $\d$-separated, so the condition A($w$) is
satisfied.

It is convenient to reformulate the condition (B1)($w$)$\&$(B0)($w$).

Introduce three new conditions:
\been
\item[$\b_1(w)$:]
For every $j$ we have:
$w$ descends $C_\d(\indl j;u_1)$ to $\phi^{-1}(a_{j+1})$
without intersecting $\cup_{\ind p\geq j+1} D_\d(p)$.

\item[$\b_2(w)$:]
For every $j$ we have:
$w$ lifts $C_\d(\indl {n-j-2};-u_0)$ to $\phi^{-1}(a_{j+2})$
without intersecting $\cup_{\ind p\leq j+1} D_\d(p)$.

\item[$\b_3(w)$:]
For every $j$ the sets
$R_j^+(w)=\stind {(-v)}b{a_{j+2}} (C_\d(\indl j;u_1))$

and

$R_j^-(w)=\stind {v}a{a_{j+2}} (C_\d(\indl {n-j-2};-u_0))$

are disjoint.
\enen

Note that $\b_1(w)$ and $\b_2(w)$ are in a sense dual to each other,
that is: $\b_2(w)=\b_1(-w)$.

\bele
Assume that A($w$) holds. Then
(B1)($w$)$\&$(B0)($w$)$\Leftrightarrow
\b_1(w)\&\b_1(w)\&\b_3(w)$
\enle
\Prf $\Rightarrow$

The condition $\b_1$ follows from (B0), $\b_2$ follows from (B1).
To prove $\b_3$, note that if
$R_j^+(w)\cap R_j^-(w)\not=\emp$, then there is
a $(-v)$-trajectory joining
$x\in C_\d(\indl j;u_1)$ with
$y\in C_\d(\indl {n-j-2}; -u_0)=
\dow\sm B_\d(\indl j; u_0)$.
And this contradicts (B1).

$\Leftarrow$

 $\b_1$ says that for every $j$ we have:
$$
\bigg(\bigcup_{\ind p\geq j+1} D_\d(p,-v)\bigg)
\cap\daw\sbs
\daw\sm C_\d(\indl j;u_1)=B_\d(\indl {n-2-j}; -u_1).
$$
That is,
$D_\d(\indl {n-1-j}; -v)\cap \daw\sbs B_\d(\indl {n-2-j};-u_1)$.
Further, $\b_3$ implies that for every $j$
every $v$-trajectory starting in
$C_\d(\indl {n-j-2}; -u_0)$
and reaching $\daw$, cannot intersect $\daw$ at a point of
$C_\d(\indl j;u_1)$. Therefore it intersects
$\daw$ at a point of
$B_\d(\indl {n-2-j};-u_1)$.
Therefore (B0)($w$) holds. The proof of (B1)($w$) is similar. $\qs$

Returning to the proof of $C^0$-openness of $(\gC)$, note that the
conditions
$\b_1(w)$ and $\b_2(w)$ are obviously open.
We know that $\b_1(v)$ and $\b_2(v)$ hold, therefore
for every $j$ the sets
$R_j^+(v), R_j^-(v)$
are compact. Since they are disjoint (by $\b_3(v)$),
we can choose two disjoint open subsets
$U_j, V_j$ of
$\phi^{-1}(a_{j+2})$ \sut
~$R_j^+(v)\sbs U_j, R_j^-(v)\sbs V_j$. Then for every
$w$
sufficiently $C^0$-close to $v$ and for every
$j$
we have (see \cite{pastpet}, Corollary 5.6):
$$
\stind {(-w)}b{a_{j+2}} \(C_\d(\indl j;u_1)\)\sbs U_j, \quad
\stind wa{a_{j+2}}\( C_\d(\indl {n-2-j}; -u_0)\)\sbs V_j
$$
and $\b_3(w)$ holds. Therefore
$\b_1(w)\&
\b_2(w)\&
\b_3(w)$
is $C^0$ open and we have proved $C^0$ openness of $(\gC)$.
\pa
{\it $C^0$ density.}
\pa
The proof is very close to that of $C^0$ density of (RP).
We shall not give it here in details; we just indicate that during
 the proof
of $C^0$ density of (RP)
we have constructed a perturbation $v'$ of a given gradient $v$,
which satisfies actually the condition $(\gC)$.
(Indeed, return to the paragraph just after (\ref{f:compo}).
The diffeomorphism $\Phi^{-1}_1$
sends the whole of $\daw\sm A_{n-2-m}^+(\d)$ to the set
$A^-_m(\d)$.
Since $u_0$ and $u_1$ are $\d$-separated, the set
$(\daw)^{\{\leq m\}}$ is in $\daw\sm A_{n-2-m}^+(\d)$
and we have checked the property (B1). )$\qs$

\subsection{ Homological gradient descent(second version)}
\lb{su:hgdsv}

In this subsection we  define  a version of the homological
gradient descent operator from \su~ \ref{su:rshgd}.
This version is better suited for the study of Morse-type
filtrations on cobordisms, which is the subject of
the next section.
Assume that $v$ satisfies condition $(\CC)$
\wrt~ $\d,\phi_0, \phi_1, u_0, u_1$.
We have the corresponding filtrations
$(\daw)^{\{\leq s\}},
(\dow)^{\{\leq s\}}$.
Denote $\phi^{-1}(\a)$ by $V_\a$.

\begin{theo}\label{t:consh}
For every $s$  there is a homomorphism
$$\HH_s(-v): H_*(V_b^{\{\leq s\}}, V_b^{\{\leq s-1\}} )
\to
H_*(V_a^{\{\leq s\}}, V_a^{\{\leq s-1\}} )
$$
with the following properties:
\begin{enumerate}\item[1)]
Let $N$ be an oriented \sma~ of $V_b$,
\sut~$N\sbs V_b^{\{\leq s\}}$
and $N\sm\Int V_b^{\{\leq s-1\}}$ is compact.
Then the manifold $N'=\stind {(-v)}ba (N)$ is in $V_a^{\{\leq s\}}$ and
$N'\sm \Int V_a^{\{\leq s-1\}}$
is compact and the fundamental class of $N'$ modulo
$V_a^{\{\leq s-1\}}$ equals to
$\HH_s(-v)\([N]\)$.
\item[2)]
There is an $\e>0$ \sut~for every $f$-gradient $w$ with
$\Vert w-v\Vert\leq \e$ and every $s$
we have:
$\HH_s(-v)=\HH_s(-w)$.
\end{enumerate}
\end{theo}

\Prf
1) Construction of $\HH_s(-v)$.

The homomorphism $\HH_s(-v)$ will be defined as the composition of two
homomorphisms: $\HH 1_s(-v)$ and $\HH 0_s(-v)$.
(Intuitively, $\HH 1$ corresponds to the descent from the level
$b$ to the
level $a_{s+1}$,
and $\HH 0$ corresponds to the descent from
$a_{s+1}$ to $a$.)
We shall denote the maps
$\stind {(-v)}b{a_{s+1}}$
and
$\stind {(-v)}{a_{s+1}}a$
by $\stexp {(-v1)}$, and respectively $\stexp {(-v0)}$. Every $(-v)$-
trajectory
starting in $\vbs$ reaches $\vass$.
Therefore $\vbs$ is in the domain of
$\stexp {(-v1)}$
and
$\stexp {(-v1)}$
defines a homeomorphism of the pair
$(\vbs, \vbsm)$ to its image, which is a pair of compact
 subsets of $V_{a_{s+1}}$.

\begin{defi}
\lb{d:ykv}
Let $k$ be an integer, $0\leq k\leq n$ .
Denote by $Y_k(v)$ the set of all
$y\in\phi^{-1}([a,a_{k+1}])$
\sut~either

i) $\g(y, \cdot; -v)$ reaches $V_a$ and intersects it at a point
$z\in \vakm$ or

ii) $\lim_{t\to\infty}\g(y,t;-v)=p$, where $p\in S(f)$.$\qt$
\end{defi}

In other words, $Y_k(v)=
T(V_a^{\{\leq k-1\}}, -v\!\mid\! W')\cup D(-v\mid W')$,
where $W'=\phi^{-1}\([a,a_{k+1}]\)$.\label{wprim}
The set $Y_k(v)$ is compact
(by \cite{pastpet}, Lemma 2.7).
Set

\begin{equation}\lb{f:triodin}
A =Y_{s+1}(v)\cap V_{a_{s+1}},\quad
B =Y_s(v)\cap V_{a_{s+1}},\quad
C =\(\cup_{\ind p\leq s} D(p,-v)\)\cap V_{a_{s+1}}
\end{equation}

Then
$C$ is closed and
$\bar C\sbs \Int B$ and $A\sps B \sps \Int B\sps C$.
Further, the map $\stexp {(-v1)}$ defines a continuous map of pairs
$$\( V_b^{\{\leq s\}}, V_b^{\{\leq s-1\}})\to (A,B)
$$
and the map
$\stexp {(-v0)}$
defines a continuous map of pairs
$$
(A\sm C, B\sm C)\to(V_a^{\{\leq s\}}, V_a^{\{\leq s-1\}})
$$

Now denote by $Exc$ the excision isomorphism
$H_*(A\sm C, B\sm C)\to H_*(A,B)$
and set
$$\HH_s(-v)=\(\stexp{(-v0)}\circ Exc^{-1}\circ\stexp{(-v1)}\)_*$$

2) Proof of the properties of $\HH_s(-v)$.

The proof is reduced to the properties of the operator $H(v)$
associated to
ranging systems (see \su~\ref{su:rshgd}).

Namely we construct a ranging system
$\RR=\ran$ for $(\phi,v)$, \sut~
the \ho~$H(v)$, associated to
$\RR$ is isomorphic to $\HH(-v)$. To
explain the construction of $\RR$ we make
a simplifying assumption:
$\phi$ has only two critical values:
one in $[a, a_{s+1}]$ and the other in $[a_{s+1}, b]$.
(In the end of the proof we indicate how to
generalize it to the general situation.)
We set
\begin{gather*}
A_b=V_b^{\{\leq s-1\}}, A_a= V_a^{\{\leq s-1\}},\\
 B_a=V_a^{\{\geq s+1\}},  B_b=V_b^{\{\geq s+1\}}
\end{gather*}
Denote $a_{s+1}$ by $\tau$ for brevity.
To define  $A_\tau, B_\tau $,
let $(\a_i),(\b_i)$ be ordering sequences for $\phi_1$, resp. $\phi_0$,
and choose $\e>0$ so small that:

\begin{gather*}
D_\d(\indl {n-s-2}, -u_1)\sbs \phi_1^{-1}(]\a_{s+1}+\e, \a_n]),\\
\a_s+\e<\a_{s+1},\\
\stv (\phi^{-1}_1([\a_0, \a_s+\e]))\sbs\phi_0^{-1}([\b_0, \b_s[)\\
\st v (\phi_0^{-1}([\b_{s+1}, \b_n]))\sbs
\phi_1^{-1}(]\a_{s+1}+\e, \a_n])
\end{gather*}

Denote by $L_s(v)$ the set of all points
$y\in\phi^{-1}\([a_{s+1}, b]\)$, \sut
~either
\been\item[i)]
$\g(y,\cdot;v)$ reaches $V_b$ and intersects it at a
point $z\in\phi^{-1}\([\a_{s+1}+\e, \a_n]\)$, or
\item[ii)] $y\in D(p,v)$, where $p\in S(f)$.
\enen

Now set
$$
A_\tau= \stexp {(-v1)} (\phi_1^{-1}([\a_0, \a_0+\e])),
\quad
B_\tau = L_s(v)\cap V_\tau
$$

It follows from $(\gC)$ that
$\{(A_\l, B_\l)\}_{\l\in\{a,\tau, b\}}$
is a ranging system for $(\phi,v)$.

There are inclusions of pairs
\begin{gather*}
j_b:(V_b^{\{\leq s\}}, V_b^{\{\leq s-1\}})\sbs (V_b\sm B_b, A_b)\\
j_a:(V_a^{\{\leq s\}}, V_a^{\{\leq s-1\}})\sbs (V_a\sm B_a, A_a)
\end{gather*}

and these inclusions induce isomorphisms in \hog.

Moreover, checking through the definition of
the \ho~$H(v)$, associated with
$\ran$, it is easy to see, that
we have
$H(v)\circ(j_b)_*=(j_a)_*\circ\HH(-v)$.
Now the properties of $\HH(-v)$ follow immediately from
those of $H(-v)$.

In the case when there is more than one critical value
in $[a_{s+1}, b]$ and
in $[a,a_{s+1}]$, the argument is
similar:
one constructs a ranging system $\RR=\ran$, \sut~the operator
$H(v)$, associated to $\RR$ is isomorphic to $\HH(-v)$.
In this ranging system we can define
$A_\l, B_\l$ by the same formulas as above for $\l=a,b,\tau$.
To define $A_\l, B_\l$ for other values of $\l$, choose first
 a set $\L=
\{\l_i\}$
of regular values
\sut~ $\l_0=a, \l_N=b, \l_i=\tau$ and for every $i$ there is only one
critical value of $\phi$ in
$[\l_i, \l_{i+1}]$.
For $i\leq j\leq N$
we define 
$A_{\l_j}$
as some thickening of 
$\stind {(-v)}b{\l_j}\(V_b^{\{\leq s-1\}})$.
The set 
$B_{\l_j}$
for $i\leq j\leq N$
is defined similarly  to $L_s(v)$
(see the definition above).
The definition of the sets
$A_\l, B_\l$
for the other values of $\l$,
as well as the proof that
$H(v)$ is isomorphic to
$\HH(-v)$ will be left to the reader.
$\qs$

\subs{Equivariant homological gradient descent}
\lb{su:ehgd}

Let $p:\wh W\rTo^H W$ be a regular covering of $W$
 with the structure group $H$.
For $A\sbs W$ we denote $p^{-1}(A)$ also by $\wh A$.
The shift along the $(-v)$-trajectories defines a \dfm~
$\stexp {(-\wh v)}:\wh {(\dow \sm K_1)}\to \wh {(\dow\sm K_0)}$
(see the page
\pageref{stind} for the definition of $K_1, K_0$)
which commutes with the right action of $G$. As before we abbreviate
$\stexp {(-\wh v)}(X \sm \wh K_1)$
to
$\stexp {(-\wh v)}(X)$.
Assume that $v$ satisfies condition
$(\gC)$ \wrt~$\d, \phi_0, \phi_1, u_1, u_0$. Denote $\phi^{-1}(\a)$ by
$V_\a$.
Then we have a corresponding equivariant version of homological gradient
 descent operator.
The properties of this operator are the contents of the next theorem.
The proof is completely similar to the proof of
Theorem \ref{t:consh} and will be omitted.

\beth\lb{t:consheqv}
For every $s$ there is a \ho~
of right $\ZZZ H$-modules
$$
\wh\HH_s(-v): H_*(\wh V_b^{\{\leq s\}}, \wh V_b^{\{\leq s-1\}} )
\to
H_*(\wh V_a^{\{\leq s\}}, \wh V_a^{\{\leq s-1\}} )
$$
having the following property:

Let $N$ be an oriented submanifold of $\wh V_b$,
\sut~$N\sbs \wh V_b^{\{\leq s\}}$ and
$    N\setminus \text{\rm Int}~
    \wh V_b^{\{\leq s-1\}} $ is compact.

Then $N'=\stind {(-\wh v)}ba (N)$ is a
submanifold of $\dow\setminus B_a$
  such that
$N'\sbs \wh V_a^{\{\leq s\}}$ and
$     N' \sm V_a^{\{\leq s-1\}} $     is compact,
and the fundamental class of $N'$ modulo
$V_a^{\{\leq s-1\}} $ equals to $\wh \HH_s(-v) ([N])$. $\qs$
\enth

\subs{Cyclic cobordisms and the condition $(\gC\YY)$}
\label{su:cyc}

In this section we develop some techniques which will be used
in \S\S 5-7. 

\begin{defi}
A {\it cyclic cobordism}
is a riemannian cobordism $W$ together
with
an isometry $\Phi:\dow\to \daw$.
\end{defi}
Let $f:W\to[a,b]$ be a Morse function
on a cyclic \co~ $W$ of dimension $n$,
$v$ be an $f$-gradient. We say that $v$ satisfies
$(\gC\YY)$ if $v$ satisfies condition $(\gC)$
(see \su~\ref{su:cc}),
and, moreover, the Morse functions $\phi_1, \phi_0$ and their gradients
$u_1, u_0$ from the definition \ref{d:cc}
can be chosen so that
$\phi_1\circ\Phi=\phi_0, u_1=(\Phi)_*(u_0)$.
The set of all $f$-gradients satisfying
$(\gC\YY)$
will be denoted by
$\GgCY(f)$. The intersection
$\GgCY(f)\cap\GT(f)$ will be denoted by
$\GgCYT(f)$.

\beth\lb{t:cyc}
$\GgCY(f)$ is open and dense in $\G(f)$ \wrt~  $C^0$ topology.
Moreover, if $v_0$ is any $f$-gradient then one can choose a $C^0$
 small perturbation $v$
of $v_0$ \sut~ $v\in\GgCY(f)$ and
$v=v_0$ in a \nei~ of $\pr W$.
\end{theo}
\Prf It goes exactly as the proof of Theorem \ref{t:cc}.
To prove $C^0$-openness of the set
$\GgCY(f)$ we do not need to change whatever it be in the first
part of the proof of Theorem \ref {t:cc}.
For the proof of $C^0$ density return again to
the proof of $C^0$ density of (RP). (Theorem \ref{t:rpdense}).
All what we have to notice is that we can choose
$u_1=\Phi_*(u_0), \xi_1=\Phi_*(\xi_0),
\phi_0=\phi_1\circ\Phi$. $\qs$

For an $f$-gradient $v$ satisfying $(\gC\YY)$
we can turn the operator of homological gradient descent
to an endomorphism of an abelian group.
Namely, set
$h_s(-v)=\HH_s(-v)\circ\Phi_*$;
then
$h_s(-v)$ is an endomorphism
of
$H_*(\mws, \mwsm)$.

\newpage

\section{Morse-type filtrations of cobordisms}
\label{s:mtfc}

Let $W$ be a cobordism, $\phi:W\to[a,b]$ be an ordered Morse
 function on
$W, v$ be a $\phi$-gradient.
Assume that $W$ is riemannian, and that $v$ is $\d$-separated.
We have considered in the subsection
\ref{su:hlfc}
the handle-like filtrations of $W$, associated to $v$.
One can reconstruct the \hog~ of $(W,\dow)$ from
the \hog~ of the successive factors of each of these filtrations.
But we can not reconstruct the absolute
\hog~ of $W$ from these factors.
(Indeed, consider the following example:
$W=V\times [0,1]$. Here all the terms of the filtration
are \hot~ equivalent to $V$, and the factors are contractible.)
In this section we introduce another filtration of $W$, which is
in a sense adjusted to
given filtrations of $\dow$ and $\daw$, and from
which
one can reconstruct the \hog~ of $W$.
In order to construct such filtration one must
 impose the condition $(\gC)$ on $v$.

\subsection{Definition of Morse-type filtration}
\label{su:dmtf}

Let
$W$ be a riemannian cobordism of dimension $n$,
$f:W\to[a,b]$ be a Morse function on   $W$
and $v$ be an $f$-gradient.
Assume that $v$ satisfies condition $(\gC)$.
Recall from Subsection \ref{su:cc}
that the condition $(\gC)$ requires the existence of
ordered
 Morse functions
$\phi:W\to[a,b], \phi_1:\pr_1 W\to\RRR, \phi_0:\dow\to\RRR$,
a $\d$-separated $\phi_0$-gradient $u_0$ and a $\d$-separated
$\phi_1$-gradient $u_1$, satisfying the conditions
(A), (B1), (B0) from
Subsection \ref{su:cc}.
We obtain therefore three handle-like filtrations
associated with the ordered Morse functions $\phi_0, \phi_1$ and $\phi$:
 the filtration
$\pws$ of $\daw$, the filtration $\mws$ of $\dow$ and the filtration
$\ws$ of $W$.

Recall the set
$$
Y_k(v)=
\bigg(T(\mwkm, -v)\cup D(-v)\bigg)\cap \phi^{-1}\([a, a_{k+1}]\)
$$
 from the previous section
(definition \ref{d:ykv}).
In the present section we deal only with one $f$-gradient, so we shall
 abbreviate
$Y_k(v)$ to $Y_k$.

\bede\lb{d:urez}
Let $A\sbs W$ and $c,d\in[a,b]$.
We denote
$A\cap\phi^{-1}\([c,d]\)$
by
$A\scc {[c,d]}$. For $\m\in[a,b]$
we denote
$A\cap\phi^{-1}(\mu)$
by$A\scc{\mu}$. $\qt$
\end{defi}

Set
\begin{equation}
Z_k=
(T(\pwkm, v))\scc{[a_{k+1},b]}
\end{equation}

Since  $(-v)$  descends
$\pwkm$ to  $\phi^{-1}(a_{k})$
(this follows from $(\gC)$)
the set $Z_k$ is homeomorphic to the product
$\pwkm\times [0,1]$.

Set for $k\geq 0$:
\begin{equation}
\begin{gathered}
\wmok=
\mwk\cup\pwk\cup Y_k\cup Z_k\\
\Wmok=
\pwk\cup \(Y_{k+1}\cap \phi^{-1}\([a, a_{k+1}]\)\)\cup Z_k
\end{gathered}\lbl{f:wkdef}
\end{equation}
Both $\{\Wmok\}_{k\in\NNN}$ and $\{\wmok\}_{k\in\NNN}$
are filtrations of $W$,
and $\wmok\sbs\Wmok$.
The set $\Wmok$ is obtained from $\wmok$ by adding
all the points in $\waa {k+1}$, lying on $v$-trajectories
starting in
$\mwk\sm\mwkm$.

(To visualize $\Wmok$ and $\Wmokm$
look at  Figures 1 and 2 below.)

\begin{setlength}{\unitlength}{0.15cm}

\begin{picture}(100,100)
\linethickness{0.07mm}

%ramka

\put(0,0){\line(1,0){45}}
\put(0,50){\line(1,0){45}}
\put(0,0){\line(0,1){50}}
\put(45,0){\line(0,1){50}}

%urovni
\put(0,30){\line(1,0){45}}
\put(0,12){\line(1,0){45}}
\put(0,20){\line(1,0){45}}

%zhirnye urovni
\linethickness{0.6mm}
\put(5,30){\line(1,0){8.33}}

\put(31.67,30){\line(1,0){8.33}}

\put(5,0){\line(1,0){35}}

\put(5,50){\line(1,0){35}}

\linethickness{0.07mm}

%vert skobki s formulami
\put(67,0){\line(1,0){2}}
\put(67,30){\line(1,0){2}}
\put(69,0){\line(0,1){30}}
\put(70,15){$W^{\{\leq k\}}$}

%ischo vert skobki s formulami
\put(50,0){\line(1,0){2}}
\put(50,20){\line(1,0){2}}
\put(52,0){\line(0,1){20}}
\put(53,10){$W^{\{\leq k-1\}}$}

%Vert cherta urovni f
\put(85,-5){\line(0,1){60}}
\put(85,0){\circle*{0,75} }
\put(85,12){\circle*{0,75} }
\put(85,20){\circle*{0,75} }
\put(85,30){\circle*{0,75} }
\put(85,50){\circle*{0,75} }

%urovni f
\put(87,0){$a$ }
\put(87,12){$a_{k-1}$ }
\put(87,20){$a_{k}$ }
\put(87,30){$a_{k+1}$ }
\put(87,50){$b$ }

%nadpis' urovni f
\put(80, 70){value of $\phi$}

%podpis'
\put(30,-20){FIGURE 1: $W^{\langle k\rangle}$}

%horiz skobki s formulami
\put(15,60){\line(0,1){2}}
\put(30,60){\line(0,1){2}}
\put(15,62){\line(1,0){15}}
\put(18,63){$\partial_1 W^{\{\leq k-2\}}$}

%escho horiz skobki s formulami
\put(10,70){\line(0,1){2}}
\put(35,70){\line(0,1){2}}
\put(10,72){\line(1,0){25}}
\put(18,73){$\partial_1 W^{\{\leq k-1\}}$}

%escho escho horiz skobki s formulami
\put(5,80){\line(0,1){2}}
\put(40,80){\line(0,1){2}}
\put(5,82){\line(1,0){35}}
\put(18, 83){$\partial_1 W^{\{\leq k\}}$}

%tolstye tochki
\put(5,50){\circle*{0,75} }
\put(10,50){\circle*{0,75} }
\put(15,50){\circle*{0,75} }
\put(30,50){\circle*{0,75} }
\put(35,50){\circle*{0,75} }
\put(40,50){\circle*{0,75} }

%nizhnyaya banka

\linethickness{0.6mm}
\put(5,0){\line(0,1){30}}
\linethickness{0.07mm}

\put(10,0){\line(0,1){20}}

\put(35,0){\line(0,1){20}}

\linethickness{0.6mm}
\put(40,0){\line(0,1){30}}
\linethickness{0.07mm}

%verhnyaya banka
\put(10,50){\line(1,-6){3,33}}
\put(10.025,50){\line(1,-6){3,33}}
\put(10.05,50){\line(1,-6){3,33}}
\put(10.075,50){\line(1,-6){3,33}}
\put(10.1,50){\line(1,-6){3,33}}
\put(10.125,50){\line(1,-6){3,33}}
\put(10.15,50){\line(1,-6){3,33}}
\put(10.175,50){\line(1,-6){3,33}}
\put(10.2,50){\line(1,-6){3,33}}

\put(15,50){\line(1,-6){5}}

\put(25,20){\line(1,6){5}}

\put(31.67,30){\line(1,6){3,33}}
\put(31.695,30){\line(1,6){3,33}}
\put(31.72,30){\line(1,6){3,33}}
\put(31.745,30){\line(1,6){3,33}}
\put(31.77,30){\line(1,6){3,33}}
\put(31.795,30){\line(1,6){3,33}}
\put(31.82,30){\line(1,6){3,33}}
\put(31.845,30){\line(1,6){3,33}}
\put(31.870,30){\line(1,6){3,33}}

\end{picture}

\end{setlength}

\vskip1.5in

\begin{setlength}{\unitlength}{0.15cm}

\begin{picture}(100,100)
\linethickness{0.07mm}

%ramka

\put(0,0){\line(1,0){45}}
\put(0,50){\line(1,0){45}}
\put(0,0){\line(0,1){50}}
\put(45,0){\line(0,1){50}}

%urovni
\put(0,30){\line(1,0){45}}
\put(0,12){\line(1,0){45}}
\put(0,20){\line(1,0){45}}

%zhirnye urovni
\linethickness{0.6mm}
\put(10,20){\line(1,0){10}}

\put(25,20){\line(1,0){10}}

\put(10,0){\line(1,0){25}}

\put(10,50){\line(1,0){25}}

\linethickness{0.07mm}

%vert skobki s formulami
\put(67,0){\line(1,0){2}}
\put(67,30){\line(1,0){2}}
\put(69,0){\line(0,1){30}}
\put(70,15){$W^{\{\leq k\}}$}

%ischo vert skobki s formulami
\put(50,0){\line(1,0){2}}
\put(50,20){\line(1,0){2}}
\put(52,0){\line(0,1){20}}
\put(53,10){$W^{\{\leq k-1\}}$}

%Vert cherta urovni f
\put(85,-5){\line(0,1){60}}
\put(85,0){\circle*{0,75} }
\put(85,12){\circle*{0,75} }
\put(85,20){\circle*{0,75} }
\put(85,30){\circle*{0,75} }
\put(85,50){\circle*{0,75} }

%urovni f
\put(87,0){$a$ }
\put(87,12){$a_{k-1}$ }
\put(87,20){$a_{k}$ }
\put(87,30){$a_{k+1}$ }
\put(87,50){$b$ }

%nadpis' urovni f
\put(80, 70){value of $\phi$}

%podpis'
\put(30,-20){FIGURE 2: $W^{\langle k-1\rangle}$}

%horiz skobki s formulami
\put(15,60){\line(0,1){2}}
\put(30,60){\line(0,1){2}}
\put(15,62){\line(1,0){15}}
\put(18,63){$\partial_1 W^{\{\leq k-2\}}$}

%escho horiz skobki s formulami
\put(10,70){\line(0,1){2}}
\put(35,70){\line(0,1){2}}
\put(10,72){\line(1,0){25}}
\put(18,73){$\partial_1 W^{\{\leq k-1\}}$}

%escho escho horiz skobki s formulami
\put(5,80){\line(0,1){2}}
\put(40,80){\line(0,1){2}}
\put(5,82){\line(1,0){35}}
\put(18, 83){$\partial_1 W^{\{\leq k\}}$}

%tolstye tochki
\put(10,50){\circle*{0,75} }
\put(15,50){\circle*{0,75} }
\put(30,50){\circle*{0,75} }
\put(35,50){\circle*{0,75} }

%nizhnyaya banka

\put(5,0){\line(0,1){30}}

\linethickness{0.6mm}
\put(10,0){\line(0,1){20}}
\linethickness{0.07mm}

\linethickness{0.6mm}
\put(35,0){\line(0,1){20}}
\linethickness{0.07mm}

\put(40,0){\line(0,1){30}}

%verhnyaya banka
\put(10,50){\line(1,-6){5}}

\put(15,50){\line(1,-6){5}}
\put(15.025,50){\line(1,-6){5}}
\put(15.05,50){\line(1,-6){5}}
\put(15.075,50){\line(1,-6){5}}
\put(15.1,50){\line(1,-6){5}}
\put(15.125,50){\line(1,-6){5}}
\put(15.15,50){\line(1,-6){5}}
\put(15.175,50){\line(1,-6){5}}
\put(15.2,50){\line(1,-6){5}}

\put(25,20){\line(1,6){5}}
\put(25.025,20){\line(1,6){5}}
\put(25.05,20){\line(1,6){5}}
\put(25.075,20){\line(1,6){5}}
\put(25.1,20){\line(1,6){5}}
\put(25.125,20){\line(1,6){5}}
\put(25.15,20){\line(1,6){5}}
\put(25.175,20){\line(1,6){5}}
\put(25.2,20){\line(1,6){5}}

\put(30,20){\line(1,6){5}}

\end{picture}
\end{setlength}

\vskip1.5in

\bere
\lb{r:infde}
\been\item
Here is an informal  description of
these filtrations.  If we consider the descending
discs $D(p,v)$ of critical points of $f$ of index $k$
as the "$k$-cells" of $W$, then $\wmok$
contains all the cells of
$\dow, \pr_1 W$ and  $W$ of dimensions $\leq k$, and the tracks
of the cells of $\pr_1 W$ of dimension $\leq k-1$.
The sets $\wmok, \Wmok$ are   thickenings of this "skeleton".
The condition $(\gC)$
guarantees that the boundary of the "cells" corresponding
 to $D(p,v), p\in S(\phi)$
belongs to the union of the thickened cells
of dimension $<\ind p$.
\item
The two filtrations are  \hot~ equivalent, see
Lemma \ref{l:qu} below.
The set $\wmok$ was introduced in \cite{pator},
and denoted there by $\Wmok$.
 The filtration $\{\Wmok\}_{k\in\NNN}$
has some technical advantages
over $\wmok$;
that is why we have  changed our notation
and why we shall consider mainly the filtration $\Wmok$
in the present paper.
\enen
\enre

We set $\wmok=\Wmok=\emp$ for $k<0$.

\bele\lb{l:qu}
The inclusion
$(\wmok, \wmokm)\rInto(\Wmok,\Wmokm)$
is a \hot~ equivalence.
\enle
\Prf
We have
\begin{equation}
\Wmok\sm\wmok
=
(Y_{k+1}\cap \wk)
\sm
(Y_k\cup\mwk)
\end{equation}
The function $\phi$ has no critical points in the domain
$A=\phi^{-1}([a,a_{k+1}])\sm Y_k$.
Therefore there is $T>0$, \sut~
for every~ $x\in A$
we have:
$\g(x,T;-v)\in \dow\cup Y_k$.
For $x\in W$ denote by $\tau(x)$
the moment when
$\g(x,\cdot;-v)$ reaches $\dow$
(if the trajectory never reaches
$\dow$, set $\tau(x)=\infty$.)
Set
\begin{gather}
\Womok  =\Wmok\cap \phi^{-1}([a,a_{k+1}])\\
\womok= \wmok\cap \phi^{-1}([a,a_{k+1}])
\end{gather}

Define now
a deformation
$H_t$ of the space
$\Womok  $
to itself
(where $ t\in[0,T]      $)
setting
\begin{equation}
H_t(x)=
\begin{cases}
\g(x,t;-v) & \mbox{ if } t\leq \min(T,\tau(x))  \\
\g(x,\min(T,\tau(x));-v) & \mbox{ if } t\geq \min(T,\tau(x))
\end{cases}
\end{equation}

It is easy to see that
$H_0=\id, H_T(\Womok)\sbs \womok$
and that
$\womok$ is $H_t$-invariant for all $t\geq 0$. Therefore
the inclusion
$\womok\rInto\Womok$
is a \hot~ equivalence.

Now it is not difficult to extend the deformation
$H_t$
to the whole of
$\Wmok$
(set $H_t(x)=x$
for
$x\in\phi^{-1}([a_{k+1}+\e,b])$
and glue appropriately two maps in the band
$\fifi {a_{k+1}}{a_{k+1}+\e}$.
The extended deformation will
provide the \hot~ inverse to the inclusion
$\wmok\rInto\Wmok$.

A homotopy inverse to the inclusion cited in the statement of the lemma
is constructed in the same way
(only one may need to choose $T$ larger).
We leave the details to the reader.
$\qs$

\bere
\label{re:shift}
The method used in the proof will be referred to as
{\it $(-v)$-shift}.
We shall use it several times in this section.
In particular, in the proof of the Lemma above
we have applied 
the $(-v)$-shift 
 in the domain
$\phi^{-1}([a,a_{k+1}])\sm Y_k$.      
\enre

\subsection{Homology of $\(\Wmok, \Wmokm\)$}
\label{su:hwkwkm}

Here we show that $\Wmok$ form a cellular filtration
of $W$, and  compute the \hog~ of
$\(\Wmok, \Wmokm\)$. The remark \ref{r:infde}
suggests that for $*\not= k$
the group
$H_*(\Wmok, \Wmokm)$
vanishes and that for $*=k$ this group is free
abelian with the base formed by the
\hog~ classes of the
corresponding cells.
In this subsection we show that it is indeed so.

To compute $H_*(\Wmok, \Wmokm)$
we shall first define  four homomorphisms with values  in
$H_k(\Wmok, \Wmokm)$.
We shall prove in Theorem \ref{t:mtf} that their images
 form a direct sum decomposition of this group.
\been
\item
For every $s$ the sets
$\mws, \pws$
are subsets of
$\Wmos$. Therefore there are homomorphisms, induced by the
corresponding inclusions:
\begin{gather}
C_k(u_0)=H_k\bigg(\mwk, \mwkm\bigg)\rTo^{\JJ_0} H_k(\Wmok, \Wmokm)
\quad\lbl{f:j0}\\
C_k(u_1)=H_k\bigg(\pwk, \pwkm\bigg)\rTo^{\JJ_1} H_k(\Wmok, \Wmokm)
\lbl{f:j1}
\end{gather}

\item
Every $(-v)$-trajectory starting in $\pwkm$
reaches $\phi^{-1}(a_k)$.  Therefore there is a continuous map
$$
G_k:\pwkm\times[a_k,b]\to \phi^{-1}\([a_k, b]\)
$$
which is a homeomorphism onto its image and for every $x$
the curve $t\mapsto G_k(x,t)$
is a reparameterized $v$-trajectory, and $\phi(G_k(x,t))=t$.
Set $I_k=[a_k,b]$.
We obtain a map of pairs
$$
G_k:\bigg( \pwkm, \pwkmm\bigg)\times (I_k,\pr I_k)\to (\Wmok, \Wmokm)
$$

Define a \ho~
\begin{equation}
\SS:C_{k-1}(u_1)\to H_*(\Wmok, \Wmokm)\lbl{f:s}
\end{equation}
by
the formula
$\SS(x) =(G_k)_*(x\otimes \iota)$
(where $\otimes$ is the K\"unneth product and $\iota$
is the fundamental class of $(I_k,\pr I_k)$).

\item
Let $p\in S_k(\phi)$.
Recall  that we denote by $d(p)$ the embedded disc $D(p,v)\cap
\phi^{-1}\([a_k,b]\)$.
The condition $(\gC)$ implies that
$\pr d(p)\sbs \Wmokm$.
We obtain thus a \ho~
\begin{equation}
C_k(v)=H_k\(\wk, \wkm\)\rTo^ \II H_k(\Wmok, \Wmokm)\lbl{f:i}
\end{equation}

\enen
\beth\lb{t:mtf}
\been\item
$H_*(\Wmok, \Wmokm)=0\mx{ if } *\not= k$
\item
The homomorphism
\begin{equation}
L_k=(\JJ_1, \JJ_0, \II,\SS ):
C_k(u_1)\oplus C_k(u_0)\oplus C_k(v)\oplus
C_{k-1}(u_1)\to H_k(\Wmok, \Wmokm)     \lbl{f:mtf}
\end{equation}
is an isomorphism.
 \enen
\enth
\bere\lb{r:mtf}
 This theorem implies in particular that
$H_k(\Wmok, \Wmokm)     $ is a free abelian group
of rank
$\sharp S_k(\phi_1)+
\sharp S_k(\phi_0)+
\sharp S_k(\phi)+
\sharp S_{k-1}(\phi_1)$

\enre
\Prf
We fix the value of $k$ for the proof. Set
$$
a_k=\l, a_{k+1}=\mu,
Z=T(\pwkm, v)\scc{[\l,b]}
$$
Further, set
\begin{gather}
\D=\cup_{p\in S_k(\phi)} d(p),      \lbl{f:defR}\\
R=Z\cup\D\cup \pwk \cup Y_{k+1}\scc{\l},     \lbl{f:defS}            \\
S=\Wmokm\scc{[\l, b]}
\end{gather}

so that $S\sbs R\sbs\phi^{-1}\([\l,b]\)$.
The next lemma reduces the study of the \hog~
of $(\Wmok, \Wmokm)$
to the study of the intersection of
this pair with
$\phi^{-1}\([\l,b]\)$.

\bele\label{l:2incl}
The inclusions of pairs
\begin{gather}
(R,S)\rInto (\Wmok\scc{[\l,b]} , \Wmokm\scc{[\l,b]})
\rInto (\Wmok, \Wmokm) \quad \lbl{f:1incl}   \\
(\D\cup \(Y_{k+1}\scc{[a_k]}\)  , Y_{k}\scc{[a_k]})
\rInto
  (\Wmok\scc{[\l,\mu]} , \Wmokm\scc{[\l,\mu]})\quad   \lbl{f:2incl}
\end{gather}
induce isomorphisms in \hog.
\enle
\Prf
We shall prove that (\ref{f:2incl})
induces an isomorphism in \hog~. The proof for the
(\ref{f:1incl})
is similar.
Using the $(-v)$-shift in the domain
$
\Wmok\scc{[\l,\mu]}\sm\Wmokm\scc{[\l,\m]}
$
it is easy to prove that the inclusion
$$
\(Y_k\scc{[\l,\m]}\cup Y_{k+1}\scc{\l}, Y_k\scc{\l})
\rInto
\(\Wmok\scc{[\l,\m]}, \Wmokm\scc{[\l,\m]}\)
$$
is a \hog~equivalence.
The intersection of the sets
$  Y_k\scc{[\l,\m]}   $
and
$Y_{k+1}\scc{\l}$
is in $Y_k\scc{\l}$, thus the quotient
$
 \(Y_k\scc{[\l,\m]}       \cup Y_{k+1}\scc{\l} \)
/
Y_k\scc{\l}
$
is the wedge
$$
 Y_k\scc{[\l,\m]}/Y_{k}\scc{\l}
\bigvee
 Y_{k+1}\scc{\l}
/Y_k\scc{\l}
$$
Use again the $(-v)$-shift in the domain
$\phi^{-1}([\l,\m])\sm
 Y_k\scc{[\l,\m]}       $
and then excision to see
that
$Y_k\scc{[\l,\m]}/Y_{k}\scc{\l}
\rInto
\phi^{-1}([\l,\m])/\phi^{-1}(\l)$
is a \hog~equivalence.
The quotient
$
\(\D\cup
 Y_{k+1}\scc{\l}\)
/Y_k\scc{\l}
$
is homeomorphic to the wedge
$$
(\D/\pr\D)\vee
(Y_{k+1}\scc{\l}\)
/Y_k\scc{\l})
$$
and the standard Morse-theoretic argument identifying the \hot~type of
$
(\phi^{-1}([\l,\m]), \phi^{-1}(\l))
$
with that of $(\D,\pr\D)$
finishes the proof. $\qs$

So it remains to  compute the \hot~type of $R/S$. Set
$$
\rho=R\scc{\l},\quad \s=S\scc{\l}
$$
Note that $R$ is obtained from $S$
by attaching four compact subsets:
$\phi_1^{-1}\([\a_k,\a_{k+1}]\),
G_k(\pf {k-1}{k})\times [\l,b]),
\D, \rho$.
Every two of these subsets
intersect by a subset of $S$,
therefore
$$
R/S=
\(\pwk/\pwkm\)\vee
\Sigma\(\pwkm/ \pwkmm\)
\vee
\((\D\cup\rho)/\s\)
$$
where $\Sigma$ is suspension.
Now we shall describe $(\D\cup\rho)/\s$.
Since $\D\cap\rho\sbs\s$,
we have a homeomorphism
\begin{equation}
(\D/\pr \D )\vee(\rho/\s)\rTo^{(i_1, i_2)}_{\approx}
\(\D\cup\rho\)/\s\lbl{f:bouq}
\end{equation}
The space $\D/\pr\D$ is homeomorphic to
the wedge of all the spaces
$d(p)/\pr d(p)$
and $\rho/\s$ is homeomorphic
(via the map $\stind {(-v)}{\l}a$)
to
$\( \mwk /\mwkm\)$.
Thus we obtain the following homeomorphism:
\begin{equation}
\begin{gathered}
R/S\approx
\(\pwk/\pwkm\)\vee
\Sigma\(\pwkm/\pwkmm\)\vee\\
\(\bigvee_{p\in S_k(\phi)} d(p)/\pr d(p)\)\vee
\(\mwk/\mwkm\)
\end{gathered}
   \lbl{f:rk/sk}
\end{equation}
The \hog~of this last wedge is isomorphic to
$$
C_k(u_1)\oplus C_k(u_0)\oplus C_k(v)\oplus C_{k-1}(u_1)$$
It is clear that the images of the direct summands of
this group in $H_k(\Wmok, \Wmokm)$ equal to
$\Im \JJ_1$, resp. $\Im \JJ_0, \Im\II, \Im \SS$.
$\qs$

\bere\lb{re:bypro}
Note here for the further use a byproduct of our proof:
the inclusion
\begin{equation}
(\D/\pr\D)\vee(\rho/\s)\rInto^j
\Wmok\scc{[\l,\mu]}/ \Wmokm\scc{[\l,\mu]}
\end{equation}
is a \hog~equivalence,
and therefore it induces an isomorphism
\begin{equation}
J:C_*(v)\oplus C_*(u_0)\rTo^\approx
H_*(\Wmok\scc{[\l,\mu]} , \Wmokm\scc{[\l,\mu]})  \lbl{f:bypro}
\end{equation}
\enre

\subsection{Boundary operators, associated with Morse-type filtrations}
\lb{su:boamtf}

We have seen in the previous subsection, that $\Wmok$ form
a cellular filtration of $W$. Since $\Wmok=\emp$
for $k<0$,
the \hog~of $W$ can be reconstructed from
the chain complex, associated with the filtration.  Namely, let
$E_k=H_k(\Wmok, \Wmokm)$. The Proposition 1.3 of \cite{dold}, Ch. 5
says that $H_*(W)$ is isomorphic
to the \hog~of the complex $(E_*, d_*)$, where $d_k:
E_k\to E_{k-1}$ is the boundary operator of the exact
sequence of the triple
$\(\Wmok, \Wmokm, \Wmokmm\)$.
In this subsection
we compute the matrix of this boundary operator \wrt~ the
direct sum
decomposition of $E_*$ provided by Theorem \ref{t:mtf}.
Set $d'_s=L_s^{-1}\circ D_s\circ L_{s+1}$.
Then $d'_s$ is given by its $(4\times 4)$-matrix of group
\ho s.
To describe the components of this matrix, recall that
$C_*(u_1), C_*(u_0), C_*(v)$ are themselves chain complexes,
 and denote the corresponding boundary operators
by $\pr_*^{(1)}, \pr_*^{(0)}, \pr_*$.
\bepr\lb{p:matr}
The matrix of $d'_{k+1}$ is
\begin{equation}
\left(
\begin{matrix}
\pr_{k+1}^{(1)} & 0                &         0       &\Id \\
0               & \pr_{k+1}^{(0)}  & P_k              & -\HH_k(-v) \\
0               & 0                & \pr_{k+1}         & -N_k   \\
0                &  0              &   0              & -\pr_k^{(1)}
\end{matrix}
\right)\lbl{f:matr}
\end{equation}
\enpr

Here $\HH_k(-v)$ is the operator of homological
gradient descent, (see Theorem \ref{t:consh})
and $P_k, N_k$ are
some homomorphisms, which will be defined   during the proof.

\Prf
We shall use in the proof the terminology of Subsections 
\ref{su:dmtf} and
\ref{su:hwkwkm}.

\pa
{\it  The first two columns}
\pa
They are obtained from the definition of $\JJ_0, \JJ_1$ by
the functoriality of the boundary operators associated with filtrations.

\pa
{\it  The third column}
\pa
By Lemma \ref{l:2incl}
it suffices to calculate the \hog~ class of $\pr d(p)$
in
$
\( \Wmok\scc{[\l,b]},
\Wmokm\scc{[\l,b]}
\).
$
Since $\pr d(p)$
belongs to $\phi^{-1}(\mu)$, it suffices to
calculate
the \hog~ class
of $\pr d(p)$ in
$
\( \Wmok\scc{[\l,\mu]},
\Wmokm\scc{[\l,\mu]}
\).
$
The \hog~  of this pair is
isomorphic to
$C_k(v)\oplus C_k(u_0)$
(by Remark \ref{re:bypro}).
The projection of
$[\pr d(p)]$ to the first component of this direct sum
is easily identified with
$\pr_{k+1}([p])$.
Denote by $P_k$
the projection onto the second component, and
the calculation of
the third column is finished.

\pa
{\it Homomorphism $N_k$}
\pa

Let $x\in C_k(u_1)= H_*(\pwk, \pwkm)$.
\lb{homhk}

Choose a singular chain $\bar x$, representing $x$.
The condition $(\gC)$
implies that $(-v)$
descends $\bar x$ to $\phi^{-1}(\mu)$;
consider
the descended singular chain
$\bar{\bar x}$
as a cycle in
the pair
$$
\( \Wmok\scc{[\l,\mu]},
\Wmokm\scc{[\l,\mu]}
\)
$$
We denote by $N_k(x)$
the projection of
$J^{-1}([x])$
to the first direct summand $C_*(v)$
of
(\ref{f:bypro}).

Recall that $C_k(v)$
is the free abelian group
generated by
$[p], p\in S_k(\phi)$. Write
$N_k(x)=\sum_{p\in S_k(\phi)}\langle x,p\rangle[p]$.
Then the map
$x\mapsto \langle x,p\rangle$
is a \ho~ $C_k(u_1)\to\ZZZ$.
This \ho~will be useful in the sequel, and
now we shall give an interpretation of
the number
$\langle x,p\rangle$
when the following restriction on $x$ holds:
\pa
$x$ is represented by an oriented \sma~ $X$ of $\daw$
belonging to
$\Int \pwk$
and transversal to
$D(p,-v)\cap\daw$.
\pa
Consider the cooriented $(n-k-1)$-dimensional \sma
$B_p=D(p,-v)\cap\daw$
of $\daw$.
The condition $(\gC)$ implies that
$B_p\sbs \phi_1^{-1}([\a_k, \a_n])$.
We have the following formula for the algebraic intersection
index
of $X$ and $B_p$:
\begin{equation}
\langle x,p\rangle = X\krest B_p  \lbl{f:kre1}
\end{equation}
Thus
\begin{equation}
N_k(x)=\sum_{p\in S_k(\phi)} (X\krest B_p)[p] \lbl{f:kre2}
\end{equation}

{\it The fourth column}
\pa
Let $x\in C_k(u_1)= H_k(\pwk, \pwkm)$.
For notational simplicity we shall assume that
$x$ is represented by an oriented submanifold $X$ of $\pr_1 W$,
 $\dim X=k$, \sut~
$X\sbs \pwk, \pr X\sbs \pwkm$.
Then it follows from the definition
of $\SS(x)$, see Subsection \ref{su:hwkwkm}, that
$d_{k+1}(\SS(x))$
is represented by
the following sum
of oriented singular manifolds with boundary
:
\begin{equation}
X-G_k(\pr X\times[a_{k+1}, b])- G_k(X\times\{a_{k+1}\})\lbl{f:ds}
\end{equation}
The first term corresponds to the top of the 4th column.
The second and the third terms
of
(\ref{f:ds})
do not represent  cycles of $(\Wmok, \Wmokm)$.
We can repair this as follows.
Every $(-v)$-\tr~  starting in $\pr X$ reaches  $\phi^{-1}(a_k)$,
and we can write
\begin{equation}
G_k(\pr X\times [a_{k+1},b])=
G_{k}\(\pr X\times [a_k,b])
 -
G_{k}\(\pr X\times [a_k,a_{k+1}])
\lbl{f:dsds}
\end{equation}
and then $d_{k+}(\SS(x))$ is represented by the sum
\begin{equation}
X-
G_{k}\(\pr X\times [a_k,b])
+ Y  \lbl{f:dsdsds}
\end{equation}
where
$Y=
G_{k}\(\pr X\times [a_k,a_{k+1}])
-
G_{k}\(\pr X\times \{a_{k+1}\}).
$
Now the two last  terms in the righthand side of
(\ref{f:dsdsds})
 represent cycles in
$(\Wmok, \Wmokm)$.
The
 \hog~ class of
$G_{k}\(\pr X\times [a_k,b]\) $
equals obviously to $\SS(\pr^{(1)}_kx)$.
It is not difficult to check that
$\LL_k^{-1}(Y)$
equals to the element
$-N_k(x)-\HH_k(-v)(x)$
and that finishes the computation of the fourth column.
$\qs$

\subs{Morse-type filtrations of regular coverings of
cobordisms}
\lb{su:fccrc}

In the preceding subsection we have seen that
it is possible to calculate the \hog~
$H_*(W)$
from a chain complex $E_*$, associated to the filtration $\Wmok$.
In the present subsection we refine this result and
show how to recover the simple \hot~type of $W$
from the homological
data associated to the filtration.

Let $p:\wi W\to W$ be a regular covering of $W$ with structure
 group $H$.
For $A\sbs W$ we shall denote $p^{-1}(A)$ by $\wi A$.

The Morse type filtration $\Wmok$ constructed in the previous
subsection
gives rise to the filtration
$\tiWmok=p^{-1}(\Wmok)$ of $\wi W$.
An immediate generalization of the results of the previous
 subsection leads to the following result, describing the
structure of relative
\hog~
$H_*(\tiWmok, \tiWmokm)$.

\bepr\lb{p:mtfcov}
\been
\item $H_s(\tiWmok, \tiWmokm)=0$, if $s\not= k$.
\item Set $\wi E_k= H_k(\tiWmok, \tiWmokm)$.
Then $\wi E_k$ is a free $\ZZZ H$-module
and there is an isomorphism of $\ZZZ H$-modules
\begin{equation}
\wi L_k:\wi C_k(u_1)\oplus \wi C_k(u_0)
\oplus\wi C_k(v)\oplus \wi C_{k-1}(u_1)\to \wi E_k
\lbl{f:dirscov}
\end{equation}
\item
Let $\wi d_{k+1}:\wi E_{k+1}\to\wi E_k$
be the boundary operator in the exact sequence of the triple
$(\wi W^{\langle k+1\rangle}, \tiWmok, \tiWmokm)$.
Then the matrix of the \ho~
$\wi d'_{k+1}=
(\wi L_k)^{-1}\circ
\wi d_{k+1}\circ
\wi L_{k+1}$
is

\begin{equation}
\left(
\begin{matrix}
\wi\pr_{k+1}^{(1)} & 0                &         0       &\Id \\
0               & \wi\pr_{k+1}^{(0)}  & \wi P_k        & -\wi\HH_k(-v) \\
0               & 0                & \wi\pr_{k+1}         & -\wi N_k   \\
0                &  0              &   0              & -\wi\pr_k^{(1)}
\end{matrix}
\right)\lbl{f:matrcov}
\end{equation}
\end{enumerate}

(Here $\wi\HH_k(-v):\wi C_k(u_1)\to\wi C_k(u_0)$ is the operator
of homological gradient descent
associated to the covering $\wi W\to W$, see
\su~\ref{su:ehgd}.)
\enpr
\Prf
The proof runs parallel to the proof of Proposition \ref{p:matr}.
We shall develop here only the points which are of use in 
the sequel, namely
the construction of chains representing
generators of $\wi E_k$, and the construction and properties of the \ho
~$\wi N_k$.

\pa
{\it Generators of $\wi E_k$}
\pa
To start, we choose for every point
$x\in S(\phi)\cup S(\phi_0)\cup S(\phi_1)$
a lifting $\wi x$ of $x$ to $\wi W$.
Then for every $r\in S_k(\phi)$
a lifting
$d(\wi r)$
of the disc $d(r)$ to $\wi W$
is defined.
We have obviously
$d(\wi r)\sbs \tiWmok, \pr d(\wi r)\sbs \tiWmokm$.
We define
$\wi L_k([\wi r])$ to be the image in
$H_*(\tiWmok, \tiWmokm)$
of the fundamental class of the pair
$(d(\wi r), \pr d(\wi r))$.

Similarly we define $\wi L_k$ on the elements
$[\wi q], [\wi t]$, where $q\in S_k(\phi_1), t\in S_k(\phi_0)$.

Proceeding to the image of the generators of the fourth
direct summand,
note that the lifting $\wi q$ defines also a lifting $T_k$ 
of the manifold
$G_k(d(q)\times I_k)$ to
$\wi W$.
The space $T_k$
is an oriented topological \ma~with boundary,
thus the fundamental class in
$H_k(T_k,\pr T_k)  $
is defined.
Denote its image in
$H_k(\tiWmok, \tiWmokm)$
by
$[[\wi q]]$
and
define $\wi L_k$
on the component $\wi C_{k-1}(u_1)$ of the lefthand side of
(\ref{f:dirscov})
setting $\wi L_k([\wi q])=[[\wi q]]$.

\pa
{\it Homomorphism $\wi N_k$}
\pa
Let $x\in \wi C_k(u_1)= H_k(\tipwk, \tipwkm)$.
Choose a singular chain $\bar x$ representing $x$.
The condition $(\gC)$
implies that $(-v)$ descends $\bar x$ to $\phi^{-1}\(a_{k+1}\)\wi{}$.
Consider the descended singular chain $\bar{\bar x}$ as a 
cycle in the pair
$$
\bigg(\tiWmok\scc{[\l,\m]},
\tiWmokm\scc{[\l,\m]}\bigg)$$
Similarly to Lemma \ref{l:2incl} the inclusion
$$
((\D\cup Y_{k+1}\scc{a_k} )\wi{} ,
\wi Y_{k}\scc{a_k})
\rInto
  (\tiWmok\scc{[\l,\mu]} , \tiWmokm\scc{[\l,\mu]})
$$
induces an isomorphism in \hog:
$$\wi J: \wi C_k(v)\oplus \wi C_{k-1}(u_0)\rTo
H_*(\tiWmok\scc{[\l,\mu]} , \tiWmokm\scc{[\l,\mu]})$$
Denote by $\wi N_k(x)$
the projection of
$\wi J^{-1}([\bar{ x}])$
to the first summand.
We have an equivariant analog of the formulas (\ref{f:kre1}) and
(\ref{f:kre2}).
Namely, for every $p\in S_k(\phi)$
choose a lifting $\wi p\in\wi M$.
Recall that the \hog~ classes $[\wi p]$
of the pairs $(d(\wi p), \pr d(\wi p))$
form a base of the module $\wi C_k(v)$. Write
$$
\wi N_k(x)=\sum_{p\in S_k(\phi)} [\wi p]
\langle\langle
x,p\rangle
\rangle
$$
with
$\langle\langle
x,p\rangle
\rangle\in\ZZZ H$; then the map
$x\to
\langle\langle
x,p\rangle
\rangle
$
is a \ho~of $\ZZZ H$-modules
$\wi C_k(u_1)\to\ZZZ H$.
Assume now that $x$ is represented by a lifting $\wi X$ to $\tidaw$
of an
 oriented \sma~ $X$ of $\daw$
belonging to
$\Int \pwk$
and transversal to
$B_p=D(p,-v)\cap\daw$.
We have a lifting $\wi B_p$
of the \sma~ $B_p\sbs \daw$, and for every $g\in H$
the intersection index
$\wi X\krest \wi B_p g\in\ZZZ$
is defined.
We have:
\begin{equation}
\langle\langle
x,p\rangle
\rangle=
\sum_{g\in H}
\bigg(\wi X\krest \wi B_pg\bigg) g
\lbl{f:kreeq1}
\end{equation}
and thus
\begin{equation}
\wi N_k(x)=
\sum_{p\in S_k(\phi), g\in H}
\bigg(\wi X\krest \wi B_pg\bigg) [\wi p] g\lbl{f:kreeq2}\qquad\qs
\end{equation}

Thus  the chain complex
$\{...\to\wi E_k\rTo^{\wi d_k}\wi E_{k-1}\to...\}$
is a free $\ZZZ H$-complex.
The point (1)
of the above proposition shows that the filtration
$\tiWmok$
is cellular, therefore
$H_*(\wi E_*)\approx H_*(\wi W)$
(and it is not difficult to show that this isomorphism
can be chosen as to commute with the $\ZZZ H$-action).

We can say more.
Choose a $C^1$ triangulation
$\D$ of $W$
and lift it to $\wi W$ to obtain a $H$-invariant
triangulation
of $\wi W$.
Then the chain simplicial complex $C_*^\D(\wi W)$ of $\wi W$
is a free chain complex of finitely generated
$\ZZZ H$-modules. Moreover, choose some
liftings of simplices of $W$ to $\wi W$, and
 $C_*^\D (\wi W)$ becomes a based chain complex.

On the other hand, if we choose orientations
of descending discs, and liftings to $\wi W$
of all the critical points of $\phi_0, \phi_1, \phi$,
then we obtain a $\ZZZ H$-base in each of the chain complexes
$\wi C(u_1), \wi C_*(u_0), \wi C_*(v)$
and (via the isomorphism $\wi L_k$)
a free $\ZZZ H$-base in $\wi E_*$.

To state the next proposition we need a definition.

\begin{defi}\lb{d:complex}
\been\item
Let $R$ be a ring.
An {\it $R$-complex}
is a finite chain complex of finitely
 generated free based right $R$-modules.
\item
Two $\ZZZ H$-complexes $C_*, D_*$ are said to be simply \hot~
  equivalent if there is a \hot~
equivalence
$\phi:C_*\to D_*$
\sut~ the torsion $\tau(\phi)\in Wh(\ZZZ H)$
vanishes.
\enen
\end{defi}

\bepr\lb{p:shecob}
The $\ZZZ H$-complexes
$C_*(\wi W)$ and $\wi E_*$ are simply \hot~ equivalent.
\enpr
\Prf
For the case of closed manifolds (that is $\pr W=\emp$)
this proposition is proved in \cite{patou}, Appendix, Theorem A.5.
In the general case the proof is similar. We shall just
give the outline of the proof.
Choose a $C^1$ triangulation of $W$ so that it
satisfies the following condition $(\TT)$:

\pa
$(\TT)$:\lb{condT}

All the subspaces $\dow,\daw$ are
simplicial subspaces of $W$.
All the level surfaces
$\phi^{-1}(a_i),
\phi_1^{-1}(\a_i),
\phi_0^{-1}(\b_i)$
and the discs
$d(p), d(q), d(r), \mbox{  where }
p\in S(\phi), q\in S(\phi_0), r\in S(\phi_1)$ are
simplicial subcomplexes of
$W$.
For every $k$ the subspaces $\Wmok$
and the subspaces $R,S$ defined in (\ref{f:defR})
and (\ref{f:defS})
are simplicial subcomplexes of $W$.
\pa

(The existence of such triangulation is a standard
 exercise in the theory of triangulating of
manifolds,
see \cite{patou}, p.330 -- 332 for a detailed proof of a
similar result).

In what follows we assume that the reader is familiar with
 the paper \cite{patou}.
To adjust our
present terminology  to that of  of \S 3 of \cite{patou},
abbreviate
$C_*^\D(\wi W)$ to
$C_*$,
and denote $C_*^\D(\tiWmok)$ by $C_*^{(k)}$;
then the $\ZZZ H$-complex $\wi E_*$ defined above is just
$C_*^{gr}$
in the terminology of \cite{patou}.
The proposition \ref{p:mtfcov}
implies that the filtration
$C_*^{(k)}$
is {\it nice}
(see \cite{patou}, page 310 for the definition of nice filtration).
Therefore (\cite{patou}, Corollary 3.4), there is
a canonical
chain \hot~equivalence
$\xi:C_*^{gr}\to C_*$.
The \hot~ class of this chain \hot~ equivalence is uniquely determined
by the following requirement:
\lb{chheq}
\been \item[F1]
$\xi(C_k^{gr})\sbs C_*^{(k)}$ for every $k$,
\item[F2]
The \ho~
$C_k^{gr}\to H_k(C_*^{(k)}, C_*^{(k-1)})$
induced by $\xi$ equals to
the identity.
\enen

Define a $\ZZZ G$-complex $ \mu(k)_*$
setting
\begin{equation}
\begin{cases}
\mu(k)_j=0 & \mx{ if } j\not= k\\
\mu(k)_k
=H_k(C_*^{(k)}, C_*^{(k-1)}) & {}
\end{cases}
\end{equation}

The chain map $\xi$ gives rise to chain \hot~  equivalences
\begin{equation}
\xi_k:\mu(k)_*\to C_*^{(k)}/C_*^{(k-1)}\quad\lbl{f:heqeq}
\end{equation}
One proves by an easy induction argument that if
all maps $\xi_k$
are simple \hot~ equivalences, then
$\xi$ is also a simple \hot~
equivalence.

Thus it remains to show that every $\xi_k$
is a simple \hot~equivalence.
This is done in three steps. 

First step: one shows that the inclusion
$(R,S)\rInto (\Wmok, \Wmokm)$
induces a \she~  of $\ZZZ H$-complexes
\begin{equation}
j:C_*(\wi R, \wi S)\to
C_*(\tiWmok, \tiWmokm)=
C_*^{(k)}/C_*^{(k-1)}
\end{equation}
(the argument is parallel to the proof of Lemma \ref{l:2incl}).
\pa
Second step:
The map
$\xi_k$ factors through
a \hot~equivalence
$\xi_k':\mu(k)_*\to C_*(\wi R, \wi S)$.
The \hog~of $L_*=C_*(\wi R, \wi S)$
vanishes except in the dimension $k$ and the module
$H_k(\wi R, \wi S)\approx \wi E_k$
is free.
Thus, introducing the trivial filtration
$\{0, L_*\}$
in $L_*$
we shall obtain a chain \hot~equivalence
$\xi''_k:\mu(k)_*\to L_*$
(by \cite{patou}, Corollary 3.4)
and it is not difficult to see that $\xi''_k$
is homotopic to $\xi'_k$.
\pa
Third step:
One proves that $\xi''_k$ splits as a sum of four \hot~equivalences
corresponding to the splitting of $R/S$
into the wedge of four summands of
(\ref{f:rk/sk}),
and one proves that each of these equivalences is simple.
$\qs$
\newpage

\section{Algebraic lemmas.}
\lb{s:aprel}

We shall use here the terminology from
Subsection
\ref{su:algk1}
and we introduce some more.
We shall say that $\phi$ is an {\it $U$-\she}
if  $\tau(\phi|U)$
 vanishes.
If the value of $U$ is clear from the context
we abbreviate
{\it $U$-simple \hot~equivalence}
to
{\it simple \hot~ equivalence}.

Note that if each $\phi_i:F_i\to D_i$ is an \iso,
then $\tau(\phi)=\sum_i(-1)^i[\phi_i]$
\footnote{We adopt here the sign convention of
\cite{milnWT}, so that our torsion $\tau(C_*)$
differs by sign from that of Turaev's one
(\cite{tur}).}

\subs{ A lemma on cone-like chain complexes}
\lb{su:lclcc}

Let $C_*$ be an $R$-complex, let $d_k:C_k\to C_{k-1}$
be its boundary \ho s.
Set $E_k=C_k\oplus C_{k-1}$.
Assume that for each $k$ an \iso~
$A_k:C_k\to C_k$
and a \ho~
$d'_k:C_k\to C_{k-1}$
are given. Let
$\DD_k:E_k\to E_{k-1}$
be the \ho, with the  matrix
$\begin{pmatrix}
d_k   & A_{k-1}  \\
0 & d'_{k-1}
\end{pmatrix}$

Assume that $\DD_{k-1}\circ\DD_k=0$, so that
$(E_*, \DD_*)$
is a chain complex.

\bele\lb{kontor}
\been
\item
$E_*$ is chain \hot~equivalent to 0.
\item
$\tau(E_*)=\sum_i(-1)^{i+1}[A_i]$
\enen
\enle
\Prf
Consider the cone $L_*$ of the identity
map $\Id:C_*\to C_*$;
we have:
$L_k=C_k\oplus C_{k-1}$,
and the matrix of the boundary operator
$\D_k:L_k\to L_{k-1}$
is
$
\begin{pmatrix}
d_k   & \Id  \\
0 & -d_{k-1}
\end{pmatrix}$.

We have of course $\tau(L_*)=0$.
Now the condition for $\DD_k$ to be boundary operator
is equivalent to
the conjunction of two following conditions:
\been\item
$d'_k$ is boundary operator
\item
$d_kA_k+A_{k-1}d'_k=0$ for every $k$.
\enen
Define a map $\l:E_*\to L_*$
by
$(c_k, c_{k-1})\mapsto (c_k, A_{k-1}c_{k-1})$
A simple computation shows that $\l$ is an \iso
~of complexes. Therefore $E_*$ is chain \hot~equivalent to 0, and
$\tau(L_*)=\tau(E_*)+\tau(\l)=\tau(\l)$.
But
$\tau(\l)=\sum_i(-1)^{i+1}[A_i] \qs$

\subs{The ring $A[t]$}
\label{su:rgra}

Let $A$ be an arbitrary commutative ring with a unit, let $R=A[t]$.
Consider the embedding
$j:A\rInto A[t]$
and the projection
$\pi:A[t]\rOnto A[t]/tA[t] =A$.
We have $\pi j=\id$.
If $C_*$ is an $R$-complex, then  the complex
$ C_*^-=C_*\tens{A} R$
is called
{\it $R$-extended from $C_*$}
(or just {\it extended} if \noconf)
\footnote{ The reason for the notation $C_*^-$
will be clear in the next section, see (\ref{f:minus}).}
If $f:C_*\to D_*$
is a map of $A$-complexes, then the map
$f\otimes\id: C_*\otimes R\to D_*\otimes R$
is called {\it $R$-extended from $f$}.
We shall keep for $f\otimes \id$ the same notation
$f$, if \noconf.
Note that
$C_*= C_*^-/t C_*^-$, thus
$H_*(C^-_*)=0$ \ifff~$H_*( C_*)=0$.
If  $U$ is  a subgroup of $A^\bu$, then
 the map
$j_*:K_1(A|U)\to K_1(R|U)$
is a split injection, the right inverse
being given by
$\pi_*: K_1(R|U)\to K_1(A|U)$.
The following lemma is now obvious.

\bele
The complex $ C_*^-$ is $U$-simply \hot~ equivalent to zero
\ifff~ $C_*$ is $U$-simply \hot~ to zero. $\qs$
\enle

\beco\lb{co:ext}
Let
\begin{equation}
0\rTo M_* \rTo^{\phi} N_*\rTo P_*\rTo 0 \lbl{tp}
\end{equation}
be an exact sequence of $R$-complexes \sut~$P_*$
is extended. Then $\phi$ is a $U$-\she~
\ifff~ $\phi/t$ is an $U$-\she.
\enco
\Prf
The exact sequence (\ref{tp})
splits as an exact sequence of $R$-modules, therefore
it remains exact
after
tensor multiplication by $R/tR$.
Further, the \hog~
of $P_*$ vanishes \ifff~ the \hog~ of $P_*/tP_*$
vanishes, and thus $\phi/t$
is a \hot~ equivalence \ifff~ $\phi$
is.
Moreover, the torsion of $\phi$, resp. of $\phi/t$
equals to the torsion of
$P_*$
(resp. $P_*/tP_*$). But $\tau(P_*)=j_*(\tau(P_*/tP_*))$
since $P_*$ is extended. $\qs$

\newpage

\section{Condition $(\gC\CC)$, Morse-type filtrations for
 Morse maps $M\to S^1$, and $C^0$-generic
rationality of Novikov incidence coefficients}
\lb{s:ccmtfrat}
\pa
\subs{Condition $(\gC\CC)$ and $C^0$-generic
rationality of Novikov incidence coefficients}
\lb{su:ccrat}

In this section we proceed at last to Morse maps
$M\rTo^f S^1=\RRR/\ZZZ$. We shall use here the terminology of
\ref{su:simnt} and we introduce still some more.

Assume that $f$ is primitive, that is
$f_*:H_1(M)\to H_1(S^1)=\ZZZ$ is epimorphic.
To simplify the notation we shall assume that
$1\in S^1$ is a regular value for $f$.
Denote $f^{-1}(1)$ by $V$.
Recall that $\CC:\bar M\to M$ is
the infinite cyclic covering,
associated to $f$, and $F:\bar M\to\RRR$ is a lifting of $f$.
Denote $F^{-1}(\a)$ by $V_\a$, and
$F^{-1}\([0,1]\)$
by $W$, and
$F^{-1}\(]-\infty,1]\)$
by $V^-$.
The cobordism $W$ can be thought of as the result of
cutting $M$ along $V$.
Our choice of the generator $t$ of the structure group
of $\CC$ implies that $V_\a t=V_{\a-1}$.
Denote $Wt^s$ by $W_s$; then $\bar M$
is the union
$\cup_{s\in\ZZZ} W_s$, the neighbor copies
$W_{s+1}$ and $W_s$
 intersecting by
$V_{-s}$.
For any $k\in\ZZZ$ the restriction of $\CC$ to $V_k$
is a diffeomorphism
$V_k\to V$.
Endow $M$ with an arbitrary riemannian metric and lift it to
a $t$-invariant riemannian metric on $\bar M$.
Now $W$ is a riemannian cobordism, and actually it is
a cyclic cobordism \wrt~the isometry
$t^{-1}:\dow=V_0\to\daw=V_1$.
We shall say that $v$ satisfies condition
$(\gC\CC)$
($\gC$ for {\it cellular} and $\CC$ for {\it circle})
if the $(F\vert W)$-gradient $v$ satisfies the condition $(\gC\YY)$
from \su~\ref{su:cyc}.
The set of $f$-gradients $v$ satisfying $(\gC\CC)$
will be denoted by
$\GG\gC\CC(f)$.
It follows from the theorem \ref{t:cyc} that the set
$\GCCT(f)=
\GCC(f)\cap\GT(f)$
is $C^0$-open-and-dense in $\GT(f)$.

Let $v\in\GCC(f)$.
The condition $(\gC\YY)$ provides a Morse function
$\phi_1$ on $V_1$
together with its gradient $u_1$,
and  a Morse function
$\phi_0=\phi_1\circ t^{-1}$ on $V_0$
together with its gradient $u_0=(t)_*(u_1)$.
For every $k\in\ZZZ$ we obtain also an ordered Morse function
$\phi_k=\phi_0\circ t^{k}:V_k\to\RRR$
 and a $\phi_k$-gradient
$u_k=(t^{-k})_*(u_0)$.

\beth\lb{t:cccgr}
Let $v\in \GCCT(f)$. Then every Novikov incidence coefficient
is a rational function of $t$.
\enth
\Prf
Let $r,s\in S(f), \ind r=l+1=\ind s+1$.
We can assume that
the liftings
$\bar r, \bar s$ of the critical points $r,s$ to $\bar M$
are in $W$.
Then $n_k(r,s;v)=0$
for $k< 0$.

Consider the oriented submanifold
$S_{\bar r}=D(\bar r,v)\cap V_0$
of $V_0$.
The condition $(\gC)$ implies that $S_{\bar r}\sbs (V_0)^{\{\leq l\}}$
and $S_{\bar r}\sm \Int (V_0^{\{\leq l-1\}})$ is compact,
so there is the fundamental class
$[S_{\bar r}]$
of $S_{\bar r}$
in the group
$H_l(V_0^{\{\leq l\}}, V_0^{\{\leq l-1\}})$.
Then applying Theorem \ref{t:consh} 1)
we deduce by induction that for every $k$ the submanifold
$X_k=\stind {(-v)}0{-k} ( S_{\bar r} )$ of $V_{-k}$ is in
$V_{-k}^{\{\leq l\}}$
and
$X_k\sm \Int (V_{-k}^{\{\leq l-1\}})$
is compact.
Moreover the \hog~class of $X_k$ in
$H_*(V_{-k}^{\{\leq l\}}, V_{-k}^{\{\leq l-1\}})$
equals to
$\HH_l^k(-v)\([S_{\bar r} ]\)$.

Set $Y_k=B_{\bar st^{k+1}} = D(\bar st^{k+1}, -v)\cap V_{-k}$.
Then $Y_k$ is
 a cooriented embedded $(n-l-1)$-dimensional \sma~ of $V_{-k}$
(where $n=\dim M$) and
$Y_k\sbs \phi_{-k}^{-1}\([\a_l, \a_{n+1}]\)$.
Since $v$ satisfies \TA, we have:
$X_k\pitchfork Y_k$
and
$n_{k+1}(r,s;v)=X_k\krest Y_k$.
Applying the formula
(\ref{f:kre1}) to the \cob~ $W_{k+1}$
we obtain:
\begin{equation}
n_{k+1}(r,s;v)=
\langle [X_k], \bar s t^{k+1}\rangle
\end{equation}

Denote $t^{-1}\HH_l(-v)$ by $h$; then we have
\begin{equation}
n_{k+1}(r,s;v)=
\langle
h^k([S_{\bar r}])t^k, \bar s t^{k+1}\rangle=
\langle h^k([\bar S_r]), \bar s t\rangle
\end{equation}
Thus we have the following formula for the total Novikov
incidence coefficient:
$$
n(r,s;v)= n_0(r,s;v)+t\sum_k\xi(A^k x)t^k
$$
with
$A=h, x=[S_{\bar r}], \xi=
\langle\cdot, \bar s t\rangle$.
Now return to (\ref{f:lalg})
to see that $n(r,s;v)$
is rational. $\qs$

\subsection{ Novikov incidence coefficients associated to a regular
covering}
\lb{su:nicarc}

Here we deal with a version of Novikov complex
associated to a regular covering of the
manifold.
We assume the terminology of \su~\ref{su:simnt} and \su~
\ref{su:algk1}.
Let $f:M\to S^1$ be a Morse map, non homotopic to zero.
Let $\PP:\wi M\to M$
be a regular covering with structure group $G$.
As always in this paper we assume that $G$ is abelian.
Assume that $f\circ \PP$ is homotopic to zero;
then the covering $\PP$ is factored through
$\CC$ as follows:

\begin{diagram}[LaTeXeqno]
\wi M & \rTo^{p}       &   \bar M         \\
      &  \rdTo_{\PP}    &   \dTo_\CC      \lb{f:covs} \\
       &                 &   M             \\
\end{diagram}

Therefore there is a lifting of $f:M\to S^1$ to
a Morse function $\wi F:\wi M\to\RRR$.
The  \ho~$f_*:\pi_1M\to\ZZZ$ can be factored through a \ho
~$\xi:G\to\ZZZ$.
Denote $\Ker\xi$ by $H$ and $\xi^{-1}(]-\infty, 0])$ by $G_-$.
Recall that to define the ordinary version of the Novikov complex
(\su~\ref{su:simnt})
we need to choose orientations of the descending discs of
critical points of $f$.
Choose in addition for each critical point
$x\in S(f)$
a lifting $\wi x$ of $x$ to $\wi M$.
Let $v\in \GT(f)$.
Lift $v$ to a $G$-invariant vector field on $\wi M$
(which will be denoted by the same letter $v$).
Let $x,y\in S(f), \ind x=\ind y+1$.
Denote by $\G(x,y,g)$ the set of $(-v)$-orbits, joining
$\wi x$ to $\wi yg$. Let $\g\in \G(x,y,g)$.
Then as usual a sign $\ve(\g)\in\{-1, 1\}$
is defined.
Set
$\wi n(x,y;v)=
\sum_{g\in G}\bigg(\sum_{\g\in\G(x,y,g)}\ve(\g)\cdot g\bigg)$.
Then $\wi n(x,y;v)$
is an element of $\wh{\wh\L}$, andit is not difficult to prove that
  $\wi n(x,y;v)\in\Lxi$.
We shall keep the same name "Novikov incidence coefficient" for
this element, since \noconf.
Let $\wi\CC_p(v)$ be the free right $\Lxi$-module, freely generated
by $S_p(f)$.
Define a $\Lxi$-\ho
~$\wi\DD_p:\wi\CC_p(v)\to \wi\CC_{p-1}(v)$ setting
$$
\wi\DD_p(x)=\sum_{y\in S_{p-1}(f)} y\cdot\wi n(x,y;v)
$$
It turns out that $\wi\DD_p\circ \wi\DD_{p+1}=0$;
therefore we obtain a based chain complex $\wi\CC_*(v)$, called
{\it Novikov Complex}.
Moreover, there is a chain \hot~equivalence
$\phi:\wi\CC_*(v)\to C_*^\D(\wi M)\tens{\L}\Lxi$,
\sut
~the image of $\tau(\phi)$
in $K_1(\Lxi|U_\xi)$
vanishes (\cite{patou}, Th. 2.2).

In view of Theorem \ref{t:cccgr} it is natural to suppose that if
$v\in\GCCT(f)$, then the Novikov incidence
coefficients
are not merely power series, but
"rational functions of several variables".
That is, one expects
$\wi n(x,y;v)\in\L_{(\xi)}$.
It is really the case:

\beth\lb{t:greqv}
Let $v\in\GCCT(f)$.
Then all the Novikov incidence coefficients
are in $\Lx$.
\enth
\Prf It is an equivariant version of the proof of
Theorem \ref{t:cccgr}.
We shall just indicate the main changes to make.
Let $r,s\in S(f), \ind r=\ind s+1=l+1$.
We can assume that all the liftings of the critical points of $f$
to $\wi M$ are in $\wi W$.
Thus we obtain a lifting $X'$ of the manifold $X$ to $\wi V_1$.
For $a\in S(f)$ set $\bar a=p(\wi a)$.
Choose an element $\theta\in G$, \sut~$\xi(\theta)=-1$.
Every critical point of $F:\bar M\to\RRR$
is of the form
$\bar bt^k$
with $k\in\ZZZ, b\in S(f)$.
Lift the critical point
$\bar bt^k\in \bar M$
to
the point $\wi b\theta^k\in\wi M$.
Then we obtain for every $k$ a lifting
$Y_k'$ of $Y_k$ to $\wi V_{-k}\sbs \wi M$, \sut
~$Y_k'=Y_0'\theta^k$.
Similarly to the end of \su
~\ref{su:cyc}
define an endomorphism
$\wi h_l(-v)$ of the free right $\ZZZ H$-module
$H_*(\wi V_0^{\{\leq l\}},
\wi V_0^{\{\leq l-1\}})$
by
$\wi h_l(-v)=\wi\HH_l(-v)\circ(\theta^{-1})_*$.
Since we consider only one $f$-gradient, we shall
abbreviate $\wi h_s(-v)$ to $\wi h_s$.

Set$$
\wi n_k(r,s;v)=
\sum_{g:\xi(g)=-k}\bigg( \sum_{\g\in\G(r,s;g)}\ve(\g)g\bigg)
$$
Then
$\wi n(r,s;v)=
\sum_{k\geq 0} \wi n_k(r,s;v)
.$
We have:
\begin{gather}
\wi n_{k+1}(r,s;v)=\sum_{g\in H}(\wi X_k\krest \wi B_{\wi sg\theta^{k+1}}
)g\theta^{k+1}
=\\
\sum_{g\in H}(\wi X_k\theta^{-k}\krest \wi B_sg\theta)g\theta^{k+1}=
\langle\langle
[\wi X_k\theta^{-k}], \wi s\theta\rangle\rangle\theta^{k+1}
\end{gather}
Applying the theorem
\ref{t:consheqv}
we obtain by induction that
$[\wi X_k\theta^{-k}]
=
\wi h^k([\wi X_0])$.
Thus
\begin{equation}
\wi n(r,s;v)=
\wi n_0(r,s;v)+
\sum_{k\geq 0}
\langle
\langle
\wi h^k([\wi X_0]),\wi s\theta
\rangle
\rangle
\theta^{k+1}
\lbl{f:novinc}
\end{equation}
and this last element belongs to $\Lx$
which can be proved by the suitable  generalization
of (\ref{f:lalg}) $\qs$

\bere\label{re:nics}
Here is one more formula for \nics, which will be useful in the sequel.
It follows from
(\ref{f:kreeq2}).
\begin{equation}
\sum_s\wi n_{k+1}(r,s;v)[\wi s]=
\wi N_l(\wi h^k([\wi X_0])\theta^{k+1}    \lbl{f:nics}
\end{equation}
\enre

\subsection{Morse-type filtration,
 associated with a Morse map $M\to S^1$}
\lb{su:mtfmm}

We keep here the terminology of two previous sections.
Let $v\in \GCCT(f)$.
We associate to $v$ the corresponding filtration
$\Wmok$ of $W$, and we introduce a $t$-invariant
filtration in $V^-$
setting
$\vk=
\bigcup_{s\geq 0}t^s \Wmok$.
It induces in turn a $ G_-$-invariant filtration $\tivk$
of $\wi V^-=
p^{-1}(V^-)$
(see the beginning of \su~\ref{su:nicarc}
for the definition of $p$, $G_-$ etc.)
In this subsection we shall reconstruct the \sht~
of $\wi M$ from the data associated to the filtration
$\vk$.
(Of course there is already a Morse-theoretic procedure of this kind:
just take a Morse function
$g:M\to\RRR$
and associate to $g$ and to the covering $\PP$
the corresponding Morse complex.
The procedure which we propose here is better suited to Novikov
\hog~ of $M$)
The next proposition
continues the lines of
Theorem \ref{t:mtf}, Proposition \ref{p:matr} and
Proposition \ref{p:mtfcov}.
We need some more terminology.
Set
\begin{equation}
\wi C_*^-(u_1)=\wi C_*(u_1)\tens{\ZZZ H}\ZZZ G_-, \quad
\wi C_*^-(v)=\wi C_*(v)\tens{\ZZZ H}\ZZZ G_-   \lbl{f:minus}
\end{equation}
These are $\ZZZ G_-$-complexes.
Extend the $\ZZZ H$-\ho s
$
\wi\pr_*:     \wi C_*(v)\to\wi C_{*-1}(v),
\wi\pr_*^{(1)}:    \wi C_*(u_1)\to\wi C_{*-1}(u_1)
$
by $\ZZZ G_-$-linearity
to $\ZZZ G_-$-\ho s
$
\wi C_*^-(v)\to\wi C_{*-1}^-(v)$,
resp.
$
\wi C_*^-(u_1)\to\wi C_{*-1}^-(u_1).
$
We keep for the extended \ho s the same notation
$\wi\pr_*^{(1)},
\wi\pr_*$.
Choose any $\theta\in G_{(-1)}$.
We have an identification
$\ZZZ G_-=\ZZZ H[\theta]$.
We shall identify the $\ZZZ H$-submodule
$\wi C_*(u_1)\theta$ with $\wi C_*(u_0)$.
The \ho s
$\wi \HH_k(-v): \wi C_*(u_1)  \to
\wi C_*(u_0)$
and
$\wi P_k:
\wi C_*(v)\to
\wi C_*(u_0)
$
can be considered as \ho s with values in $\wi C^-_{*}(u_1)$.
Extend them by $\ZZZ G_-$-linearity to the whole of
$\wi C_*(u_1)  , \wi C_*(v)$
(keeping the same notation for the extensions).
Set
$\wi h_k=\HH_k(-v)\theta^{-1}:
\wi C_*(u_1)
\to
\wi C_*(u_1)
$

\bepr\lb{p:mtftub}
\been
\item $H_s(\tivk, \tivkm)=0$, if $s\not= k$.
\item Set $\wi  \EE_k= H_k(\tivk, \tivkm)$.
Then $\wi\EE_k$ is a free $\ZZZ G_-$-module
and there is an isomorphism of $\ZZZ G_-$-modules
\begin{equation}
\wi \LL_k:\wi C_k^-(u_1)
\oplus\wi C_k^-(v)\oplus \wi C_{k-1}^-(u_1)\rTo^{\approx} \wi \EE_k
\lbl{f:trimem}
\end{equation}
\item
Let $\wi D_{k+1}
:\wi\EE_{k+1}\to\wi\EE_k$ be
 the boundary operator in the exact sequence of the triple
$(\tivkp, \tivk, \tivkm)$.
Then the matrix of
$\wi D_{k+1}'=
\wi \LL_k^{-1}\circ
\wi D_{k+1}\circ
\wi \LL_{k+1}$
\wrt~the decomposition of the left hand side of
(\ref{f:trimem})
is

\begin{equation}
\left(
\begin{matrix}
\wi\pr_{k+1}^{(1)} &  \wi P_{k+1}    &  \Id-\theta\wi h_k(-v) \\
0               & \wi\pr_{k+1}        & \wi N_k   \\
0                &  0                    & -\wi\pr_k^{(1)}
\end{matrix}
\right)\lbl{f:trimat}
\end{equation}
\end{enumerate}

\enpr

\Prf
We shall give  the proof for a particular case when
$\wi M=\bar M$ and thus $p$ is the identity map.
The general case  differs from this
particular one only by notational complications:
add tildas and replace $t$ by $\theta$.

Recall the filtration $\wmok$ of $W$ introduced in (\ref{f:wkdef}). Set
$V_{\prec k\succ}^- =
\cup_{s\geq 0} t^s \wmok$.
It is not difficult to prove that the inclusions
\begin{equation}
(V_{\prec k\succ}^-  , V_{\prec k-1\succ}^-)
\rInto
(V_{\prec k\succ}^-  , V_{\langle k-1\rangle}^-)
\rInto
(\vk, \vkm)          \lbl{f:triinc}
\end{equation}

are homotopy equivalences
(the argument used for the proof of the Lemma \ref{l:qu}
applies here also).
Consider the subspace
\begin{equation}
W^{[k]}=\pwk\cup Y_k\cup Z_k
=\wmok\sm
(\mwk\sm\mwkm)
\end{equation}
of $W$. We have obviously
$$
\Wmokm\sbs\wkwk\sbs\wmok
$$
therefore
there is an inclusion of pairs
\begin{equation}
(\wkwk, \Wmokm)
\rInto^{\l}
(\vkvk, \vkm)
\end{equation}

It is clear that
\begin{gather}
H_k(\wkwk, \Wmokm)\approx
C_k(u_1)\oplus C_k(v)\oplus C_{k-1}(u_1)\\
H_*(\wkwk, \Wmokm)=0 \mbox{ if } *\not= k
\end{gather}
\bele
The \ho~
\begin{equation}
H_*(\wkwk, \Wmokm)\rTo^{\l_*}H_*(\vkvk, \vkm)
                                  \end{equation}
induced by the inclusion $\l$ is injective and
$H_*(\vkvk, \vkm)$
is the $\ZZZ[t]$-module, extended from $\Im \l_*$.
\enle
\Prf
Note that
$$
\wkwk\cap \phi^{-1}\(]0, a_k[)
=
\Wmokm\cap\phi^{-1}(]0, a_k[)
$$
Thus by excision we can get away the set
$\phi^{-1}\(]0, a_k[\)$
from both
$\wkwk$ and $\Wmokm$.
Applying this procedure to
$t^sW$ for every $s$, it is easy to
see that
$H_*(\vkvk, \vkm)$
is a $\ZZZ[t]$-extended module from its $\ZZZ$-submodule
$$
H_*\bigg(\wkwk\cap
\phi^{-1}\([a_k,b]\),
\Wmokm\cap
\phi^{-1}\([a_k,b]\)
\bigg)
$$
This last module is isomorphic to
$H_*(\wkwk, \Wmokm)$
and the lemma is proven.
$\qs$

Now the expression
(\ref{f:trimat})
is obtained directly from
(\ref{f:matr}); just note that in the space
$V^-$
the subspace $\dow$ coincides with $\daw\cdot t. \qs$

Thus $\eey$
 is a free chain complex over $\ZZZ G_-$.
The filtration $\tivk$
of
$\wi V^- $
induces the corresponding filtration of the singular
chain complexes.
This filtration is nice (by points 1) and 2) of
Proposition \ref{p:mtftub}),
therefore
the \hog
~$H_*(\wi V^-)$
is isomorphic to that of $\wi\EE_*$.
We can also partially recover the simple \hot~type
of
the $\ZZZ G$-complex
$C_*(\wi M)$
from the complex $\wi\EE_*$. To explain what it means
choose any smooth triangulation
 of $M$, then $\wi M$
obtains a smooth $G$-invariant triangulation.
Choose some liftings of simplices of $M$ to $\wi M$
and obtain a base of $C_*^\D(\wi M)$.
Thus the simplicial
chain complex
$\cmm$
is a $\ZZZ G$-complex.
Note further that $\wi\EE_*$ obtains a natural base via the isomorphism
$\wi\LL_k$
from (\ref{f:trimem})
thus
$\wi\EE_*$
is a $\ZZZ G_-$-complex.

\bepr\lb{p:simhot}
The
$\ZZZ G$-complexes
$\cmm$
and
$\eey\tens{\ZZZ G_-}\ZZZ G$
 are simply \hot~ equivalent
(that is, there is a \hot~ equivalence $\phi$       between them
with torsion
$\tau(\phi)\in Wh(G)= \ove{K}_1(\ZZZ G| G)$
equal to zero).
\enpr
\Prf
Choose a triangulation of $M$ so that $V$
is a simplicial subcomplex, and so that the resulting triangulation
of $W$ satisfies the condition $(\TT)$ from the page
\pageref{condT}.

Then
$$\cmm=\cvm\tens{\ZZZ G_-} \ZZZ G$$
 Thus it suffices to prove that there is a
chain
\hot~     equivalence of $\ZZZ G_-$-complexes
\begin{equation}
\psi:\eey\rTo\cvm
\end{equation}
\sut~ that $\tau(\psi)$
vanishes in $K_1(\ZZZ G_-|H)$
(recall from
that such chain \hot~ equivalences $\psi$ are called
$H$-simple).
Note that each
$\tivk$
is a $\ZZZ G_-$-invariant simplicial subcomplex
of $\tivm$.
Thus
$ C_*^\D(\tivk)$
form
 a filtration of
$C_*^\D(\tivm)$.
The part 2 of Proposition \ref{p:mtftub}
implies that this filtration is nice.
Thus there is a canonical chain \hot~ equivalence
\begin{equation}
\psi:\eey\to\cvm
\end{equation}
uniquely determined by the conditions
F1) and F2) (page \pageref{chheq}).
Define a chain complex $\wi\EE(k)_*$
by
\begin{equation}
\wi\EE(k)_*=
\begin{cases}
0 & \mbox{ if } *\not= k \\
\wi\EE_k  & \mbox{if } *=k
\end{cases}
\end{equation}
Arguing by induction, it suffices to prove
that every chain map
\begin{equation}
\psi_k:\wi\EE(k)_*\rTo C_*^\D(\tivk / \tivkm)
\end{equation}
is an $H$-simple \hot~ equivalence.

\bele\lb{l:chinc}
The chain maps
\begin{equation}
C_*^\D
   (\wi V_{\prec k\succ}^-  ,\wi  V_{\prec k-1\succ}^-)
\rOnto
C_*^\D
(\wi V_{\prec k\succ}^-  , \wi V_{\langle k-1\rangle}^-)
\rInto
C_*^\D
(\tivk, \tivkm)          \lbl{f:chtriinc}
\end{equation}

induced by the inclusions from (\ref{f:triinc})
are $H$-simple \hot~ equivalences.
\enle
\Prf
It suffices to prove that for every $k$ the inclusion
$$
J:C_*^\D( \tivkvk)\rInto
C_*^\D(\tivk)
$$
is an $H$-simple \hot~ equivalence.
The  quotient $\ZZZ G_-$-complex
$C_*(\tivk/ \tivkvk)$
is $\ZZZ G_-$-extended. (Note that $\ZZZ G_-=(\ZZZ H)[\theta]$.)
Thus it suffices (by Corollary \ref{co:ext})
to prove that
$J/\theta$ is an $H$-simple \hot~ equivalence.
The chain map $J/\theta $
coincides with the following chain \hot~equivalence,
induced by the inclusion:
\begin{equation}
C_*^\D(\tiwmok,
(\pr_0\wi W)^{\{\leq k\}}   )
\rInto^{J'} C_*^\D( \tiWmok, (\pr_0\wi W)^{\{\leq k\}}   )
\end{equation}
The proof that $J'$ is a simple \hot~equivalence
 is an exercise in the theory of simple \hot~types and is left
to the reader
(Indication:
note that $\Wmok$
is obtained from $\wmok$
by adding a subset  homeomorphic to the cylinder
$\phi_0^{-1}([\a_k, \a_{k+1}])\times [0,1]$.)
$\qs$

The map $\psi_k$
factors obviously through a chain map
\begin{equation}
\psi_k':
\wi \EE(k)_*\to C_*^\D(\tivkvk/\tivkm)      \lbl{f:psiprim}
\end{equation}
and it suffices to prove that $\psi'_k$
is a simple \hot~ equivalence.
The map $\psi_k$
is determined up to chain \hot~by the condition that for every
$x\in \wi\EE(k)_*$
the \hog~class of
$\psi_k(x)$ is $x$ itself.
Thus we can consider that for every $p\in S_k(\phi)$
the element
$\psi_k'([\wi p])$
is equal
to the fundamental class of
$[d(\wi p), \pr d(\wi p)]$
in the pair
$$
(\wmok, \Wmokm)\sbs (\vkvk, \tivkm)
$$
So we can assume that $\psi'_k$
is an extended chain map and in order to check the vanishing of
$\tau(\psi'_k|H)$
it suffices to
check the vanishing of
$\tau(\psi'_k/\theta| H)$.

The proof of this last assertion is similar to the proof of
Proposition \ref{p:shecob} and will be omitted.
$\qs$

                                                      \pa

\subsection{Change of base}
\lb{su:chbase}

In this subsection we continue the study of the chain complex
$\wi\EE_*$
introduced in the previous subsection.
We are going to compare the three chain complexes:
$\wi\EE_*, \cmm\tens{\L}\Lx$
and
$\wi\CC_*(v)$.
We proved in \cite{patou}
that  for every $f$-gradient $v$ satisfying \TA~
the $\Lxi$-complexes
$C_*^\D(\wi M)\tens{\L}\Lxi$
and
$\wi \CC_*(v)$
are \hot~ equivalent, and one can choose
this \hot~ equivalence $\phi$ in such a way that
$\tau(\phi|U_\xi)\in K_1(\Lxi|U_\xi)$
vanishes.
(The torsion of $\phi$ in the group
$\ove{K_1}(\Lxi)$
does not vanish in general, as we shall see in the sequel.)
Now we shall refine this result assuming that
$v\in\GG\CC\gC\TT(f)$.
Using the formula (\ref{f:trimem})
we shall obtain in this and the following subsection
a chain \hot~equivalence
$$
\wi \CC_*(v)\rTo^\psi \wi\EE_*\tens{\ZZZ G_-}\Lx
$$
of $\Lx$-complexes
(recall that $\wi \CC_*(v)$
is defined over $\Lx$ by Theorem \ref{t:greqv})
\sut~its torsion $\tau(\psi|G)$
is explicitly computable
in terms of the \ho s $\wi h_k$.

Set $\wi{\wi\EE}_*=\wi\EE_*\tens{\ZZZ G_-}\Lx$.
We shall make a change of base in $\wi{\wi\EE}_*$
so that the Novikov complex $\wi\CC_*(v)$
will appear as a free subcomplex generated by some base elements.
It turns out that the quotient complex is contractible.

The elements of the new
base are introduced in the definition \ref{d:novbase}; for each
critical point $p\in S_k(\phi)$ there is
a base element $\{\wi p\}$.
The geometric sense of these elements can be described as follows:
The generator $[\wi p]$
corresponds to the part
$D(\wi  p,v)\cap(\phi\circ p)^{-1}([a_k, a_{k+1}])$
of the descending disc, and the generator $\{\wi p\}$
corresponds to the totality of the descending disc $D(\wi p,v)$.

{\bf Terminological remarks}

\been\item
In order not to overburden the notation we shall
{\it identify}
$\wi\EE_k$
with
$\wi C_k^-(u_1)
\oplus\wi C_k^-(v)\oplus \wi C_{k-1}^-(u_1)$.
The base elements in these three direct summands will be denoted
as
$[\wi p]$, resp. $[\wi r]$, resp. $[[\wi q]]$,
where
$p\in S_k(\phi_1), r\in S_k(f), q\in S_{k-1}(\phi_1)$.
So these base elements (in this order) form a free base
in $\wi\EE_k$
denoted by $\BB_k$.
The identification of the $\ZZZ G_-$-submodule
$\wi C_k^-(u_1)
\oplus  0      \oplus      0\sbs \wi\EE_k$
with the submodule
$0
\oplus 0 \oplus \wi C_{k}^-(u_1)\sbs \wi\EE_{k+1}$
will be denoted by $\s_k$.
\item
We shall keep the notation $\wi D_*$
for the boundary operator as well as for the boundary operator in
$\eti$, and in $\wi\EE_*\tens{\L}\Lxi$
\item
We shall often
suppress the indices $k+1,k$
from the notation whenever \noconf.
Thus we write for example
$\wi\pr^{(1)}$
instead of
$\wi\pr^{(1)}_{k+1}$
an
$\wi\pr^{(1)}_k$.
Similarly, we shall often suppress the symbol $v$ in the notation
since in rest of the section we consider only one
$f$-gradient.
The definition below
 is the first example of the abbreviated notation.

\item
Set
$$
\ctiu=\wi C_*(u_1)\tens{\ZZZ H}\Lx
$$
$$
\ctiv=\wi C_*(v)\tens{\ZZZ H}\Lx
$$
these modules are direct summands of $\eti$.
\item
We identify $\Lx$
with its image in $\Lxi$ under the canonical inclusion.
Similarly we identify $\eti$
with its image in
$\wh\EE_*=\wi\EE_*\tens{\L}\Lxi$
\enen
\pa

\bede\lb{d:novbase}
Let $r\in S_k(f)$.
Set
\begin{equation}
\{\wi r\}=
[\wi r]
-
\sum_{j=0}^{\infty} \sigma \wi h^j (\wi P([\wi r])) \theta^j
\end{equation}
\end{defi}

Note first of all that the element
$\sum_{j=0}^\infty \wi h^j(a)\theta^j$
belongs to
$\ctiu$ for every $a\in \wi C_*(u_1)$. Thus $\{\wi r\}$
is a well defined element in $\eti$.
Note further that the element
$\{\wi r\}-[\wi r]$
is a linear combination (with $\Lx$-coefficients)
of elements $[[\wi q]]$, therefore the elements
$[\wi p], \{\wi r\},  [[\wi q]]$
(here
$r\in S_k(\phi),
q\in S_{k-1}(\phi_1),
p\in S_k(\phi)$)
form a base $\BB_k'$ in $\etik$, and
the matrix of passage from $\BB_k$ to $\BB_k'$
is a product of elementary matrices.

Denote by $\wi\CC'_k(v)$ the free $\Lx$-submodule
of $\etik$, generated by the elements
$\{\wi p\}, p\in S_k(\phi)$.
Define a \ho~
$\d:
\wi\CC'_k(v)\to
\wi\CC'_{k-1}(v)$
setting
$$
\d_k(\{\wi r\})=
\sum_{r'\in S_{k-1}(\phi)}
\{\wi r'\}   \cdot  \wi n(r,r';v)
$$
(Thus the graded module $\wi\CC'_*(v)$, endowed with the differential
$\d_*$ is isomorphic to Novikov complex $\wi\CC_*(v)$.)

\bepr\lb{p:novdif}
The matrix of the boundary operator
$\wi D_{k+1}: \etikp\to\etik$
\wrt~ the bases $\BB_{k+1}', \BB'_k$
is:
\begin{equation}
\left(
\begin{matrix}
\wi\pr_{k+1}^{(1)} & 0           &   Id-\wi h_{k+1}\theta \\
0               & \d_{k+1}     &   \eta_{k+1} \\
0                &  0         &   \D_k
\end{matrix}
\right)\lbl{f:novdif}
\end{equation}

(Here $\eta_{k+1}$ and $\D_k$
are some \ho s, we do not compute them for the moment.)
\enpr
\Prf The first and the third columns are easy.
To proceed further, we need more definitions.
For an element
$x\in\wh\EE_*$
we shall need to consider "the part of $x$ between the levels
$a$ and $b$ of the function $\wi F$."
This and similar notions are introduced in the following definition.

\bede\lb{d:height}
\been\item
{\it A monomial } of $\Lxi$
is an element of the form
$ng, n\in\ZZZ, g\in G$.
{\it A monomial of $\wh\EE_*$}
is an element of $\wh\EE_*$
of the form $x\cdot ng$, where $n\in\ZZZ, g\in G,$
and $x$ is one of the generators of the base $\BB_*$

\item
The {\it height $ht(\l)$}
of a monomial $\l=ng$ of $\Lxi$
where $n\not=0$ is by definition
$\xi(g)$. The height $ht(l)$
of a monomial $l$ of $\wh\EE_*$
is by definition the subset of $\RRR$
defined by
\begin{equation}
ht(l)=
\begin{cases}
\xi(g)+[0,1], \mbox{ if } l=[\wi r]\cdot n g, n\not=0, g\in G \\
\xi(g)+[0,1], \mbox{ if } l=[[\wi r]]\cdot n g, n\not=0, g\in G \\
\xi(g)+1, \mbox{ if } l=[\wi p]\cdot n g, n\not=0, g\in G
\end{cases}     \end{equation}

\item
Let $A\sbs\RRR$.
For $\l\in\Lxi, \l=\sum_{g\in G} n_g g$
set
$\l_A=\sum_{n_g\not= 0, \xi(g)\in A} n_g g$

For $x\in \Lxi$ define $x_A$ as the sum of all
 the monomials $\mu$ from the decomposition of $x$
\wrt~the base $\BB_k$, satisfying
$ht(\mu)\sbs A$.

\item
For $A=\{a\}$ we shall also write $x_a$ instead of $x_A$.

\item
For $n\in\ZZZ, \l\in\wh\EE_*$ we denote by $\Sigma_n\l$
the element
$(\wi D(\l_{[n,\infty[}))_n  $

\enen

\end{defi}

\bele
\lb{l:st1}
Let $n\geq 0, r\in S_k(\phi)$.
Then
\begin{equation}
\Sigma_n(\{\wi r\})
=
\wi h^n(\wi P([\wi r])\theta^n
\end{equation}
\begin{equation}
\{\wi r\}_{[-n-1, -n]}= -\s(\Sigma_{-n}(\{\wi r\}))
\end{equation}
\enle
\Prf
Obvious $\qs$
\bele\lb{l:st2}
Let $A\in\ZZZ G, \supp A\sbs \xi^{-1}([n,\infty[)$,
where
$n\in\ZZZ$.
Let $R=\{\wi r\} A$. Let $m\leq n$.
Then
\begin{equation}
R_{[m-1,\infty[}=
R_{[m,\infty[}-
\s(\Sigma_m(R))
\end{equation}
\enle
\Prf
The assertion is reduced by  linearity to the case
$A=g$ with $g\in G$, and this case is obtained easily from Lemma
\ref{l:st1}.
$\qs$
\bele\lb{l:st3}
For every $n\geq 0$:
\begin{equation}
(\d\rrr)_{[-n-1,\infty[}
=
(\d\rrr)_{[-n,\infty[}-
\s(\Sigma_{-n}(\d\rrr))
+
\wi N\wi h^n\wi (P([\wi r]))\theta^n  \lbl{f:st3}
\end{equation}
\enle
\Prf
We have
\begin{equation}
(\d\rrr)_{[-n-1,\infty[}=
\sum_{r'\in S_{k-1}(\phi)}\rrr\cdot(\wi n(r,r'))_{[-n,\infty[}
+
\sum_{r'\in S_{k-1}(\phi)}\rrr\cdot(\wi n(r,r'))_{-n-1}
\lbl{f:dvachl}
\end{equation}
By Lemma \ref{l:st2}
the first term of the righthand side (\ref{f:dvachl})
equals to the sum of the two first terms
in the righthandside
of
(\ref{f:st3}).
The second term
equals to
$\wi N\wi h^n (\wi P([\wi r]))\theta^n$
by
(\ref{f:nics}). $\qs$

\bele\lb{l:dd}
Let $n\geq 0$. Then
\begin{equation}
\wi D(\{\wi r\}_{[-n,\infty[})
=
(\delta\rrr)_{[-n,\infty[}+
\wi h^n(\wi P([\wi r]))\theta^n  \lbl{f:dd}
\end{equation}
\enle
\Prf Induction in $n$. The case $n=0$ is obvious.
Assume that the formula (\ref{f:dd}) holds for some integer $n$.
Then we have in particular
\begin{equation}
\wi D(\d\rrr_{[-n,\infty]})=
-\wi\pr^{(1)}\wi h^n(\wi P([\wi r]))\theta^n
=
\Sigma_{-n}(\d\rrr) \lbl{f:ddd}
\end{equation}
We have:
\begin{equation}
\widetilde{D}(\{\wi r\}_{[-n-1,-n]})=
\bigg(
-\wi h^n\wi P([\wi r])
+
\wi h^{n+1}\wi P([\wi r])\theta
+
\wi N\wi h^n
\wi P([\wi r])+
\s\wi\pr^{(1)}\wi h^n
\wi P([\wi r])
\bigg)
\theta^n
\end{equation}

and (by \ref{l:st3}):
\begin{equation}
(\d\rrr)_{[-n-1,\infty[}-
(\d\rrr)_{[-n,\infty[}=
\wi N \wi h^n
\wi
P([\wi r])
\theta^{n}+
\s \wi\pr^{(1)} \wi h^n
\wi
P([\wi r])
\theta^n \lbl{f:pochti2}
\end{equation}

Now just combine the above formulas to get
the assertion of the Lemma at the rank $n+1$.
$\qs$

Now the computation of the second column of the matrix
of $\wi D_{k+1}$
is finished, since Lemma \ref{l:dd}
implies
$\wi D\{\wi r\}=\d\{\wi r \}$. $\qs$

\subsection{Homotopy equivalence
$\wi\CC_*(v)   \to         \eti     $
and its torsion}
\lb{su:torexpl}

Thus we have an inclusion
$
I:
\wi\CC_*(v)\rInto \eti$
and the quotient is a free $\Lx$-complex
$\QQ_*$
where
$\QQ_k=\ctiuk\oplus \ctiukm$.
The matrix of the boundary operator
$\nabla_k:\QQ_k\to\QQ_{k-1}$
is
\begin{equation}\left(
\begin{matrix}
\wi\pr_{k}^{(1)} &  1-\wi h_k\theta \\
0               & \D_{k-1}
\end{matrix}\right)
\end{equation}
Note that the \ho~
$1-\wi h_k\theta:
\ctiuk
\to
\ctiuk
$
is an isomorphism
since $\det(1-\wi h_k\theta)$
is of the form
$1+\theta\xi$
with $\xi\in\ZZZ H$,
and such elements are invertible in $\Lx$.
Thus we can apply
Lemma \ref{kontor}
to deduce
that
$I$ is a chain \hot~equivalence
and that its torsion
$\tau(I)\in\ove{K_1}(\Lx)$
equals to
$\sum_i(-1)^{i+1}\a_i$
where
$\a_i$ is the class of the isomorphism
$(1-\wi h_i\theta)$.
Due to the special form of the \ho~
$(1-\wi h_k\theta)$
we can simplify still more.
In the next lemma
$A$ is a commutative ring with unit,
$R=A[t],
S=\{P\in R\mid P=1+tQ(t)\}, \wi R= S^{-1}R$,
 $L$
is a free \fg~ $R$-module, and
$\wi L$ stands for
$S^{-1}L$.

\bele\lb{l:dett}
Let $\xi:L\to L$
be a \ho.
Then the class of the isomorphism
$1+\xi\theta:\wi L\to \wi L$
in the group
$\ove{K}_1(\wi R)$
equals to the unit
$\det(1+\xi\theta)$
of the ring $\wi R$.
\enle
\Prf
Let $n=\rk L$.
Consider the class of invertible $(n\times n)$-matrices
$M$ over $\wi R$,
satisfying the property that
$[M]=[\det M]$
in
$\ove{K_1}(\wi R)$.
The upper triangular matrices are obviously in this class.
Note also that this class is closed under the elementary operations.
Thus it suffices to prove that our matrix
$(1+\xi\theta)$
can be reduced to an upper triangular matrix in
$\wi R$
by elementary operations. This is easy and will be left to the reader.
$\qs$

Thus we have constructed an inclusion $I:\wi C_*(v)\rInto\eti$
\sut~it is a chain \hot~equivalence
and its torsion
$\tau(I)\in \ove{K}_1(\Lx)$
satisfies
\begin{equation}
\tau(I)
=
\bigg[\prod_k \bigg(\det(1-\wi h_k\theta)\bigg)^{(-1)^{k+1}}\bigg]
\lbl{f:last1}
\end{equation}

Note that we have chosen the bases in $\wi \CC_*(v)$ and in $\eti$
in a special way:
we constructed them
from one and the same family of liftings of the critical points
of $f$ to $\wi M$.
Moreover, we required that the liftings of the critical points
 of $f$ to $\wi M$
belong to $\wi W$.

If we impose no restrictions on choosing the liftings of the
critical points to $\wi M$, we can
keep track only
of the torsions in $\ove{K_1}(\Lx|G)$
and we arrive (with the help of Proposition \ref{p:simhot})
at the following corollary.

\beco\lb{co:last}
There is a chain \hot~ equivalence
$
\phi:\wi\CC_*(v)\to\CC_*^\D(\wi M)\tens{\L}\Lx$
with
\begin{equation}
\tau(\phi|G)=
\bigg[\prod_k  \det(1-\wi h_k\theta)^{(-1)^{k+1}  }
\bigg]    \lbl{f:last2}
\end{equation}
$\qs$
\end{coro}
\section{Proof of the main theorem}
\label{s:proof}

We keep here  the terminology of the previous section.
The set of all Kupka-Smale $f$-gradients $v$, satisfying $(\gC\YY)$
will be denoted by $\GKSC(f)$,
and we shall show that this set satisfies the conclusions of the
theorem.
Indeed, the first item follows since $(\gC\YY)$
is an $C^0$-open-and-dense
condition, and since $\GKSC(f)$
is dense in $\GG\KK\SS(f)$.
The first part of the item  2)
follows from Theorem \ref{t:greqv}.
It remains to prove that
$\zeta_L(-v)      \in \Lx$ and that
the image of
$\zeta_L(-v)      $
in
$\ove{K}_1(\Lx|G)$
equals to the element

$$
\bigg[\prod_m    \det(1-\wi h_m\theta)^{(-1)^{m+1}      }
\bigg]
$$

To do this we need one more description of the operator $\wi h_m$.

\subsection{Homological gradient descent (third version)}
\label{su:hgd3}
Recall that we have a Morse function
$\phi_1:V_1\to\RRR$
and its $\d$-separated gradient
$u_1$.
For $p\in S_m(\phi_1)$ set
$$
R(p)=
\(D_\d(p,u_1)\cup V_1^{\{\leq m-1\}} \) / V_1^{\{\leq m-1\}}
$$
This space is \hot~equivalent to the sphere $S^m$.
For $k\in\ZZZ$
set
$$
Q_k^{\lceil m\rceil}
=
\wi V_k^{\lceil m\rceil}(\d)/\wi V_k^{\{\leq m-1\}}
$$
(where as usual $\wi A$ stands for
$\PP^{-1}(A)$).
The group $H$ acts on
$Q_k^{\lceil m\rceil}$
leaving invariant the point
$\a_k=[\wi V_k^{\{\leq m-1\}}]$.
Denote by $R(\wi p)$ the lifting of $R(p)$
to
$
Q_0^{\lceil m\rceil}
$
corresponding to the lifting
$\wi p\theta$
of $p$.
Then
$$
Q_k^{\lceil m\rceil}
=
\bigvee_{\substack{p\in S_m(\phi_1)\\
         h\in H, \xi(h)=k}}
 R(\wi p)\cdot h
$$

The generator of the $m$th \hog~group of $R(\wi p)$
corresponding to the chosen orientation of
$D(p,v)$
will be denoted by $\rho(\wi p)$.
The $\ZZZ H$-module
$H_*(
Q_k^{\lceil m\rceil})
$
is then free with the base
$\{\rho(\wi p)\}_{p\in S_m(\phi_1)}$.
It is not difficult to show that the gradient descent map
$\stexp {(-v)}$
defines a continuous map
$v\da
:
Q_k^{\lceil m\rceil}
\to
Q_{k-1}^{\lceil m\rceil}
$
which is $H$-equivariant, sends $\a_k$ to $\a_{k-1}$,
and \sut~the following diagram
is commutative

\begin{diagram}[LaTeXeqno]
H_*(\wi V^{\{ \leq m \}}_k  ,
\wi V^{\{ \leq m-1   \}}_k )
&
\rTo^{\wi\HH(-v)}
&
H_*(\wi V^{\{\leq  m \}}_{k-1}  ,
\wi V^{\{ \leq m-1   \}}_{k-1} )\\
\uTo_{\approx}
&
&
\uTo_{\approx}          \lbl{f:hcont}  \\
H_*(
Q_k  ^{\lc m\rc})
&
\rTo^{(v\da)_*}
&
H_*(
Q_{k-1} ^{\lc m\rc})
\end{diagram}

 We leave to the reader the proof of the
existence of $v\da$
as an exercice in homological gradient descent theory.

\subsection{Computation of the $\zeta$-function}
\label{su:compu}

Let
$Cl^{[s]}(-v)$
be the subset of all the closed orbits of $(-v)$
passing through a point of
$V^{\lceil\leq s\rceil}(\d)\sm V^{\{\leq s-1\}}$.
It follows from the property $(\gC\YY )$
that the set $Cl(-v)$
is the disjoint union of its subsets
$Cl^{[s]}(-v)$.
Set
 \begin{equation}
\eta_s(-v)=
\sum_{\g\in Cl^{[s]}(-v)} \ve(\g) \frac {\pi([ \g])}{m(\g)}, \lbl{f:etas}
 \end{equation}

 Since
 $$
 \ln\det (1-\theta\wi h_m)=
 -\sum_{k>0}
 \frac {\Tr (\theta\wi h_m)^k}k
 $$
it suffices to prove that for every $s$ we have:

 \begin{equation}
 (-1)^s\sum_{k>0}\frac {\Tr (\theta\wi h_s)^k}k
=
\eta_s(-v)\lbl{f:etaeta}
 \end{equation}

So we fix $s$ up to the end of the proof.
We shall abbreviate
$Cl^{[s]}(-v)$
to
$Cl$.

The equality (\ref{f:etaeta}) will be proved by
 translating our data to the language of
 the fixed point theory.

We shall say that a point $a\in\beta\sm\{\o\}$
is a {\it $G$-fixed-point} of $(v\da)^k$, if
$(v\da)^k(a)=a\cdot g, g\in G$.
The element $g\in G_{(-k)}$
is uniquely determined by $a$ and will be denoted by $g(a)$.
The set of all
$G$-fixed points  of $(v\da)^k$
will be denoted by $GF(k)$.
The set of all
$G$-fixed points  of $(v\da)^k$
with given $g(a)=g$
will be denoted by $GF(k,g)$.
Thus
$GF(k)=
\sqcup_{g\in G_{(-k)}}GF(k,g)$.
By analogy with the standard fixed point theory we  define the
multiplicity
$\mu(a)$ and the index $\ind a$ for every $a\in GF(k)$.

Let $a\in GF(k)$. Let $a_i=(v\da)^i(a)$
and let $\wh a_i$ be (the
unique) point in $\b$ belonging to the $G$-orbit of $a_i$.
The set of all $\wh a_i, i\in\NNN$
will be called {\it quasiorbit of $a$}, and denoted by
$\QQ (a)$; it is a finite set of cardinality $\frac k{\mu(a)}$.

Now let
$a\in GF(k)$.
Consider the integral curve
$\g$ of $(-v)$ in $\wi M$, such that $\g(0)=a$; then for some
$T>0$ we have
$$
\g(T)=
(v\da)^k(a)=a\cdot g(a)\in \wi V_{-k}
$$
The map
$\PP\circ\g: [0,T]\to M$
 is then a closed orbit. Thus we obtain a map
$\G:GF(k)\to Cl$
whose image is exactly the subset
$Cl_k\sbs Cl$, consisting of all
$\g\in Cl$
with
$f_*([\g])=-k$.

For every $\g\in Cl_k$
the set $\G^{-1}(\g)$
is the quasiorbit
of some $a\in GF_k$.
(Note that the set
$\PP(\G^{-1}(\g))$
is the intersection of the orbit $\g$ with $V$.)
Further, for every
$a\in \G^{-1}(\g)$
we have:
$g(a)=\pi([\g]),
\ve(\g)=\ind a$
and
$\m(\g)=\m(a)$.

 Thus the following power series
 \begin{equation}
 \nu(-v) =
 \sum_{k=1}^\infty \frac 1k \sum_{a\in GF_k}
(\ind a) g(a)
 \label{f:nunu}
 \end{equation}

equals to
$\eta_s(-v)$.

To prove our theorem it remains to show that
$\nu(-v)=
 (-1)^s\sum_{k>0}\frac {\Tr (\theta\wi h_s)^k}k
$. This follows obviously form the next lemma.

\bele
\lb{l:kappa}
\begin{equation}
\Tr (\theta\wi h_s)^k
=(-1)^s\sum_{a\in GF(k)}  (\ind a) g(a)
\end{equation}
\enle
\Prf
Set
$\beta=Q_0^{\lceil s \rceil}$.
The set $GF(k,g)$
is exactly the fixed point set of the following composition:
\begin{equation}
\beta\rTo^{(v\da)^k} B_{(-k)}= \bigvee_{h\in G_{(-k)}}\beta\cdot h
\rTo^{\pi_g}\beta
\end{equation}
where $\pi_g$ is the map which sends every component of the wedge
to $\o$, except the component $\beta g$, and this one is sent to
$\beta$ via the map $x\mapsto x g^{-1}$.

The point $\o\in\b$
is a fixed point of this map, and its index equals to 1.
Apply now
the Lefschetz-Dold fixed point formula
and the proof of the Lemma and of the Main theorem
is over.
$\qs$

\newpage
\label{refer}

\end{document}